\documentclass[12pt,leqno]{article}
\usepackage{amssymb}
\usepackage[mathscr]{eucal}
\usepackage{amsmath,amssymb,latexsym,theorem,bbm}
\usepackage{color,url}

\newcommand{\SC}{\scriptstyle}

\newcommand{\CC}{\mathsf{C}}
\newcommand{\DD}{\mathsf{D}}

\newcommand{\NN}{\mathbb{N}}

\newcommand{\RR}{\mathbb{R}}

\newcommand{\ZZ}{\mathbb{Z}}

\newcommand{\ba}{{\boldsymbol{a}}}
\newcommand{\bA}{{\boldsymbol{A}}}

\newcommand{\bb}{{\boldsymbol{b}}}
\newcommand{\bB}{{\boldsymbol{B}}}

\newcommand{\bc}{{\boldsymbol{c}}}
\newcommand{\bC}{{\boldsymbol{C}}}

\newcommand{\be}{{\boldsymbol{e}}}
\newcommand{\bE}{{\boldsymbol{E}}}
\newcommand{\bI}{{\boldsymbol{I}}}

\newcommand{\bm}{{\boldsymbol{m}}}

\newcommand{\bM}{{\boldsymbol{M}}}
\newcommand{\bp}{{\boldsymbol{p}}}
\newcommand{\bP}{{\boldsymbol{P}}}

\newcommand{\bQ}{{\boldsymbol{Q}}}

\newcommand{\bR}{{\boldsymbol{R}}}
\newcommand{\bcR}{{\boldsymbol{\cR}}}
\newcommand{\bS}{{\boldsymbol{S}}}

\newcommand{\bT}{{\boldsymbol{T}}}
\newcommand{\bu}{{\boldsymbol{u}}}

\newcommand{\bv}{{\boldsymbol{v}}}
\newcommand{\bV}{{\boldsymbol{V}}}

\newcommand{\bx}{{\boldsymbol{x}}}
\newcommand{\bX}{{\boldsymbol{X}}}

\newcommand{\bz}{{\boldsymbol{z}}}

\newcommand{\bU}{{\boldsymbol{U}}}

\newcommand{\bgamma}{{\boldsymbol{\gamma}}}
\newcommand{\bxi}{{\boldsymbol{\xi}}}

\newcommand{\bbeta}{{\boldsymbol{\beta}}}
\newcommand{\bvare}{{\boldsymbol{\vare}}}

\newcommand{\bzero}{{\boldsymbol{0}}}

\newcommand{\cA}{{\mathcal A}}
\newcommand{\cB}{{\mathcal B}}

\newcommand{\cD}{{\mathcal D}}

\newcommand{\cF}{{\mathcal F}}

\newcommand{\cM}{{\mathcal M}}
\newcommand{\bcM}{\boldsymbol{\cM}}

\newcommand{\cP}{{\mathcal P}}

\newcommand{\cR}{{\mathcal R}}
\newcommand{\cS}{{\mathcal S}}
\newcommand{\cU}{{\mathcal U}}

\newcommand{\bcU}{\boldsymbol{\cU}}

\newcommand{\cX}{{\mathcal X}}
\newcommand{\bcX}{\boldsymbol{\cX}}

\newcommand{\cY}{{\mathcal Y}}
\newcommand{\cZ}{{\mathcal Z}}
\newcommand{\cW}{{\mathcal W}}
\newcommand{\bcW}{\boldsymbol{\cW}}

\newcommand{\bcY}{\boldsymbol{\cY}}

\newcommand{\bcZ}{\boldsymbol{\cZ}}

\newcommand{\dd}{\mathrm{d}}

\newcommand{\slu}{{\SC\mathrm{lu}}}

\newcommand{\INARp}{\textup{INAR($p$)}}

\DeclareMathOperator*{\esssup}{ess\,sup}

\newcommand{\EE}{\operatorname{\mathbb{E}}}
\newcommand{\PP}{\operatorname{\mathbb{P}}}
\newcommand{\OO}{\operatorname{O}}

\newcommand{\var}{\operatorname{Var}}
\newcommand{\cov}{\operatorname{Cov}}

\newcommand{\vare}{\varepsilon}

\renewcommand{\mid}{\,|\,}
\newcommand{\bmid}{\,\big|\,}

\renewcommand{\leq}{\leqslant}
\renewcommand{\geq}{\geqslant}

\newcommand{\skor}[1]{\stackrel{\cS_{#1}}{\longrightarrow}}
\newcommand{\stoch}{\stackrel{\PP}{\longrightarrow}}
\newcommand{\distr}{\stackrel{\cD}{\longrightarrow}}
\newcommand{\distrf}{\stackrel{\cD_f}{\longrightarrow}}
\newcommand{\distre}{\stackrel{\cD}{=}}
\newcommand{\mean}{\stackrel{L_1}{\longrightarrow}}

\newcommand{\lu}{\stackrel{\slu}{\longrightarrow}}
\newcommand{\as}{\stackrel{{\mathrm{a.s.}}}{\longrightarrow}}
\newcommand{\ase}{\stackrel{{\mathrm{a.s.}}}{=}}

\newcommand{\bbone}{\mathbbm{1}}
\newcommand{\ns}{{\lfloor ns\rfloor}}
\newcommand{\nt}{{\lfloor nt\rfloor}}
\newcommand{\nT}{{\lfloor nT\rfloor}}
\newcommand{\nr}{{\lfloor nr\rfloor}}
\newcommand{\proofend}{\hfill\mbox{$\Box$}}

\newcommand{\bh}{{\boldsymbol{h}}}
\newcommand{\bw}{{\boldsymbol{w}}}
\newcommand{\bN}{{\boldsymbol{N}}}

\numberwithin{equation}{section}

\theoremstyle{change} \theorembodyfont{\em}
\newtheorem{Lem}{Lemma.}[section]
\newtheorem{Thm}[Lem]{Theorem.}
\newtheorem{Pro}[Lem]{Proposition.}
\newtheorem{Cor}[Lem]{Corollary.}

\theorembodyfont{\rm}
\newtheorem{Rem}[Lem]{Remark.}

\begin{document}

\begin{center}
{\bfseries\Large Asymptotic behavior of some strongly critical\\[-2mm] decomposable 3-type Galton--Watson processes\\[1mm]
 with immigration} %having triangular offspring mean matrix}% with diagonal entries 1}
\vspace*{3mm}

 {\sc\large
  M\'aty\'as $\text{Barczy}^{*,\diamond}$ \text{and}
  \ D\'aniel $\text{Bezd\'any}^{**}$ }

\end{center}

\vskip0.2cm

\noindent
 * HUN-REN--SZTE Analysis and Applications Research Group,
   Bolyai Institute, University of Szeged,
   Aradi v\'ertan\'uk tere 1, H--6720 Szeged, Hungary.

\noindent
 ** Bolyai Institute, University of Szeged,
    Aradi v\'ertan\'uk tere 1, H--6720 Szeged, Hungary.

\noindent E-mails: barczy@math.u-szeged.hu (M. Barczy),
                   bezdany@server.math.u-szeged.hu (D. Bezd\'any).

\noindent $\diamond$ Corresponding author.

\renewcommand{\thefootnote}{}
\footnote{\textit{2020 Mathematics Subject Classifications\/}: 60J80, 60F17. }
\footnote{\textit{Key words and phrases\/}:
  multi-type Galton-Watson process with immigration, critical, decomposable, asymptotic behavior, squared Bessel process.}
\footnote{M\'aty\'as Barczy was supported by the project TKP2021-NVA-09.
Project no.\ TKP2021-NVA-09 has been implemented with the support
 provided by the Ministry of Culture and Innovation of Hungary from the National Research, Development and Innovation Fund,
 financed under the TKP2021-NVA funding scheme.
D\'aniel Bezd\'any was supported by the 2025-2026 grant of the M\'oricz Doktorandusz Alap\'itv\'any.}

\vspace*{0.1cm}

\begin{abstract}
We study the asymptotic behavior of a critical decomposable 3-type Galton-Watson process
with immigration when its offspring mean matrix is triangular with diagonal entries 1.
It is proved that, under second or fourth order moment assumptions on the offspring and immigration distributions,
a sequence of appropriately scaled random step processes formed from such a Galton-Watson
process converges weakly. The limit process can be described using independent squared Bessel processes $(\cX_{t,1})_{t\geq0}$, $(\cX_{t,2})_{t\geq0}$, and $(\cX_{t,3})_{t\geq0}$,
the linear combinations of the integral processes of $(\cX_{t,1})_{t\geq0}$ and $(\cX_{t,2})_{t\geq0}$, and possibly the 2-fold iterated integral process of $(\cX_{t,1})_{t\geq0}$.
The presence of the 2-fold iterated integral process in the limit distribution is a new phenomenon
in the description of asymptotic behavior of critical multi-type Galton-Watson processes with immigration.
Our results complete and extend some results of Foster and Ney (1978) for some strongly
critical decomposable 3-type Galton-Watson processes with immigration.
\end{abstract}

\tableofcontents

\section{Introduction}\label{section_intro}

Investigation of the asymptotic behavior of branching processes has a long tradition and it is an active area of research.
In our paper, we focus on critical decomposable 3-type Galton-Watson processes with immigration (GWI processes).
For an overview of the history of results on the asymptotic behavior of critical GWI processes,
 see the introduction of Barczy et al.\ \cite{BarBezPap2}.
In what follows, we only mention and summarize those earlier results that are directly connected to our present paper.

Foster and Ney \cite[Theorems 4 and 5]{FosNey} proved limit theorems for some special strongly critical decomposable
multi-type GWI processes.
In case of a $3$-type GWI process $(\bX_k)_{k\geq 0} = ((X_{k,1},X_{k,2},X_{k,3}))_{k\geq 0}$
with a lower triangular offspring mean matrix $(a_{i,j})_{i,j=1}^3\in[0,\infty)^{3 \times3}$,
where $a_{i,i}=1$ for each $i\in\{1,2,3\}$ and $a_{i+1,i}>0$ for each $i\in\{1,2\}$,
Foster and Ney \cite[Theorem 4]{FosNey} showed that
$(n^{-1} X_{n,1}, n^{-2} X_{n,2}, n^{-3}X_{n,3})$ converges in distribution as $n\to\infty$,
 and they also characterized the limit distribution by its Laplace transform, which contains
 an integral of some function of a solution of a differential equation (see \cite[formula (9.4)]{FosNey}).
In our Theorem \ref{main_4}, we will extend this result of Foster and Ney \cite{FosNey} by proving weak convergence of
a sequence of appropriately scaled random step processes formed
from $(\bX_k)_{k\geq 0}$, and we characterize the limit process as a pathwise unique strong
 solution of a system of stochastic differential equations (SDEs).
We extend the results of Foster and Ney \cite[Section 9]{FosNey}
 in another respect as well, namely, using their notations,
 in the case of $p=N=3$, we allow that the off-diagonal entries $a_{2,1}$ and $a_{3,2}$
 of the offspring mean matrix $(a_{i,j})_{i,j=1}^3$ be $0$ as well, and in all the possible cases, we prove not only weak
 convergence of one-dimensional distributions, but weak convergence of a sequence of appropriate random step processes
 formed from the branching processes in question as well, see Theorems \ref{main_1}, \ref{main_2} and \ref{main_3}.

In Barczy et al.\ \cite{BarBezPap2}, we described the asymptotic behavior of a critical decomposable 2-type GWI process.
Under second or fourth order moment assumptions on the offspring and immigration distributions, we proved that
 a sequence of appropriately scaled random step processes formed from a critical decomposable
 2-type GWI process converges weakly.
The limit process can be described using one or two independent squared Bessel processes and possibly the unique stationary
 distribution of an appropriate single-type subcritical GWI process.
Our results completed and extended the results of Foster and Ney \cite{FosNey} for some strongly critical decomposable 2-type GWI process.
At the end of the introduction, we highlight the novelties of our present paper compared to Barczy et al.\ \cite{BarBezPap2}.

Decomposable Galton--Watson processes, with or without immigration, arise naturally as stochastic models for structured populations.
Consider a geographically structured population, where each individual is located in one of two distinct regions (such as two islands), 
 and define the type of an individual as its location.
Suppose that newborn individuals of type $1$ either stay at their birth region, or migrate to the other region,
all individuals born in the second region stay there (they do not migrate), and at each step, immigration from outside the population may occur to the regions.
This population may be modeled by a decomposable $2$-type GWI process, provided that 
 the offspring and immigration distributions depend on the regions on which the individuals are located, 
 and the immigrants join, respectively.
Jagers \cite{Jag} also pointed out that the reproduction of biological populations consisting of two types of
individuals often displays an irreversibility property described above in the sense that individuals
of one type might give birth to descendants of both kinds, whereas those of the
other type can have descendants only of their own kind. For example, if human diploid cells in
a tumour are considered the first type in the cell population, and cells of higher diploidity are
considered the second type, then, provided that endomitosis (a process where chromosomes
duplicate but the cell does not subsequently divide, causing higher ploidity) is possible, the
population of cells in this tumour has the irreversibility property in question.

Decomposable multi-type Galton--Watson processes, with or without immigration, have found other applications as well.
For example, Balelli et al.\ \cite{BalMilWai} introduced some (not necessarily critical) decomposable multi-type Galton--Watson processes to analyze 
 the interactions between division, mutation and selection in a simplified evolutionary model, 
 assuming that the population observed can be classified into fitness levels.
In their Propositions 2 and 3, the offspring mean matrix has two block matrices on its diagonal,
one of them is the $2\times 2$ identity matrix and the other one is the sum of a constant multiple of
an identity matrix and a constant multiple of the transition probability matrix
of an appropriate Markov chain. 
Recently, Sagitov and St{\aa}hlberg \cite{SagStaahl} have used
 a supercritical decomposable $4$-type GWI process (with Bernoulli immigration) to model the outcomes of
 polymerase chain reaction (PCR) bar encodings. 
For
%the
 applications of decomposable multi-type continuous time branching processes, see, e.g., Coldman and Goldie \cite{ColGol}
 and Lee and Yang \cite{LeeYan}.
 
The paper is organized as follows.
In Section \ref{Section_multi-type_GWI} we recall the notion of multi-type GWI processes, some of their basic properties (e.g., the mean function), their classification, and the notion of decomposable GWI processes.
Section \ref{Section_conv_results} contains our main results on the asymptotic behavior of critical decomposable 3-type GWI processes,
 see Theorem \ref{main_1} and Theorems \ref{main_2}--\ref{main_4}.
The investigation of such processes can be reduced to four cases
 presented in \eqref{tablazat_esetek} according to the form of the offspring mean matrix.
Under second or fourth order moment assumptions on the offspring and immigration distributions, in the above mentioned four cases,
 we describe the limit behavior of a sequence of appropriately scaled random step processes
 formed from a critical decomposable 3-type GWI process.
The limit process can be described using independent squared Bessel processes $(\cX_{t,1})_{t\geq0}$, $(\cX_{t,2})_{t\geq0}$, and $(\cX_{t,3})_{t\geq0}$,
the linear combinations of the integral processes of $(\cX_{t,1})_{t\geq0}$ and $(\cX_{t,2})_{t\geq0}$, and possibly the 2-fold iterated integral process of $(\cX_{t,1})_{t\geq0}$.
The presence of the 2-fold iterated integral process in the limit distribution is a new phenomenon
in the description of asymptotic behavior of critical multi-type GWI processes.
Section \ref{Section_prel_proof} contains preliminaries for the proofs of the main results,
and, in particular, in Proposition \ref{Pro_Xk_exp},
we describe the asymptotic behavior of the expected value of the Galton--Watson process in question.
Sections \ref{Proof1}--\ref{Proof4} are devoted to the proofs, which are based on limit theorems for sequence of martingale
 differences, delicate applications of the continuous mapping theorem and fine moment estimations for the multi-type GWI process in question.
The overall structure of the proofs of Theorem \ref{main_1} and Theorems \ref{main_2}--\ref{main_4} is the same.
Roughly speaking, first we show the joint convergence of the coordinate processes
 which correspond to those types that can not be produced by any other type.
Then, using a decomposition of the process and a version of the continuous mapping theorem,
 we describe the asymptotic behavior of the remaining coordinate processes.
We close the paper with four appendices.
Appendix \ref{appInt} contains some formulae for sums of some weighted values of a function defined on the set of nonnegative integers,
 which are used throughout the proofs.
In Appendix \ref{GWI_moments} we present some formulae and estimates for the first, second and fourth order moments of
 the coordinates of the multi-type GWI process in question.
Appendix \ref{CMT} contains a version of the continuous mapping theorem and some related results.
In Appendix \ref{section_conv_step_processes} we recall a result about convergence of random step processes
 towards a diffusion process due to Isp\'any and Pap \cite{IspPap},
 and a result about the asymptotic behavior of a critical single-type GWI process due to Wei and Winnicki \cite[Theorem 2.1]{WW}.

Next, we summarize the novelties of the paper.
We point out that only few results are available for functional limit theorems for critical decomposable multi-type GWI processes,
 in fact, we can mention only Barczy et al.\ \cite{BarBezPap2} about critical decomposable 2-type GWI processes.
The presence of a 2-fold iterated integral process of a squared Bessel process in the limit distribution in Theorem \ref{main_4}
 of the present paper is a new phenomenon compared to the limit distributions that can be found in Barczy et al.\ \cite{BarBezPap2}.
Furthermore, Lemma \ref{ConvIteratedIntegral} may be of interest on its own right, since
  it may be used to describe the joint distribution of a continuous stochastic process, its integral process
 and its 2-fold iterated integral process, provided that the continuous stochastic process in question
 is the weak limit of a sequence of some c\`{a}dl\`{a}g stochastic processes.

In the present paper, we restrict our investigation to decomposable $3$-type GWI processes with (lower) triangular offspring matrices
 where each diagonal entry is $1$ (i.e. the strongly critical case).
As a future topic, instead of the asymptotic behavior of the processes considered here,
 we plan to investigate the asymptotic behavior of a general strongly critical decomposable $p$-type GWI process, where $p\in\NN$.
In the present paper we chose to focus only on the $3$-type (lower) triangular case with each diagonal entry $1$,
 since it is not clear to us how the limit distributions might look like in the general $p$-type case, and we expect that
 the moment estimates and continuous mapping type results found in Appendices \ref{GWI_moments} and \ref{CMT}
 should be generalized further.
We also mention that our choice of the offspring mean matrix is not very far away from the one that are used
 in biological applications (see, e.g., Balelli et al.\ \cite[Propositions 2 and 3]{BalMilWai}).

\section{Multi-type GWI processes}\label{Section_multi-type_GWI}

Let $\ZZ_+$, $\NN$, $\RR$, $\RR_+$ and $\RR_{++}$ denote the set
 of non-negative integers, positive integers, real numbers, non-negative real
 numbers and positive real numbers, respectively.
For $x,y\in\RR$, the minimum and maximum of $x$ and $y$ is denoted by $x\wedge y$  and $x\vee y$, respectively.
For $k,m\in\ZZ_+$ with $k<m$, we define $\binom{k}{m}:=0$.
For $i,j\in\ZZ_+$, let $\delta_{i,j}:=0$ if $i\ne j$, and $\delta_{i,j}:=1$ if $i=j$ (known as the Kronecker delta).
 For functions $f:\NN\to\RR$ and $g:\NN\to\RR$, the notation $f(k) = \OO(g(k))$, $k\in\NN$, means that
 there exists a positive real number $M$ such that $\vert f(k)\vert\leq M g(k)$ for each $k\in\NN$.
More generally, for $d\in\ZZ_+$ and functions $G:\RR^{d+1}\to\RR$, $g_1,\dots,g_d:\NN\to\RR$, and $h:\NN\to\RR$, the notation
 \begin{equation}\label{help_bigO}
  G(k,\OO(g_1(k)),\dots,\OO(g_d(k)))=\OO(h(k)),\qquad k\in\NN,
 \end{equation}
means that for all $c_1,\dots,c_d\in\RR_+$,
 there exists a positive real number $M$ such that
\begin{equation*}
|G(k, c_1 g_1(k), \dots, c_d g_d(k))|\leq M h(k),\qquad k\in\NN.
\end{equation*}
For example, $\OO(1)=\OO(k)$, $k\in\NN$, because for all $c\in\RR_+$, $|c\cdot1|\leq c\cdot k$, $k\in\NN$,
but we point out that $\OO(k)=\OO(1)$, $k\in\NN$, does not hold.
The relation $\OO(1)=\OO(k)$, $k\in\NN$, fits into the framework \eqref{help_bigO} with the following choices:
 $d=1$, $G:\RR^2\to \RR$, $G(x,y)=y$, $(x,y)\in\RR^2$, $g_1:\NN\to\RR$, $g_1(k)=1$, $k\in\NN$, and
 $h:\NN\to\RR$, $h(k)=k$, $k\in\NN$.
The relation $f(k) = \OO(g(k))$, $k\in\NN$, also fits into the framework \eqref{help_bigO} with the following choices:
 $d=1$, $G:\RR^2\to \RR$ satisfying $G(k,y)=f(k)$, $(k,y)\in\NN\times\RR$, $g_1$ can be arbitrary,
 and $h:\NN\to\RR$, $h(k)=g(k)$, $k\in\NN$; or with the choices
 $d=1$, $G:\RR^2\to \RR$, $G(x,y)=y$, $(x,y)\in\RR^2$, $g_1:\NN\to\RR$, $g_1(k):=f(k)$, $k\in\NN$,
 and $h:\NN\to\RR$, $h(k)=g(k)$, $k\in\NN$.
For functions  $f:\NN\to\RR$ and $g:\NN\to\RR$, by the notation $f(k)=\OO(g(k))\to0$ as $k\to\infty$,
  we mean that $f(k)=\OO(g(k))$, $k\in\NN$, and $\lim_{k\to\infty}g(k)=0$, which imply that $\lim_{k\to\infty}f(k)=0$.

The Euclidean norm on $\RR^d$ is denoted by $\Vert\cdot\Vert$, where $d\in\NN$.
The $d\times d$ identity matrix is denoted by $\bI_d$.
For a matrix $\bB\in\RR^{d\times d}$ and $m\in\ZZ_+$, the $(i,j)$-th entry of $\bB^m$ is denoted by $b_{i,j}^{[m]}$,
 where $i,j\in\{1,\ldots,d\}$.
For a function $f:\RR\to\RR$, its positive part is denoted by $f^+$.
Every random variable will be defined on a fixed probability space $(\Omega, \cA, \PP)$.
Convergence in probability, convergence in $L_1$, convergence almost surely, equality in distribution and
 almost sure equality is denoted by $\stoch$, $\mean$, $\as$, $\distre$ and $\ase$, respectively.
For $d\in\NN$, we will use $\DD(\RR_+,\RR^d)$ for the set of c\`adl\`ag functions from $\RR_+$ to $\RR^d$,
and $\CC(\RR_+,\RR^d)$ for the space of continuous functions from $\RR_+$ to $\RR^d$.
For a function $f\in\DD(\RR_+,\RR^d)$, the size of the jump of $f$ at time $t\in\RR_+$ is denoted by $\Delta f(t)$,
that is, $\Delta f(0):=0$ and $\Delta f (t):=f(t)-f(t-)$, $t\in\RR_{++}$,
where $f(t-)$ is the left-sided limit of $f$ at the point $t\in\RR_{++}$.
For $d\in\NN$, the convergence in the sense of the Skorokhod $J_1$
topology of functions in $\DD(\RR_+,\RR^d)$,
 the weak convergence of the finite dimensional distributions of $\RR^d$-valued stochastic processes  with sample paths in $\DD(\RR_+, \RR^d)$,
and the weak convergence of $\RR^d$-valued stochastic processes with sample paths in $\DD(\RR_+, \RR^d)$
are denoted by $\skor{d}$, $\distrf$, and $\distr$, respectively
(for more details and notations, e.g., for the metric that generates the Skorokhod $J_1$ topology and
 for the definition of locally uniform convergence $\lu$, see Appendix \ref{CMT}).

First, we recall the definition and first order moment formulae of $p$-type GWI processes, where $p\in\NN$.
For each $k \in \ZZ_+$ and $i \in \{ 1, \dots, p \}$, the number of
 individuals of type $i$ in the $k^\mathrm{th}$ generation is denoted by $X_{k,i}$.
For simplicity, we suppose that the initial values are $X_{0,i} = 0$, $i \in \{ 1, \dots, p \}$.
By $\xi_{k,j,i,\ell}$ we denote the number of type $\ell$ offsprings
 produced by the $j^\mathrm{th}$ individual who is of type $i$
 belonging to the $(k-1)^\mathrm{th}$ generation.
The number of type $i$ immigrants in the $k^\mathrm{th}$ generation is denoted by $\vare_{k,i}$.
Consider the random vectors
 \[
   \bX_k := \begin{bmatrix}
             X_{k,1} \\
             \vdots \\
             X_{k,p}
            \end{bmatrix} , \qquad
   \bxi_{k,j,i} := \begin{bmatrix}
                    \xi_{k,j,i,1} \\
                    \vdots \\
                    \xi_{k,j,i,p}
                   \end{bmatrix} , \qquad
   \bvare_k := \begin{bmatrix}
                \vare_{k,1} \\
                \vdots \\
                \vare_{k,p}
               \end{bmatrix} .
 \]
Then we have
 \begin{equation}\label{MBPI(d)}
  \bX_k = \sum_{i=1}^p \sum_{j=1}^{X_{k-1,i}} \bxi_{k,j,i} + \bvare_k , \qquad
  k \in \NN,
 \end{equation}
 with $\bX_0=\bzero$ (and using the convention $\sum_{i=1}^0:=\bzero$).
Here $\big\{\bxi_{k,j,i}, \, \bvare_k : k, j \in \NN, \, i \in \{ 1, \dots, p \} \big\}$
  are supposed to be independent.
Moreover, $\big\{\bxi_{k,j,i} : k, j \in \NN\big\}$ for each
 $i \in \{1, \dots, p\}$, and $\{\bvare_k : k \in \NN\}$ are
 supposed to consist of identically distributed $\ZZ_+^p$-valued random vectors.
For notational convenience, let $\{\bxi_i : i \in \{1, \ldots, p\}\}$ and $\bvare$ be random vectors
 such that $\bxi_i \distre \bxi_{1,1,i}$ for each $i \in \{1, \ldots, p\}$ and $\bvare \distre \bvare_1$.

In all what follows, we will suppose
\begin{equation}\label{baseconds}
\EE(\|\bxi_i\|^2) < \infty, \quad i\in\{ 1,\dots, p\}, \quad\text{and}\quad \EE (\|\bvare\|^2) < \infty.
\end{equation}
Introduce the notations
 \begin{gather}\label{help_jelolesek}
  \begin{split}
 \bA := \begin{bmatrix}
          \EE(\bxi_1) & \cdots &           \EE(\bxi_p)
         \end{bmatrix} \in \RR^{p \times p}_+ , \qquad
  \bb := \EE(\bvare) \in \RR^p_+ , \\
  \bV^{(i)} := \var(\bxi_i) \in \RR^{p \times p} , \quad i\in\{1,\dots,p\},\qquad
  \bV^{(0)} := \var(\bvare) \in \RR^{p \times p}.
  \end{split}
 \end{gather}
The matrix $\bA$ and the vector $\bb$ are called the offspring mean matrix and the immigration mean vector, respectively.
Note that $\bA=(a_{i,j})_{i,j=1}^p$ with $a_{i,j} := \EE(\xi_{1,1,j,i})$, $i,j\in\{1,\ldots,p\}$,
 and we mention that some authors define the offspring mean matrix as \ $\bA^\top$.

For each $k \in \ZZ_+$, let $\cF_k^\bX := \sigma(\bX_0,\dots, \bX_k)$, where $\cF_0^\bX=\{\emptyset,\Omega\}$ due to $\bX_0=\bzero$.
By \eqref{MBPI(d)}, we get
 \begin{equation}\label{mart}
  \EE(\bX_k \mid \cF_{k-1}^\bX)
  = \sum_{i=1}^p X_{k-1,i} \EE(\bxi_i) + \bb
  = \bA \bX_{k-1} + \bb,\qquad k\in\NN .
 \end{equation}
Consequently,
 \begin{equation*}%\label{recEX}
  \EE(\bX_k) = \bA \EE(\bX_{k-1}) + \bb , \qquad k \in \NN ,
 \end{equation*}
and, since $\bX_0=\bzero$, we have
 \begin{equation*}%\label{EXk}
  \EE(\bX_k) = \sum_{j=0}^{k-1} \bA^j \bb , \qquad k \in \ZZ_+.
 \end{equation*}
By looking at \eqref{mart}, we can interpret $\bA$ as a kind of scaling matrix that determines
 how many individuals we may expect on average to appear as the offsprings of individuals of the previous generation.
It suggests that the offspring mean matrix $\bA$ plays a crucial role in the
 asymptotic behavior of the sequence $(\bX_k)_{k\in\ZZ_+}$.
A $p$-type GWI process $(\bX_k)_{k\in\ZZ_+}$ is referred to respectively as \emph{subcritical}, \emph{critical} or
 \emph{supercritical} if $\varrho(\bA) < 1$, $\varrho(\bA) = 1$
 or $\varrho(\bA) > 1$, where $\varrho(\bA)$ denotes the spectral radius of the matrix $\bA$
 (see, e.g., Athreya and Ney \cite[V.3]{AN} or Quine \cite{Q}), and it is called indecomposable or decomposable
 if $\bA$ is irreducible or reducible, respectively.
Recall that $\bA$ is called reducible if there exists a permutation matrix $\bP\in\RR^{p\times p}$
 and $q\in\{1,\ldots,p-1\}$ (and hence $p$ should be at least $2$) such that
\begin{equation}\label{reducMatrix}
\bA=\bP\begin{bmatrix}\bR&\bzero\\ \bS & \bT\end{bmatrix}\bP^\top,
\end{equation}
where $\bR\in\RR^{q\times q},\bS\in\RR^{(p-q)\times q},\bT\in\RR^{(p-q)\times(p-q)}$, and $\bzero\in\RR^{q\times(p-q)}$ is a null matrix.
The matrix $\bA$ is called irreducible if it is not reducible, see, e.g., Horn and Johnson \cite[Definitions 6.2.21 and 6.2.22]{HJ}.
We do emphasize that no 1-by-1 matrix is reducible.
If a $p$-type GWI process is decomposable, then there exists a decomposition of $\{1,\dots,p\}$ into two disjoint subsets
 $C_1$ and $C_2$ such that for each $i\in C_1$ and $j\in C_2$, the individuals of type $j$ cannot produce offspring of type $i$.
If $\bA$ takes the form \eqref{reducMatrix} with $\bP:=\bI_p$, then appropriate choices for $C_1$ and $C_2$
 are the types corresponding to $\bR$ and $\bT$, respectively.

We recall a representation of nonnegative reducible matrices, see, e.g., Foster and Ney \cite{FosNey}.
Let $p\geq2$, $p\in\NN$, and $\bA\in\RR_+^{p\times p}$.
If $\bA$ is reducible, then, by definition, we can find a permutation matrix $\bP\in\RR^{p\times p}$ and $q\in\{1,\ldots,p-1\}$
such that \eqref{reducMatrix} holds.
The matrices $\bR$ and $\bT$ are nonnegative square matrices, thus they are either irreducible or reducible.

If $\bR$ is reducible, then we can find a permutation matrix $\bP_\bR\in\RR^{q\times q}$ and $r\in\{1,\dots,q-1\}$ such that
\begin{equation*}
\bR=\bP_\bR\begin{bmatrix}\bQ&\bzero\\ \bU & \bV\end{bmatrix}\bP_\bR^\top,
\end{equation*}
where $\bQ\in\RR^{r\times r},\bU\in\RR^{(q - r)\times r},\bV\in\RR^{(q-r)\times(q-r)}$, and $\bzero\in\RR^{r\times(q-r)}$ is a null matrix.
Letting
\begin{equation*}
\bP_1:=\bP\begin{bmatrix}\bP_\bR&\bzero\\ \bzero & \bI_{p-q}\end{bmatrix},\qquad
\widetilde\bR:=\begin{bmatrix}\bQ&\bzero\\ \bU & \bV\end{bmatrix},\qquad
\widetilde\bS:=\bS\bP_\bR=:\begin{bmatrix}\bS_1 & \bS_2\end{bmatrix},
\end{equation*}
 where $\bS_1\in\RR^{(p-q)\times r}$ and $\bS_2\in\RR^{(p-q)\times (q-r)},$ we have
 $\bP_\bR \widetilde\bR \bP_\bR^\top = \bR$ and $\widetilde\bS\bP_\bR^\top = \widetilde\bS \bP_\bR^{-1} = \bS$, and hence
\begin{align*}
\begin{bmatrix}\bP_\bR&\bzero\\ \bzero & \bI_{p-q}\end{bmatrix}
\begin{bmatrix}
\widetilde\bR & \bzero\\
\widetilde\bS & \bT
\end{bmatrix}
\begin{bmatrix}\bP_\bR&\bzero\\ \bzero & \bI_{p-q}\end{bmatrix}^\top
=&\begin{bmatrix}
\bP_\bR\widetilde\bR & \bzero\\
\widetilde\bS & \bT
\end{bmatrix}
\begin{bmatrix}\bP_\bR^\top&\bzero\\ \bzero & \bI_{p-q}\end{bmatrix}\\
=&\begin{bmatrix}
\bP_\bR\widetilde\bR\bP_\bR^\top & \bzero\\
\widetilde\bS\bP_\bR^\top & \bT
\end{bmatrix}=\begin{bmatrix}\bR&\bzero\\ \bS & \bT\end{bmatrix}.
\end{align*}
Thus we have
\begin{equation*}
\bP_1\begin{bmatrix}
\widetilde\bR & \bzero\\
\widetilde\bS & \bT
\end{bmatrix}\bP_1^\top=\bP\begin{bmatrix}
\bR & \bzero\\
\bS & \bT
\end{bmatrix}\bP^\top=\bA,
\end{equation*}
that is,
\begin{equation*}
\bA=\bP_1\begin{bmatrix}
\bQ & \bzero & \bzero\\
\bU & \bV & \bzero\\
\bS_1 & \bS_2 & \bT
\end{bmatrix}\bP_1^\top.
\end{equation*}

If $\bT$ is reducible, then a similar procedure leads us to the representation
\begin{equation*}
\bA=\bP_1'\begin{bmatrix}
\bR & \bzero & \bzero\\
\bS_1' & \bQ' & \bzero\\
\bS_2' & \bU' & \bV'
\end{bmatrix}(\bP_1')^\top,
\end{equation*}
where the matrices $\bP_1'$, $\bQ'$, $\bU'$, $\bV'$, $\bS_1'$ and $\bS_2'$ are such that
\begin{equation*}
\bT=\bP_\bT\begin{bmatrix}
\bQ' & \bzero\\
\bU' & \bV'
\end{bmatrix}\bP_\bT^\top,\qquad\bP_1':=\bP\begin{bmatrix}\bI_q&\bzero\\ \bzero & \bP_\bT\end{bmatrix},\qquad
\begin{bmatrix}\bS_1' \\ \bS_2'\end{bmatrix}:=\bP_\bT^\top\bS
\end{equation*}
with an appropriate permutation matrix $\bP_\bT$.

Repeating the procedure above for each reducible matrix in the diagonal of $\bP_1^\top\bA\bP_1$ or $(\bP_1')^\top\bA\bP_1'$,
 we can find a permutation matrix $\bP_\bA\in\RR^{p\times p}$ and $N\in\NN$ such that
\begin{equation}\label{reducTeljesFelb}
\bA=\bP_\bA\begin{bmatrix}
\bA_{1,1} & \bzero & \dots & \bzero\\
\bA_{2,1} & \bA_{2,2} & \dots & \bzero\\
\vdots & \vdots & \ddots & \vdots\\
\bA_{N,1} & \bA_{N,2} & \dots & \bA_{N,N}
\end{bmatrix}\bP_\bA^\top,
\end{equation}
where $\bA_{i,j}\in\RR^{p_i\times p_j}$ for some $p_i,p_j\in\{1,\dots,p\}$, $i,j\in\{1,\dots,N\}$, $p_1+\dots+p_N=p$, and $\bA_{i,i}$, $i\in\{1,\dots,N\}$ are all irreducible.
This is true when $\bA$ is irreducible as well, as in that case $\bP_\bA=\bI_p$ and $N=1$.

Let $(\bX_k)_{k\in\ZZ_+}$ be a $p$-type GWI process such that its offspring mean matrix $\bA$ has the form given in \eqref{reducTeljesFelb} with $\bP_\bA:=\bI_p$.
In fact, one can assume that $\bP_\bA$ equals $\bI_p$, since otherwise instead of $(\bX_k)_{k\in\ZZ_+}$ one can consider the $p$-type GWI process
 $(\bP_\bA^\top\bX_k)_{k\in\ZZ_+}$, which has an offspring mean matrix $\bA' = \bP_\bA^\top\bA\bP_\bA$ having the form
  \eqref{reducTeljesFelb} with $\bP_{\bA'}:=\bI_p$.
 Then we can group the types $\{1,\dots,p\}$ into pairwise disjoint sets $C_1,\dots,C_N$, where $C_1:=\{1,\dots,p_1\}$
 and $C_{i+1}:=\{\sum_{j=1}^{i}p_j+1,\dots,\sum_{j=1}^{i}p_j+p_{i+1}\}$, $i=1,\dots,N-1$, that is,
 $C_i$ consists of the types corresponding to the matrix $\bA_{i,i}$, $i=1,\ldots,N$.
For $i,j\in\{1,\dots,p\}$, we say that type $j$ is accessible from type $i$ (in notation $i\to j$)
  if the $(j,i)$-th entry of $\bA^\ell$ is positive for some $\ell\in\NN$.
This means that individuals of type $i$ can have descendents of type $j$.
If $i\to j$ and $j\to i$, we say that the types $i$ and $j$ communicate (in notation $i \leftrightarrow j$).
The relation $\leftrightarrow$ (i.e., communication between types) is an equivalence relation,
 and the subsets $C_1,\dots,C_N$ of $\{1,\ldots,p\}$ introduced above can be considered as the partition corresponding to $\leftrightarrow$.

We say that $(\bX_k)_{k\in\ZZ_+}$ is strongly critical if $\varrho(\bA_{i,i})=1$ for each $i\in\{1,\ldots,N\}$,
 where $\bA_{i,i}$, $i\in\{1,\ldots,N\}$, are the irreducible matrices in the decomposition \eqref{reducTeljesFelb} of $\bA$ (see, e.g., Foster and Ney \cite{FosNey}).
In this paper, we examine special strongly critical, decomposable $3$-type GWI processes, with offspring mean matrix $\bA$ given in \eqref{reducTeljesFelb} with $\bP_\bA:=\bI_3$ and $N:=3$.
In this case, we have $p_1=p_2=p_3=1$, that is, $\bA_{i,i}:=[a_{i,i}]$, $i\in\{1,2,3\}$, yielding that $\bA$ is lower triangular,
 and, since $\bA\in\RR_+^{3\times3}$, we also have $a_{i,i}=|a_{i,i}|=\rho(\bA_{i,i})=1$, $i\in\{1,2,3\}$.
However, we call attention to the fact that not all the strongly critical, decomposable $3$-type GWI processes
 have this type of offspring mean matrix.
Studying the asymptotic behavior of 3-type GWI processes with such more general offspring mean matrices could be a topic of future research.

Finally, we note that Generalized Integer-valued Autoregressive (GINAR) processes of order $p\in\NN$
 (generalizations of Integer-valued AutoRegressive (INAR) processes of order $p$, first introduced by Latour \cite{Latour})
 are special cases of $p$-type GWI processes, see, e.g., Barczy et al.\ \cite[Section 3]{BarNedPap}.
The offspring mean matrix of GINAR processes has a special form, for example, 
 in case of $p=3$, its second and third rows are $(1,0,0)$ and $(0,1,0)$, respectively.
Consequently, a GINAR process of order $3$ such that $a_{1,3} = \EE(\xi_{1,1,3,1})>0$ is irreducible 
 (see, e.g., Barczy et al.\ \cite[Remark 3]{BarNedPap}), and thus it does not belong 
 to the class of branching processes investigated in this paper.
Furthermore, using again Barczy et al.\ \cite[Remark 3]{BarNedPap}, 
 the offspring mean matrix of a critical GINAR process of order $3$ such that $a_{1,3}=0$ takes the form 
 \[
 \begin{bmatrix}
  a_{1,1} & a_{1,2} & 0 \\
  1 & 0 & 0 \\
  0 & 1 & 0 \\
\end{bmatrix}, \qquad \text{where \ $a_{1,1}+a_{1,2}=1$.}
 \]
This matrix is reducible and has spectral radius $1$, yielding that such a GINAR process is critical and decomposable,
 but it is not strongly critical, and hence it does not belong to the class of branching processes that we investigate.

\section{Asymptotic behavior of random step processes formed from some strongly critical, decomposable 3-type GWI pro\-cesses}
\label{Section_conv_results}

In what follows, we consider a 3-type Galton-Watson process $(\bX_k)_{k\in\ZZ_+}$ with
 immigration starting from $\bX_0=\bzero$, and we suppose that the moment conditions \eqref{baseconds} hold
 and that the offspring mean matrix $\bA$ is lower triangular with diagonal entries $1$.
We point out that one could assume equivalently that $\bA$ is upper triangular with diagonal entries $1$.
Indeed, if $([X_{k,1}, X_{k,2}, X_{k,3}]^\top)_{k\in\ZZ_+}$ is a 3-type GWI process having an offspring
 mean matrix with $(1,2)$-, $(1,3)$- and $(2,3)$-entries $0$, then $([X_{k,3},X_{k,2},X_{k,1}]^\top)_{k\in\ZZ_+}$ is a 3-type GWI process
 having an offspring mean matrix with $(2,1)$-, $(3,1)$- and $(3,2)$-entries $0$.
 Because of this, our forthcoming Theorems \ref{main_1} and \ref{main_2}--\ref{main_4} could be reformulated for
  GWI processes with upper triangular offspring mean matrices (instead of lower triangular ones) as well.

Then the offspring mean matrix $\bA$ and the immigration mean vector $\bb$ of $(\bX_k)_{k\in\ZZ_+}$ take the following forms:
 \[
   \bA
       = \begin{bmatrix}
                 \EE(\bxi_1) & \EE(\bxi_2)  & \EE(\bxi_3)
         \end{bmatrix}
   = \begin{bmatrix}
      1 & 0 & 0 \\
      a_{2,1} & 1 & 0 \\
      a_{3,1} & a_{3,2} & 1
     \end{bmatrix}  \qquad \text{and} \qquad
   \bb = \begin{bmatrix} b_1 \\ b_2 \\ b_3 \end{bmatrix}.
 \]
Taking into account that $a_{1,2}=a_{1,3}=a_{2,3}=0$ implies that $\xi_{1,1,2,1}\ase0$, $\xi_{1,1,3,1}\ase0$ and $\xi_{1,1,3,2}\ase0$,
 Equation \eqref{MBPI(d)} with $p=3$ takes the form
 \begin{align}\label{help1}
   \begin{bmatrix} X_{k,1} \\ X_{k,2} \\ X_{k,3} \end{bmatrix}
   = \sum_{j=1}^{X_{k-1,1}}
      \begin{bmatrix} \xi_{k,j,1,1} \\ \xi_{k,j,1,2} \\ \xi_{k,j,1,3} \end{bmatrix}
     + \sum_{j=1}^{X_{k-1,2}}
        \begin{bmatrix} 0 \\ \xi_{k,j,2,2} \\ \xi_{k,j,2,3} \end{bmatrix}
     + \sum_{j=1}^{X_{k-1,3}}
        \begin{bmatrix} 0 \\ 0 \\ \xi_{k,j,3,3} \end{bmatrix}
     + \begin{bmatrix} \vare_{k,1} \\ \vare_{k,2} \\ \vare_{k,3} \end{bmatrix} , \qquad
   k \in \NN ,
 \end{align}
 with \ $[X_{0,1}, X_{0,2}, X_{0,3}]^\top=\bzero$.
Using the notion of accessibility of types of multi-type GWI processes
 (introduced in Section \ref{Section_multi-type_GWI}), the form of the offspring mean matrix $\bA$
 shows that type $1$ is not accessible from either type $2$ or type $3$, and type $2$ is not accessible from type $3$.
In other words (one can also see it from \eqref{help1}), individuals of type $3$ can only have offsprings of type $3$,
 thus they can never have descendants of type $1$ or $2$,
 and similarly, individuals of type $2$ can only have offsprings of either type $2$ or type $3$,
 yielding that they can never have descendants of type $1$.
We check that it is enough to consider only the following 4 cases
($a_{1,2}=a_{1,3}=a_{2,3}=0$ and $a_{1,1}=a_{2,2}=a_{3,3}=1$ for each case):
\begin{align}\label{tablazat_esetek}
 \begin{tabular}{|c|c|c|c|}
  \hline
  \textup{(1)} & $a_{2,1} = 0$ & $a_{3,1} = 0$ & $a_{3,2} = 0$ \\
  \hline
  \textup{(2)} & $a_{2,1} = 0$ & $a_{3,1} \in \RR_{++}$  & $a_{3,2} \in \RR_{+}$ \\
  \hline
  \textup{(3)} & $a_{2,1} \in \RR_{++}$ & $a_{3,1} \in \RR_{++}$ & $a_{3,2}=0$ \\
  \hline
  \textup{(4)} & $a_{2,1} \in \RR_{++}$ &  $a_{3,1}\in \RR_+$ & $a_{3,2} \in \RR_{++}$ \\
  \hline
 \end{tabular}\;.
\end{align}
For abbreviation, we can write the above four cases in matrix form as follows:
 \begin{align*}
   & \begin{bmatrix}
      1 & 0 & 0 \\
      0 & 1 & 0 \\
      0 & 0 & 1
     \end{bmatrix}_1 , \qquad
    \begin{bmatrix}
      1 & 0 & 0 \\
      0 & 1 & 0 \\
      ++ & + & 1
     \end{bmatrix}_2 , \qquad
    \begin{bmatrix}
      1 & 0 & 0 \\
      ++ & 1 & 0 \\
      ++ & 0 & 1
     \end{bmatrix}_3 , \qquad
   \begin{bmatrix}
      1 & 0 & 0 \\
      ++ & 1 & 0 \\
      + & ++ & 1
     \end{bmatrix}_4 .
 \end{align*}
To show that it is sufficient to consider only these 4 cases, consider a 3-type GWI process $([X_{k,1},X_{k,2},X_{k,3}]^\top)_{k\in\ZZ_+}$
whose offspring mean matrix is lower triangular with diagonal entries all $1$,
 and assume that this process does not fit into any of the cases listed in \eqref{tablazat_esetek}.
It is easy to see that (by considering which entries below the diagonal can be $0$, and noticing that all but two possibilities are covered by the cases \textup{(1)}-\textup{(4)})
  such a process must belong to one of the following two cases \textup{(a)} or \textup{(b)}, written in matrix form as
\begin{equation*}
\begin{bmatrix}
      1 & 0 & 0 \\
      0 & 1 & 0 \\
      0 & ++ & 1
     \end{bmatrix}_a,\qquad
\begin{bmatrix}
      1 & 0 & 0 \\
      ++ & 1 & 0 \\
      0 & 0 & 1
     \end{bmatrix}_b.
\end{equation*}
If the process $([X_{k,1},X_{k,2},X_{k,3}]^\top)_{k\in\ZZ_+}$ belongs to case \textup{(a)},
then the process $([X_{k,2},X_{k,1},X_{k,3}]^\top)_{k\in\ZZ_+}$ has an offspring mean matrix with
$(1,2)$-, $(1,3)$-, $(2,1)$-, $(2,3)$- and $(3,2)$-entries $0$, $(1,1)$-, $(2,2)$- and $(3,3)$-entries $1$, and $(3,1)$-entry positive,
 yielding that it belongs to case \textup{(2)}.
Similarly, if the process $([X_{k,1},X_{k,2},X_{k,3}]^\top)_{k\in\ZZ_+}$ belongs to case \textup{(b)},
then the process $([X_{k,1},X_{k,3},X_{k,2}]^\top)_{k\in\ZZ_+}$ has an offspring mean matrix with
$(1,2)$-, $(1,3)$-, $(2,1)$-, $(2,3)$- and $(3,2)$-entries $0$, $(1,1)$-, $(2,2)$- and $(3,3)$-entries $1$, and $(3,1)$-entry positive,
 yielding that it also belongs to case \textup{(2)}.

Note that \eqref{help1} readily yields that the first coordinate process $(X_{k,1})_{k\in\ZZ_+}$ of $(\bX_k)_{k\in\ZZ_+}$ satisfies
 \begin{align*}%\label{GWI_def1}
   X_{k,1} = \sum_{j=1}^{X_{k-1,1}} \xi_{k,j,1,1} + \vare_{k,1} , \qquad k \in \NN ,
 \end{align*}
 hence $(X_{k,1})_{k\in\ZZ_+}$ is a critical single-type GWI process
 (due to $\EE(\xi_{1,1,1,1})=1$).
Consequently, by a result of Wei and Winnicki \cite{WW} (see also Theorem \ref{thm:critical}), we have
 \begin{align}\label{help_Wei_Winnicki1}
   (n^{-1} X_{\nt,1})_{t\in\RR_+} \distr (\cX_{t,1})_{t\in\RR_+} \qquad
  \text{as \ $n \to \infty$,}
 \end{align}
 where the limit process $(\cX_{t,1})_{t\in\RR_+}$ is the pathwise unique strong solution of the SDE
 \[
   \dd \cX_{t,1} = b_1 \, \dd t + \sqrt{v^{(1)}_{1,1} \, \cX_{t,1}^+} \, \dd \cW_{t,1} , \qquad t\in\RR_+, \qquad \cX_{0,1} = 0 ,
 \]
 where $(\cW_{t,1})_{t\in\RR_+}$ is a standard Wiener process, $b_1= \EE(\vare_{1,1})$ and $v^{(1)}_{1,1}:=\var(\xi_{1,1,1,1})$.
The process $(\cX_{t,1})_{t\in\RR_+}$ is called a squared Bessel process.

If $a_{2,1}=0$ holds, then $\xi_{1,1,1,2}\ase0$, and
 \eqref{help1} yields that in this case the second coordinate process $(X_{k,2})_{k\in\ZZ_+}$  of $(\bX_k)_{k\in\ZZ_+}$ satisfies
 \begin{align*}%\label{GWI_def2}
   X_{k,2} = \sum_{j=1}^{X_{k-1,2}} \xi_{k,j,2,2} + \vare_{k,2} , \qquad k \in \NN .
 \end{align*}
Hence if $a_{2,1}=0$ holds, then $(X_{k,2})_{k\in\ZZ_+}$ is a critical single-type GWI process
 (due to $\EE(\xi_{1,1,2,2})=1$).

If $a_{3,1}=a_{3,2}=0$ holds, then $\xi_{1,1,1,3}\ase0$ and $\xi_{1,1,2,3}\ase0$,
 and \eqref{help1} yields that in this case the third coordinate process $(X_{k,3})_{k\in\ZZ_+}$ of $(\bX_k)_{k\in\ZZ_+}$  satisfies
 \[
   X_{k,3} = \sum_{j=1}^{X_{k-1,3}} \xi_{k,j,3,3} + \vare_{k,3} , \qquad k \in \NN .
 \]
Hence if $a_{3,1}=a_{3,2}=0$ holds, then $(X_{k,3})_{k\in\ZZ_+}$ is a critical single-type GWI process
 (due to $\EE(\xi_{1,1,3,3})=1$).

Next, we present our results on the asymptotic behavior of $(\bX_k)_{k\in\ZZ_+}$
 in the four cases (1)--(4) of its offspring mean matrix $\bA$.
The matrices $\bV^{(i)}$, $i\in\{0,1,2,3\}$ (introduced in \eqref{help_jelolesek})
 will be written in the form $\bV^{(i)}=(v_{k,\ell}^{(i)})_{k,\ell\in\{1,2,3\}}$, $i\in\{0,1,2,3\}$.

\begin{Thm}\label{main_1}
Let $(\bX_k)_{k\in\ZZ_+}$ be a 3-type GWI process such that $\bX_0=\bzero$,
 the moment conditions $\EE(\|\bxi_i\|^4) < \infty$, $i\in\{1,2,3\}$, and $\EE (\|\bvare\|^4) < \infty$
 hold, and its offspring mean matrix satisfies \textup{(1)} of \eqref{tablazat_esetek} (i.e., $\bA=\bI_3$).
Then we have
 \begin{equation}\label{Conv_X_1}
  \left(\begin{bmatrix}
          n^{-1} X_{\nt,1} \\
          n^{-1} X_{\nt,2}\\
          n^{-1} X_{\nt,3}
         \end{bmatrix}\right)_{t\in\RR_+}
  \distr
  \left(\begin{bmatrix}
          \cX_{t,1} \\
          \cX_{t,2}\\
          \cX_{t,3}
         \end{bmatrix}\right)_{t\in\RR_+}
  \qquad \text{as \ $n \to \infty$,}
 \end{equation}
 where the limit process is the pathwise unique strong solution of the SDE
 \begin{equation}\label{SDE_X_1}
  \begin{cases}
   &\dd\cX_{t,1}
   = b_1 \, \dd t
     + \sqrt{v^{(1)}_{1,1} \, \cX_{t,1}^+} \,
       \dd \cW_{t,1} , \\[1mm]
   &\dd\cX_{t,2}
   = b_2 \, \dd t
     + \sqrt{v^{(2)}_{2,2} \, \cX_{t,2}^+} \,
       \dd \cW_{t,2},\\[1mm]
   &\dd\cX_{t,3}
   = b_3 \, \dd t
     + \sqrt{v^{(3)}_{3,3} \, \cX_{t,3}^+} \,
       \dd \cW_{t,3} ,
  \end{cases}
  \qquad t \in\RR_+ ,
 \end{equation}
 with initial value $[\cX_{0,1} , \cX_{0,2} , \cX_{0,3}]^\top = \bzero$, where $(\cW_{t,1})_{t\in\RR_+}$, $(\cW_{t,2})_{t\in\RR_+}$ and $(\cW_{t,3})_{t\in\RR_+}$
 are independent standard Wiener processes yielding the independence of $(\cX_{t,1})_{t\in\RR_+}$, $(\cX_{t,2})_{t\in\RR_+}$ and $(\cX_{t,3})_{t\in\RR_+}$
 as well.
\end{Thm}

\begin{Rem}
Note that in Theorem \ref{main_1}, we do not suppose the independence of the coordinates
 of the immigration vector, but despite of this fact the coordinate processes $(\cX_{t,1})_{t\in\RR_+}$, $(\cX_{t,2})_{t\in\RR_+}$ and $(\cX_{t,3})_{t\in\RR_+}$
 of the limit process are independent.
Heuristically, one might expect this,
 as in the single-type case only the expected value of the immigration distribution appears in the corresponding limit process,
 see Theorem \ref{thm:critical}.
\proofend
\end{Rem}

\begin{Thm}\label{main_2}
Let $(\bX_k)_{k\in\ZZ_+}$ be a 3-type GWI process such that $\bX_0=\bzero$,
 the moment conditions $\EE(\|\bxi_i\|^4) < \infty$, $i\in\{1,2\}$,
  $\EE(\|\bxi_3\|^2) < \infty$, and $\EE (\|\bvare\|^4) < \infty$
 hold, and its offspring mean matrix $\bA$ satisfies \textup{(2)} of \eqref{tablazat_esetek}.
Then we have
 \begin{equation*}%\label{Conv_X_2}
  \left(\begin{bmatrix}
          n^{-1} X_{\nt,1} \\
          n^{-1} X_{\nt,2} \\
          n^{-2} X_{\nt,3}
         \end{bmatrix}\right)_{t\in\RR_+}
  \distr
  \left(\begin{bmatrix}
          \cX_{t,1} \\
          \cX_{t,2} \\
          \cX_{t,3}
         \end{bmatrix}\right)_{t\in\RR_+}
  \qquad \text{as \ $n \to \infty$,}
 \end{equation*}
 where the limit process is the pathwise unique strong solution of the SDE
 \begin{equation}\label{SDE_X_2}
   \begin{cases}
   \dd\cX_{t,1}
   = b_1 \, \dd t + \sqrt{v^{(1)}_{1,1} \cX_{t,1}^+} \, \dd \cW_{t,1} , \\[2mm]
   \dd\cX_{t,2}
   = b_2 \, \dd t + \sqrt{v^{(2)}_{2,2} \cX_{t,2}^+} \, \dd \cW_{t,2} , \\[2mm]
   \dd\cX_{t,3}
   = (a_{3,1} \cX_{t,1} + a_{3,2} \cX_{t,2}) \, \dd t ,
  \end{cases}
  \qquad t \in\RR_+,
  \end{equation}
 with initial value $[\cX_{0,1}, \cX_{0,2},\cX_{0,3}]^\top = \bzero$, where $(\cW_{t,1})_{t\in\RR_+}$ and $(\cW_{t,2})_{t\in\RR_+}$ are independent standard Wiener processes.
\end{Thm}

\begin{Thm}\label{main_3}
Let $(\bX_k)_{k\in\ZZ_+}$ be a 3-type GWI process such that $\bX_0=\bzero$,
 the moment condition \eqref{baseconds} holds, and its offspring mean matrix $\bA$ satisfies \textup{(3)} of \eqref{tablazat_esetek}.
Then we have
 \begin{equation*}%\label{Conv_X_3}
   \left(\begin{bmatrix}
          n^{-1} X_{\nt,1} \\
          n^{-2} X_{\nt,2} \\
          n^{-2} X_{\nt,3}
         \end{bmatrix}\right)_{t\in\RR_+}
  \distr
  \left(\begin{bmatrix}
          \cX_{t,1} \\
          \cX_{t,2} \\
          \cX_{t,3}
         \end{bmatrix}\right)_{t\in\RR_+}
  \qquad \text{as \ $n \to \infty$,}
 \end{equation*}
 where the limit process is the pathwise unique strong solution of the SDE
 \begin{equation}\label{SDE_X_3}
  \begin{cases}
   \dd\cX_{t,1}
   = b_1 \, \dd t + \sqrt{v^{(1)}_{1,1} \cX_{t,1}^+} \, \dd \cW_{t,1} , \\[2mm]
   \dd\cX_{t,2}
   = a_{2,1} \cX_{t,1} \, \dd t , \\[2mm]
   \dd\cX_{t,3}
   = a_{3,1} \cX_{t,1} \, \dd t ,
  \end{cases}
  \qquad t \in\RR_+ ,
   \end{equation}
 with initial value $ [\cX_{0,1}, \cX_{0,2},\cX_{0,3}]^\top = \bzero$, where $(\cW_{t,1})_{t\in\RR_+}$ is a standard Wiener process.
\end{Thm}

\begin{Thm}\label{main_4}
Let $(\bX_k)_{k\in\ZZ_+}$ be a 3-type GWI process such that $\bX_0=\bzero$,
 the moment condition \eqref{baseconds} holds, and its offspring mean matrix $\bA$ satisfies \textup{(4)} of \eqref{tablazat_esetek}.
Then we have
 \begin{equation*}%\label{Conv_X_5}
 \left(\begin{bmatrix}
          n^{-1} X_{\nt,1} \\
          n^{-2} X_{\nt,2} \\
          n^{-3} X_{\nt,3}
         \end{bmatrix}\right)_{t\in\RR_+}
  \distr
  \left(\begin{bmatrix}
          \cX_{t,1} \\
          \cX_{t,2} \\
          \cX_{t,3}
         \end{bmatrix}\right)_{t\in\RR_+}
  \qquad \text{as \ $n \to \infty$,}
 \end{equation*}
 where the limit process is the pathwise unique strong solution of the SDE
 \begin{equation}\label{SDE_X_5}
    \begin{cases}
   \dd\cX_{t,1}
   = b_1 \, \dd t + \sqrt{v^{(1)}_{1,1} \cX_{t,1}^+} \, \dd \cW_{t,1} , \\[2mm]
   \dd\cX_{t,2}
   = a_{2,1} \cX_{t,1} \, \dd t , \\[2mm]
   \dd\cX_{t,3}
   = a_{3,2} \cX_{t,2} \, \dd t ,
  \end{cases}
  \qquad t \in\RR_+,
 \end{equation}
 with initial value $[\cX_{0,1}, \cX_{0,2},\cX_{0,3}]^\top = \bzero$,
 where $(\cW_{t,1})_{t\in\RR_+}$ is a standard Wiener process.
\end{Thm}

We note that if $b_1=0$, then the limit process in Theorem \ref{main_4} is the identically zero process, that is,
 $[\cX_{t,1},\cX_{t,2},\cX_{t,3}]^\top=\bzero$, $t\in\RR_+$.
This is in accordance with part \text{(i)} of Theorem 4 in Foster and Ney \cite{FosNey} in the special case $\ba_1=b_1=0$.
In the next remark, we rewrite the limit process in Theorem \ref{main_4} in two different forms.

\begin{Rem}
The limit process in Theorem \ref{main_4} can be written in the following two forms
\begin{equation*}
\begin{bmatrix}
\cX_{t,1}\\
\cX_{t,2}\\
\cX_{t,3}
\end{bmatrix}=
\begin{bmatrix}
\cX_{t,1}\\
a_{2,1}\int_0^t\cX_{s,1}\,\dd s\\
a_{3,2}a_{2,1}\int_0^t\int_0^r\cX_{s,1}\,\dd s\,\dd r
\end{bmatrix}=
\begin{bmatrix}
\cX_{t,1}\\
a_{2,1}\int_0^t\cX_{s,1}\,\dd s\\
a_{3,2}a_{2,1}\int_0^t\int_0^t \bbone_{[0,r]}(s)\cX_{s,1}\,\dd s\,\dd r
\end{bmatrix}=
\begin{bmatrix}
\cX_{t,1}\\
a_{2,1}\int_0^t\cX_{s,1}\,\dd s\\
a_{3,2}a_{2,1}\int_0^t(t-s)\cX_{s,1}\,\dd s
\end{bmatrix}
\end{equation*}
 for $t\in\RR_+$, and
\begin{equation*}
\begin{bmatrix}
\cX_{t,1}\\
\cX_{t,2}\\
\cX_{t,3}
\end{bmatrix}
=\begin{bmatrix}
\cX_{t,1}\\
a_{2,1}\int_0^t(t-s)\,\dd\cX_{s,1}\\
\frac{a_{3,2}a_{2,1}}{2}\int_0^t(t-s)^2\,\dd\cX_{s,1}
\end{bmatrix},
\qquad t\in\RR_+.
\end{equation*}
The first formula follows from Fubini's theorem, and the second one can be checked using It\^o's formula as follows.
Namely, consider the function $f:\RR_+\times \RR\to\RR$, $f(t,x):=tx$, $t\in\RR_+$, $x\in\RR$.
Then, by It\^o's formula, we have
 \begin{equation}\label{help_Ito}
 t\cX_{t,1}=\int_0^t\cX_{s,1}\,\dd s+\int_0^ts\,\dd\cX_{s,1},\qquad t\in\RR_+.
 \end{equation}
Thus, by the SDE for $(\cX_{t,1})_{t\in\RR_+}$ (see the first equation of the SDE \eqref{SDE_X_5})
 and the fact that $\cX_{t,1} = \cX_{t,1} - 0 = \cX_{t,1} - \cX_{0,1} = \int_0^t1\,\dd\cX_{s,1}$, $t\in\RR_+$, we obtain that
 \begin{equation*}
 \int_0^t\cX_{s,1}\,\dd s=t\int_0^t1\,\dd \cX_{s,1}-\int_0^ts\,\dd\cX_{s,1}=\int_0^t(t-s)\,\dd\cX_{s,1}, \qquad t\in\RR_+,
 \end{equation*}
 as desired.
Furthermore, using again It\^o's formula at the third equality and \eqref{help_Ito} at the last equality,
 we have
 \begin{align*}
 \dd\left(\frac{1}{2} \int_0^t(t-s)^2\,\dd\cX_{s,1} \right)
   & = \dd\left( \frac{1}{2} \int_0^t (t^2 - 2ts + s^2) \,\dd \cX_{s,1} \right) \\
   & = \frac{1}{2} \dd \left( t^2 \cX_{t,1} - 2t \int_0^t s\,\dd \cX_{s,1} + \int_0^t s^2 \,\dd \cX_{s,1} \right)\\
   & = \frac{1}{2} \left( 2t \cX_{t,1} \,\dd t + t^2 \dd \cX_{t,1} - 2\left( \int_0^t s\,\dd\cX_{s,1} \right)\dd t
                        - 2t\cdot t\,\dd\cX_{t,1} + t^2 \dd\cX_{t,1} \right) \\
   & = t\cX_{t,1}\dd t - \left( \int_0^t s\,\dd \cX_{s,1} \right)\dd t
    = \left( \int_0^t \cX_{s,1} \,\dd s \right) \dd t,
    \qquad t\in\RR_+.
 \end{align*}
This implies that
 \[
  \int_0^t\int_0^r\cX_{s,1}\,\dd s\,\dd r = \frac{1}{2}\int_0^t(t-s)^2\,\dd\cX_{s,1},\qquad t\in\RR_+,
 \]
 as desired.
\proofend
\end{Rem}

In the next remark, we discuss some connections between Theorems \ref{main_1}, \ref{main_2}, \ref{main_3} and \ref{main_4}.

\begin{Rem}\label{Relax_Remark}
The conditions of Theorem \ref{main_4} may be relaxed by allowing that
 the offspring mean matrix $\bA$ satisfy any of the cases \textup{(1)}--\textup{(4)} of \eqref{tablazat_esetek},
 but if $\bA$ satisfies one of the cases \textup{(1)}--\textup{(3)} of \eqref{tablazat_esetek},
  then the second or third coordinate process of the limit process $([\cX_{t,1},\cX_{t,2},\cX_{t,3}])_{t\in\RR_+}$
  is identically zero.
Indeed, in the proof of Theorem \ref{main_4}, it is only in Step 2/(b) where we use the fact that $\bA$
satisfies case \textup{(4)} of \eqref{tablazat_esetek}. In this step, to prove convergence \eqref{case4side0},
we use that $\eta_1=1$, $\eta_2= 2$, and $\eta_3=3$
in case \textup{(4)} of \eqref{tablazat_esetek},
where $\eta_i$, $i\in\{1,2,3\}$, is defined in \eqref{etadef}.
However, the proof of convergence \eqref{case4side0} can be carried out even if
only the inequalities $\eta_1\leq 2$, $\eta_2\leq2$, and $\eta_3\leq4$ are satisfied.
By Lemma \ref{eta_values}, these inequalities for $\eta_i$, $i\in\{1,2,3\}$, hold in cases
\textup{(1)}--\textup{(3)} of \eqref{tablazat_esetek} as well, thus the conclusion of Theorem \ref{main_4} remains true under our relaxed conditions as well.
It is worth noting that if $\bA$ satisfies one of the cases \textup{(1)}--\textup{(3)} of \eqref{tablazat_esetek},
then we necessarily have $a_{3,2}a_{2,1}=0$, thus in these cases $\cX_{t,3}=0$, $t\in\RR_+$,
and if $\bA$ satisfies one of the cases \textup{(1)}--\textup{(2)} of \eqref{tablazat_esetek}, then
we have $a_{2,1}=0$, and thus $\cX_{t,2}=0$, $t\in\RR_+$, holds.
The reason for stating Theorem \ref{main_4} only when $\bA$ satisfies case \textup{(4)} of \eqref{tablazat_esetek}
is that in every other case, at least one coordinate of the limit process is identically zero.
One could similarly relax the conditions of Theorems \ref{main_1}, \ref{main_2}, and \ref{main_3} as well
(e.g. Theorem \ref{main_3} holds if the offspring mean matrix $\bA$
satisfies any of the cases \textup{(1)}--\textup{(3)} of \eqref{tablazat_esetek}).
\proofend
\end{Rem}

In the next remark, we provide a heuristic argument why the 1-fold and 2-fold iterated integral processes of $(\cX_{t,1})_{t\in\RR_+}$
 are expected to appear in the limit process in Theorem \ref{main_4}.

\begin{Rem}\label{Heuristic_Remark}
Let us suppose that the conditions of Theorem \ref{main_4} hold.
%If $(\bX_k)_{ k\in\ZZ_+}$ is a strongly critical decomposable 3-type GWI process whose offspring mean matrix is triangular and $\bX_0=\bzero$,
Then, using \eqref{help1}, we can write
\begin{align*}
 X_{k,2}
    =& \sum_{j=1}^{X_{k-1,1}}\xi_{k,j,1,2}+\sum_{j=1}^{X_{k-1,2}}\xi_{k,j,2,2}+\vare_{k,2}\\
	=&\sum_{j=1}^{X_{k-1,1}}\xi_{k,j,1,2}+\sum_{j=1}^{X_{k-1,2}}\xi_{k,j,2,2}+\vare_{k,2}+X_{k-1,1}(a_{2,1}-a_{2,1})+X_{k-1,2}(1-1)\\
	=&X_{k-1,2} + a_{2,1}X_{k-1,1} + \sum_{j=1}^{X_{k-1,1}}(\xi_{k,j,1,2}-a_{2,1})+\sum_{j=1}^{X_{k-1,2}}(\xi_{k,j,2,2}-1)+\vare_{k,2}
 \end{align*}
 for each \ $k\in\NN$.
Repeating the same procedure for $X_{k-1,2},\ldots,X_{2,1}$ (or formally, by induction) and using the fact that $X_{0,2}=0$, we get that
\begin{align*}
	X_{k,2}=&a_{2,1}\sum_{\ell=1}^{k-1}X_{\ell,1}+\sum_{\ell=1}^{k}\sum_{j=1}^{X_{\ell-1,1}}(\xi_{\ell,j,1,2}-a_{2,1})
		+\sum_{\ell=1}^{k}\sum_{j=1}^{X_{\ell-1,2}}(\xi_{\ell,j,2,2}-1)+\sum_{\ell=1}^k\vare_{\ell,2}\\
	=&a_{2,1}\int_0^{\frac{k}{n}}nX_{\ns,1}\,\dd s + E_{k,2}, \qquad  k,n\in\NN,
\end{align*}
where the error term $E_{k,2}$ is defined by
 \[
   E_{k,2}
      :=\sum_{\ell=1}^{k}\sum_{j=1}^{X_{\ell-1,1}}(\xi_{\ell,j,1,2}-a_{2,1})
	 	+\sum_{\ell=1}^{k}\sum_{j=1}^{X_{\ell-1,2}}(\xi_{\ell,j,2,2}-1)+\sum_{\ell=1}^k\vare_{\ell,2},
 \]
 and the last equality follows from \eqref{specSums_formula1} and the fact that $X_{0,1}=0$.
Similarly, using \eqref{help1}, \eqref{specSums_formula1} and the fact that $X_{0,3}=0$, we can write
  \begin{align*}
  X_{k,3}
    &=\sum_{j=1}^{X_{k-1,1}}\xi_{k,j,1,3}+\sum_{j=1}^{X_{k-1,2}}\xi_{k,j,2,3}+\sum_{j=1}^{X_{k-1,3}}\xi_{k,j,3,3}+\vare_{k,3}\\
    &=\sum_{j=1}^{X_{k-1,1}}\xi_{k,j,1,3}+\sum_{j=1}^{X_{k-1,2}}\xi_{k,j,2,3}+\sum_{j=1}^{X_{k-1,3}}\xi_{k,j,3,3}+\vare_{k,3}\\
    &\phantom{=\;} + X_{k-1,1}(a_{3,1}-a_{3,1}) +  X_{k-1,2}(a_{3,2}-a_{3,2}) +  X_{k-1,3}(1-1)\\
    &=X_{k-1,3} + a_{3,2}X_{k-1,2} + a_{3,1}X_{k-1,1} \\
    &\phantom{=\;} + \sum_{j=1}^{X_{k-1,1}}(\xi_{k,j,1,3}-a_{3,1})+\sum_{j=1}^{X_{k-1,2}}(\xi_{k,j,2,3}-a_{3,2})
       + \sum_{j=1}^{X_{k-1,3}}(\xi_{k,j,3,3}-1) +\vare_{k,3}\\
    &= a_{3,2} \sum_{\ell=1}^{k-1} X_{\ell,2} + E_{k,3}
     = a_{3,2} \int_0^{\frac{k}{n}} n X_{\ns,2}\,\dd s + E_{k,3},   \qquad  k,n\in\NN,
  \end{align*}
 where the error term $E_{k,3}$ is defined by
 \begin{align*}
  E_{k,3}&:=  a_{3,1} \sum_{\ell=1}^{k-1} X_{\ell,1}
             + \sum_{\ell=1}^k \sum_{j=1}^{X_{\ell-1,1}}(\xi_{\ell,j,1,3}-a_{3,1})
             + \sum_{\ell=1}^k \sum_{j=1}^{X_{\ell-1,2}}(\xi_{\ell,j,2,3}-a_{3,2}) \\
          &\phantom{:=\,}
             + \sum_{\ell=1}^k \sum_{j=1}^{X_{\ell-1,3}}(\xi_{\ell,j,3,3}-1)
             + \sum_{\ell=1}^k \vare_{\ell,1}.
 \end{align*}
All in all, writing $\nt$ instead of $k$ as well, we can obtain that
 \begin{align*}
 \begin{bmatrix}
 n^{-1}X_{\nt,1}\\
 n^{-2}X_{\nt,2}\\
 n^{-3}X_{\nt,3}
 \end{bmatrix}
 &=
 \begin{bmatrix}
 n^{-1}X_{\nt,1}\\
 a_{2,1}\int_0^{\frac{\nt}{n}}n^{-1}X_{\ns,1}\,\dd s\\
 a_{3,2}\int_0^{\frac{\nt}{n}}n^{-2}X_{\nr,2}\,\dd r
 \end{bmatrix}
 + \begin{bmatrix}
 0\\
 n^{-2} E_{\nt,2}\\
 n^{-3} E_{\nt,3}
 \end{bmatrix}\\
 &=\begin{bmatrix}
 n^{-1}X_{\nt,1}\\
 a_{2,1}\int_0^{\frac{\nt}{n}}n^{-1}X_{\ns,1}\,\dd s\\
 a_{3,2}a_{2,1}\int_0^{\frac{\nt}{n}}\int_0^{\frac{\nr}{n}}n^{-1}X_{\ns,1}\,\dd s\,\dd r
 \end{bmatrix}
 +\begin{bmatrix}
  0\\
  n^{-2} E_{\nt,2}\\
  a_{3,2}\int_0^{\frac{\nt}{n}}n^{-2}E_{\nr,2}\,\dd r+ n^{-3}E_{\nt,3}
 \end{bmatrix}\\
& =: \widetilde\bX_{\nt} + \bE_\nt, \qquad n\in\NN, \; t\in\RR_+.
 \end{align*}
By Theorem \ref{thm:critical}, Lemmas \ref{ConvMatrixMult} and \ref{ConvIteratedIntegral},
 and the fact that the limit process in Theorem \ref{thm:critical} has continuous sample paths almost surely,
 we have $(\widetilde\bX_\nt)_{t\in\RR_+}\distr(\bcX_t)_{t\in\RR_+}$ as $n\to\infty$,
 where $(\bcX_t)_{t\in\RR_+}$ is the limit process in Theorem \ref{main_4}.
Therefore, if the error process $(\bE_\nt)_{t\in\RR_+}$ disappeared as $n\to\infty$
in an appropriate sense, then Theorem \ref{main_4} would hold.
\proofend
\end{Rem}

\begin{Rem}
We suspect that the moment conditions in Theorems \ref{main_1} and \ref{main_2} might be relaxed to
 the moment condition \eqref{baseconds} using the method of the proof of Theorem 3.1 in Barczy et al.\ \cite{BarIspPap0}.
In fact, the fourth order moment assumptions in the proofs of Theorems \ref{main_1} and \ref{main_2} are used only for checking the
 conditional Lindeberg condition, namely, condition (iii) of Theorem \ref{Conv2DiffThm}.
For single-type critical Galton-Watson processes with immigration, a detailed exposition of a proof of the conditional Lindeberg condition
 in question under second order moment assumptions can be found, e.g., in Barczy et al.\ \cite{BarBezPap1}.
\proofend
\end{Rem}

\section{Preliminaries for the proofs}\label{Section_prel_proof}

First, we recall some results about the powers of lower triangular matrices with
 all the diagonal entries equal to $1$.

\begin{Lem}\label{matrixlemma}
Let $\bA\in\RR^{p\times p}$
 be a lower triangular matrix such that $a_{i,i}=1$, $i\in\{1,\ldots,p\}$.
Then for each $k\in\NN$, we have
\begin{equation}\label{mPow}
\bA^k=\sum_{m=0}^{p-1}\binom{k}{m}\left(\bA-\bI_p\right)^m.
\end{equation}
\end{Lem}

\noindent\textbf{Proof.}
Let $\bC := \bA-\bI_p$.
Then, by the binomial theorem, we get
 \begin{equation*}
   \bA^k=(\bC+\bI_p)^k=\sum_{m=0}^k\binom{k}{m}\bC^m, \qquad k\in\NN.
 \end{equation*}
Since $\bC$ is lower triangular and $c_{i,i}=0$ for $i\in\{1,\dots,p\}$, we have $\bC^m=\bzero$ whenever
  $m\geq p$, $m\in\NN$, so we get that
 \begin{equation*}
   \bA^k=\sum_{m=0}^{(p-1)\wedge k}\binom{k}{m}\bC^m, \qquad k\in\NN.
 \end{equation*}
This shows \eqref{mPow} taking into account the convention $\binom{k}{m}=0$ for $k,m\in\ZZ_+$ with $k<m$.
\proofend

In the next remark, we collect a few facts about the powers of the matrix $\bC=\bA-\bI_p$ defined in the proof of Lemma \ref{matrixlemma}.

\begin{Rem}\label{matrix_power_remark}
(i). We check that for each $i,j\in\{1,\dots,p\}$, we have
\begin{equation}\label{C_0}
  c_{i,j}^{[m]}=0 \qquad \text{ for $m>i-j$, $m\in\NN$.}
\end{equation}
First, suppose that $i,j\in\{1,\ldots,p\}$ are such that $i\leq j$.
Then \eqref{C_0} holds trivially for $m=1$, and if it holds for $m\in\NN$, then,
 by the rules of matrix multiplication, the nonnegativity of the entries of $\bC$,
 and the fact that $c_{i,r}=0$ for $r\geq i$, $i,r\in\{1,\ldots,p\}$,
 we get that
\begin{equation*}
c_{i,j}^{[m+1]} =\sum_{r=1}^pc_{i,r}c_{r,j}^{[m]}=\sum_{r=1}^{i-1} c_{i,r} c_{r,j}^{[m]}\leq\sum_{r=1}^{j-1} c_{i,r}c_{r,j}^{[m]}
                =\sum_{r=1}^{j-1}c_{i,r}\cdot 0=0,
\end{equation*}
where at the last but one equality, we used the induction hypothesis.
Since $c_{i,j}^{[m+1]}\geq 0$, it implies that $c_{i,j}^{[m+1]}=0$, hence we proved \eqref{C_0} for $i\leq j$, $i,j\in\{1,\dots,p\}$ and $m\in\NN$.

Next, suppose that $i,j\in\{1,\ldots,p\}$ are such that $i>j$.
Then $i-j\in\{1,\ldots,p-1\}$.
If $i-j=1$ and $m>i-j$, $m\in\NN$, then $m-1\in\NN$, and thus
\begin{equation*}
c_{i,j}^{[m]}=c_{j+1,j}^{[m]}=\sum_{r=1}^pc_{j+1,r}c_{r,j}^{[m-1]}
    =\sum_{r=1}^j c_{j+1,r} c_{r,j}^{[m-1]}
    =\sum_{r=1}^j c_{j+1,r}\cdot0=0,
\end{equation*}
 where the last but one equality follows from the previous case.
 Now, suppose that the statement holds whenever $i-j\in\{1,\dots,\ell\}$ with some $\ell\in\{1,\ldots,p-2\}$.
If $i-j=\ell+1$, then for $m\ge \ell+2$, $m\in\NN$, we have
\begin{equation*}
 c_{i,j}^{[m]}=\sum_{r=1}^pc_{i,r}c_{r,j}^{[m-1]}=\sum_{r=1}^{i-1}c_{i,r}c_{r,j}^{[m-1]}=\sum_{r=1}^{i-1}c_{i,r}\cdot0=0,
\end{equation*}
  where the last but one equality can be checked as follows.
If $r\in\{1,\ldots,i-1\}$ is such that $r\leq j$, then, by the previous case, $c_{r,j}^{[m-1]}=0$.
If $r\in\{1,\ldots,i-1\}$ is such that $r>j$, then $ 1\leq r-j\leq i-j-1=\ell< \ell+1\leq m-1$,
 yielding that $r-j\in\{1,\ldots,\ell\}$ and $m-1> r-j$.
Consequently, by our hypothesis, we have $c_{r,j}^{[m-1]}=0$.

(ii). We derive that
\begin{equation}\label{C_prods}
c_{i,j}^{[i-j]}=\prod_{r=j}^{i -1}a_{r+1,r} \qquad \text{for $i>j$, $i,j\in\{1,\dots,p\}$.}
\end{equation}
This can be easily shown by induction on the value of $i-j\in\{1,\ldots,p-1\}$.
If $i-j=1$, then $c_{i,j}^{[i-j]}=c_{j+1,j}=a_{j+1,j}$, so \eqref{C_prods} holds.
Suppose that \eqref{C_prods} holds whenever $i-j\in\{1,\ldots,\ell\}$ with some $\ell\in\{1,\ldots,p-2\}$.
If $i-j=\ell+1$, then
\begin{equation*}
 c_{i,j}^{[i-j]}=\sum_{r=1}^pc_{i,r} c_{r,j}^{[\ell]}
                =\sum_{r=1}^{i-1}a_{i,r} c_{r,j}^{[\ell]}=a_{i,i-1}c_{i-1,j}^{[i-j-1]}
  =a_{i,i-1}\prod_{r=j}^{i-2}a_{r+1,r}
  =\prod_{r=j}^{i-1}a_{r+1,r},
\end{equation*}
 since if $r\in\{1,\ldots,i-2\}$, then $r-j< (i-1)-j=\ell$, and thus, by \eqref{C_0},
  $c_{r,j}^{[\ell]}=0$ whenever $r\in\{1,\ldots,i-2\}$.
\proofend
\end{Rem}

For a lower triangular matrix $\bA\in\RR_+^{p\times p}$ such that $a_{j,j}=1$, $j\in\{1,\dots,p\}$,
let us recall the notation $\bC=\bA-\bI_p$, and for $i\in\{1,\dots,p\}$, introduce
 \begin{equation}\label{etadef}
\eta_i := \eta_i(\bA) :=\max_{m\in\{1,\dots,p\}}\{m \mid \exists\, j\in\{1,\dots,p\}: c_{i,j}^{[m-1]}>0\},
 \end{equation}
  where we recall that $c_{i,j}^{[m-1]}$ denotes the $(i,j)$-th entry of $\bC^{m-1} =(\bA-\bI_p)^{m-1}$.
 Note that due to $c_{i,i}^{[0]}=(\bI_p)_{i,i}=1$, we have $\eta_i\geq1$ for each $i\in\{1,\dots,p\}$.

In the next lemma, we derive an inequality between $\eta_i$ and $\eta_j$, $i,j\in\{1,\dots,p\}$.

\begin{Lem}\label{eta_ineq}
Let $p\in\NN$, and let $\bA\in\RR_+^{p\times p}$ be a lower triangular matrix such that $a_{i,i}=1$, $i\in\{1,\dots,p\}$.
If $c_{i,j}^{[m]}>0$ for some $i,j\in\{1,\dots,p\}$ and $m\in\{0,1,\ldots,p-1\}$,
 then $\eta_i\geq\eta_j+m$.
\end{Lem}

\noindent\textbf{Proof.}
Assume $c_{i,j}^{[m]}>0$ for some $i,j\in\{1,\ldots,p\}$ and $m\in\{0,\dots,p-1\}$.
Since $\eta_j\geq 1$, by definition \eqref{etadef}, there exists
 some $r_0\in\{1,\dots,p\}$ such that $c_{j,r_0}^{[\eta_j-1]}>0$,
 and, by the properties of matrix multiplication and the non-negativity of the coefficients $c_{\ell,k}^{[q]}$, $q\in\ZZ_+$, $\ell,k\in\{1,\ldots,p\}$,
 we get that
 \begin{equation*}
 c_{i,r_0}^{[\eta_j+m-1]}=\sum_{r=1}^pc_{i,r}^{[m]}c_{r,r_0}^{[\eta_j-1]}\geq c_{i,j}^{[m]}c_{j,r_0}^{[\eta_j-1]}>0.
 \end{equation*}
This implies $\eta_i\geq\eta_j+m$.
\proofend

In the rest of this section, we will only consider the case $p=3$.
First, we specialize formula \eqref{mPow} to the case $p=3$.
\begin{Cor}\label{Cor_Ak}
Let $\bA\in\RR^{3\times 3}$ be a lower triangular matrix such that $a_{i,i}=1$, $i\in\{1,2,3\}$.
Then for each $k\in\ZZ_+$, we have
 \[
  \bA^k
  =\begin{bmatrix}
     1 & 0 & 0 \\
     a_{2,1} & 1 & 0 \\
     a_{3,1} & a_{3,2} & 1 \\
   \end{bmatrix}^k
 =  \begin{bmatrix}
     1 & 0 & 0 \\
     k a_{2,1} & 1 & 0 \\
     \binom{k}{2}a_{3,2}a_{2,1} + ka_{3,1}& ka_{3,2} & 1 \\
   \end{bmatrix}.
 \]
\end{Cor}

\noindent{\bf Proof.}
In the considered case, with $\bC=\bA-\bI_3$, we have
\begin{equation}\label{C_hatvanyok}
\bC^0=\begin{bmatrix}
     1 & 0 & 0 \\
     0 & 1 & 0 \\
     0 & 0 & 1 \\
   \end{bmatrix},
 \qquad
\bC^1=   \begin{bmatrix}
     0 & 0 & 0 \\
     a_{2,1} & 0 & 0 \\
     a_{3,1} & a_{3,2} & 0 \\
   \end{bmatrix},\qquad
\bC^2=   \begin{bmatrix}
     0 & 0 & 0 \\
     0 & 0 & 0 \\
     a_{3,2}a_{2,1} & 0 & 0 \\
   \end{bmatrix},
\end{equation}
 and, by \eqref{mPow}, for each $k\in\NN$ we have
\begin{equation*}
\begin{bmatrix}
     1 & 0 & 0 \\
     a_{2,1} & 1 & 0 \\
     a_{3,1} & a_{3,2} & 1 \\
   \end{bmatrix}^k=\binom{k}{0}\bC^0+\binom{k}{1}\bC^1+\binom{k}{2}\bC^2,
\end{equation*}
 yielding the assertion.
\proofend

The quantities $\eta_i$, $i\in\{1,2,3\}$, defined in \eqref{etadef}, will become relevant in the moment estimations
 presented in Lemmas \ref{EEX1} and \ref{M_sum_growth}, which play crucial role
 in the proofs of Theorems \ref{main_1} and \ref{main_2}--\ref{main_4}.
Next, we calculate their exact values in the cases \textup{(1)}-\textup{(4)} of \eqref{tablazat_esetek}.

\begin{Lem}\label{eta_values}
Let $(\bX_k)_{k\in\ZZ_+}$ be a strongly critical 3-type GWI process such that $\bX_0=\bzero$, the moment condition \eqref{baseconds} holds,
 and suppose that the offspring mean matrix $\bA$ is lower triangular such that $a_{i,i}=1$, $i\in\{1,2,3\}$.
Then, in the four cases \textup{(1)}-\textup{(4)} of \eqref{tablazat_esetek},
 the quantities $\eta_i$, $i\in\{1,2,3\}$, defined in \eqref{etadef}, take the following values:
\begin{align*}%\label{eta_value_table}
 \begin{tabular}{|c|c|c|c|}
  \hline
  \textup{(1)} & $\eta_1 = 1$ & $\eta_2 = 1$ & $\eta_3 = 1$ \\
  \hline
  \textup{(2)} & $\eta_1 = 1$ & $\eta_2 = 1$ & $\eta_3 = 2$ \\
  \hline
  \textup{(3)} & $\eta_1 = 1$ & $\eta_2 = 2$ & $\eta_3 = 2$ \\
  \hline
  \textup{(4)} & $\eta_1 = 1$ & $\eta_2 = 2$ & $\eta_3 = 3$ \\
  \hline
 \end{tabular}\;.
\end{align*}
\end{Lem}

\noindent\textbf{Proof.}
By \eqref{C_hatvanyok}, it is clear that in all the four cases, $c_{1,j}^{[1]}=c_{1,j}^{[2]}=0$ for each $j\in\{1,2,3\}$,
 so $\eta_1\leq 1$, which together with $\eta_1\geq1$ implies $\eta_1=1$.

Similarly, we have that, in all the four cases, $c_{2,j}^{[2]}=0$ for each $j\in\{1,2,3\}$, therefore $\eta_2\leq 2$.
Since $c_{2,2}^{[1]}=c_{2,3}^{[1]}=0$, the value of $\eta_2$ depends on whether $c_{2,1}^{[1]}=a_{2,1}$ is positive or not.
In cases \textup{(1)}-\textup{(2)} we have $a_{2,1}=0$, and so $\eta_2=1$,
while in cases \textup{(3)}-\textup{(4)} we have $a_{2,1}>0$, thus $\eta_2=2$.

Similarly, we have that $c_{3,2}^{[2]}=c_{3,3}^{[2]}=0$ in all four cases,
 and $c_{3,1}^{[2]}=a_{3,2}a_{2,1}>0$ holds if and only if $a_{3,2}>0$ and $a_{2,1}>0$,
 which corresponds to case \textup{(4)}, thus in this case we have $\eta_3=3$, and otherwise $\eta_3\leq 2$.
In cases \textup{(2)}-\textup{(3)} we have $c_{3,1}^{[1]}=a_{3,1}>0$, thus $\eta_3\geq2$, which together with $\eta_3\leq 2$ implies $\eta_3=2$.
The only remaining case is $\textup{(1)}$, in which we have $c_{3,j}^{[1]}=0$ for each $j\in\{1,2,3\}$,
 implying $\eta_3\leq 1$, which together with $\eta_3\geq1$ yields that $\eta_3=1$.
\proofend

The explicit form of $\bA^k$, $k\in\NN$, in Corollary \ref{Cor_Ak} together with formula \eqref{mean2} enable us to describe
 the asymptotic behavior of $\EE(\bX_k)$ as $k\to\infty$, see the next proposition.

\begin{Pro}\label{Pro_Xk_exp}
Let $(\bX_k)_{k\in\ZZ_+}$ be a 3-type GWI process such that $\bX_0=\bzero$,
 and the moment condition \eqref{baseconds} holds.
If the offspring mean matrix $\bA$ of $(\bX_k)_{k\in\ZZ_+}$ satisfies \textup{(1)} of \eqref{tablazat_esetek}, then we have
 \begin{align}\label{help_Xk_exp_1}
   \lim_{k\to\infty} \begin{bmatrix}
                       k^{-1} \EE(X_{k,1}) \\
                       k^{-1} \EE(X_{k,2}) \\
                       k^{-1} \EE(X_{k,3}) \\
                     \end{bmatrix}
                     = \begin{bmatrix}
                       b_1 \\
                       b_2 \\
                       b_3 \\
                     \end{bmatrix}.
 \end{align}
If the offspring mean matrix $\bA$ of $(\bX_k)_{k\in\ZZ_+}$ satisfies \textup{(2)} of \eqref{tablazat_esetek}, then we have
 \begin{align}\label{help_Xk_exp_2}
 \lim_{k\to\infty} \begin{bmatrix}
                       k^{-1} \EE(X_{k,1}) \\
                       k^{-1} \EE(X_{k,2}) \\
                       k^{-2} \EE(X_{k,3}) \\
                     \end{bmatrix}
                     = \begin{bmatrix}
                       b_1 \\
                       b_2 \\
                       \frac{1}{2}(b_1a_{3,1} + b_2a_{3,2}) \\
                     \end{bmatrix}.
 \end{align}
If the offspring mean matrix $\bA$ of $(\bX_k)_{k\in\ZZ_+}$ satisfies \textup{(3)} of \eqref{tablazat_esetek}, then we have
 \begin{align}\label{help_Xk_exp_3}
 \lim_{k\to\infty} \begin{bmatrix}
                       k^{-1} \EE(X_{k,1}) \\
                       k^{-2} \EE(X_{k,2}) \\
                       k^{-2} \EE(X_{k,3}) \\
                     \end{bmatrix}
                     = \begin{bmatrix}
                       b_1 \\
                       \frac{1}{2}b_1a_{2,1} \\
                       \frac{1}{2}b_1a_{3,1} \\
                     \end{bmatrix}.
 \end{align}
If the offspring mean matrix $\bA$ of $(\bX_k)_{k\in\ZZ_+}$ satisfies \textup{(4)} of \eqref{tablazat_esetek}, then we have
 \begin{align}\label{help_Xk_exp_4}
 \lim_{k\to\infty} \begin{bmatrix}
                       k^{-1} \EE(X_{k,1}) \\
                       k^{-2} \EE(X_{k,2}) \\
                       k^{-3} \EE(X_{k,3}) \\
                     \end{bmatrix}
                     = \begin{bmatrix}
                       b_1 \\
                       \frac{1}{2}b_1a_{2,1} \\
                       \frac{1}{6}b_1a_{3,2}a_{2,1} \\
                     \end{bmatrix}.
 \end{align}
\end{Pro}

\noindent{\bf Proof.}
By \eqref{mean2}, we have $\EE(\bX_k)= \sum_{j=0}^{k-1} \bA^j \bb$, $k\in\NN$.
Using Corollary \ref{Cor_Ak}, it implies that
 \begin{align*}
  \EE(\bX_k)
   & = \begin{bmatrix}
      k & 0 & 0 \\
      a_{2,1}\sum_{j=1}^{k-1} j & k & 0 \\
      a_{3,2}a_{2,1} \sum_{j=1}^{k-1} \binom{j}{2} + a_{3,1}\sum_{j=1}^{k-1} j & a_{3,2}\sum_{j=1}^{k-1} j & k \\
    \end{bmatrix}
     \bb \\
  & = \begin{bmatrix}
      k & 0 & 0 \\
      a_{2,1}\frac{(k-1)k}{2} & k & 0 \\
      \frac{a_{3,2}a_{2,1}}{2} \sum_{j=1}^{k-1} (j^2-j) + a_{3,1}\frac{(k-1)k}{2} & a_{3,2}\frac{(k-1)k}{2} & k \\
    \end{bmatrix}
    \bb, \qquad k\in\NN,
 \end{align*}
 where
 \[
   \sum_{j=1}^{k-1} (j^2-j) = \frac{(k-1)k(2(k-1)+1)}{6} - \frac{(k-1)k}{2}
                            = \frac{(k-2)(k-1)k}{3}.
 \]
This yields that
 \begin{align}\label{help_Xk_exp_5}
   \EE(\bX_k)
      =  \begin{bmatrix}
            b_1k \\
            b_1a_{2,1}\binom{k}{2} + b_2k \\
            b_1a_{3,2}a_{2,1}\binom{k}{3} + (b_1a_{3,1} + b_2a_{3,2})\binom{k}{2} + b_3k \\
          \end{bmatrix},
          \qquad k\in\NN.
 \end{align}
Hence, in case of \textup{(1)} of \eqref{tablazat_esetek}, we have
 \[
   \EE(\bX_k)
      =  \begin{bmatrix}
            b_1k \\
            b_2k \\
            b_3k \\
          \end{bmatrix},
          \qquad k\in\NN,
 \]
 yielding \eqref{help_Xk_exp_1}.
In case of \textup{(2)} of \eqref{tablazat_esetek}, we have
 \[
   \EE(\bX_k)
      =  \begin{bmatrix}
            b_1k \\
            b_2k \\
            \frac{1}{2}(b_1a_{3,1} + b_2a_{3,2})k^2 + \left(b_3-\frac{b_1a_{3,1}+b_2a_{3,2}}{2}\right)k \\
          \end{bmatrix},
          \qquad k\in\NN,
 \]
 yielding \eqref{help_Xk_exp_2}.
In case of \textup{(3)} of \eqref{tablazat_esetek}, we have
 \[
   \EE(\bX_k)
      =  \begin{bmatrix}
            b_1k \\
           \frac{1}{2}b_1a_{2,1}k^2 + \left(b_2-\frac{b_1a_{2,1}}{2}\right)k\\
           \frac{1}{2}b_1a_{3,1}k^2+\left(b_3-\frac{b_1a_{3,1}}{2}\right)k
          \end{bmatrix},
          \qquad k\in\NN,
 \]
 yielding \eqref{help_Xk_exp_3}.
In case of \textup{(4)} of \eqref{tablazat_esetek}, using the formula \eqref{help_Xk_exp_5}, we have
 \[
   \EE(\bX_k)
      =  \begin{bmatrix}
            b_1k \\
           \frac{1}{2}b_1a_{2,1}k^2 + \OO(k)\\
           \frac{1}{6}b_1a_{3,2}a_{2,1}k^3+\OO(k^2)
          \end{bmatrix},
          \qquad k\in\NN,
 \]
yielding \eqref{help_Xk_exp_4}.
\proofend

We note that Proposition \ref{Pro_Xk_exp} is in fact a generalization of Theorem 3 in Foster and Ney \cite{FosNey}
 in case of $p=3$, for more details, see the last paragraph of Remark \ref{FNmoment}.
Remark also that the normalization of the coordinates of $\EE(\bX_k)$ in Proposition \ref{Pro_Xk_exp} in the cases
 \textup{(1)}--\textup{(4)} of \eqref{tablazat_esetek} are exactly the same ones as the normalizations
 for the coordinates of $\bX_{\nt}$ in Theorems \ref{main_1} and \ref{main_2}--\ref{main_4} corresponding to the cases
 \textup{(1)}--\textup{(4)} of \eqref{tablazat_esetek}.
Concerning the order of the moments of a strongly critical $p$-type GWI process starting from $\bzero$
 and with lower triangular offspring mean matrix having diagonal entries $1$, see Lemma \ref{EEX1} as well.

Next, we derive a decomposition for the $3$-type GWI process $(\bX_k)_{k\in\ZZ_+}$
starting from $\bzero$, satisfying the moment conditions \eqref{baseconds},
and having a lower triangular offspring mean matrix $\bA$ with diagonal entries $1$.
Let us introduce the sequence
 \begin{equation}\label{Mk1}
  \begin{split}
  \begin{bmatrix} M_{k,1} \\ M_{k,2} \\ M_{k,3} \end{bmatrix}
  &:= \bM_k := \bX_k - \EE(\bX_k \mid \cF_{k-1}^{\bX})
   = \bX_k - \bA \bX_{k-1} - \bb \\
  &= \begin{bmatrix}
      X_{k,1} - X_{k-1,1} - b_1 \\
      X_{k,2} - a_{2,1} X_{k-1,1} - X_{k-1,2} - b_2 \\
      X_{k,3} - a_{3,1} X_{k-1,1} - a_{3,2} X_{k-1,2} - X_{k-1,3} - b_3
     \end{bmatrix} , \qquad
   k \in \NN ,
  \end{split}
 \end{equation}
 of martingale differences with respect to the filtration
 \ $(\cF_k^\bX)_{k\in\ZZ_+}$, \ where we used \eqref{mart}.
From \eqref{Mk1}, we obtain the recursion
 $ \bX_k = \bA \bX_{k-1} + \bM_k + \bb$, $k \in \NN$,
 which together with \ $\bX_0=\bzero$ \ implies
 \begin{equation}\label{X}
   \bX_k = \sum_{\ell=1}^k \bA^{k-\ell} (\bM_\ell + \bb) , \qquad
   k \in \NN.
 \end{equation}

Finally, by \eqref{X} and Corollary \ref{Cor_Ak}, we can get a decomposition
 \begin{equation}\label{Xdeco}
  \begin{bmatrix}
   X_{k,1} \\
   X_{k,2} \\
   X_{k,3}
  \end{bmatrix}
  = \begin{bmatrix}
     X_{k,1}^{(1)} \\
     a_{2,1} X_{k,2}^{(1)} + X_{k,2}^{(2)} \\
     a_{3,2}a_{2,1} X_{k,3}^{(1)} + a_{3,1} X_{k,3}^{(2)} + a_{3,2} X_{k,3}^{(3)} + X_{k,3}^{(4)}
    \end{bmatrix} , \qquad k \in \NN ,
 \end{equation}
 where the stochastic processes $(X_{k,1}^{(1)})_{k\in\NN}$,
 $(X_{k,2}^{(1)})_{k\in\NN}$, $(X_{k,2}^{(2)})_{k\in\NN}$, $(X_{k,3}^{(1)})_{k\in\NN}$, $(X_{k,3}^{(2)})_{k\in\NN}$, $(X_{k,3}^{(3)})_{k\in\NN}$,
 and $(X_{k,3}^{(4)})_{k\in\NN}$ are given by
 \begin{align*}
    & X_{k,1}^{(1)} := \sum_{\ell=1}^k (M_{\ell,1} + b_1),\\[1mm]
    & X_{k,2}^{(1)} := \sum_{\ell=1}^k (k-\ell) (M_{\ell,1} + b_1),
     \qquad
     X_{k,2}^{(2)} := \sum_{\ell=1}^k (M_{\ell,2} + b_2),
 \end{align*}
 and
 \begin{align*}
    & X_{k,3}^{(1)} := \frac{1}{2} \sum_{\ell=1}^k (k-\ell)(k-\ell-1) (M_{\ell,1} + b_1) ,\qquad
      X_{k,3}^{(2)} := \sum_{\ell=1}^k (k-\ell) (M_{\ell,1} + b_1) = X_{k,2}^{(1)},\\[1mm]
    & X_{k,3}^{(3)} :=  \sum_{\ell=1}^k (k-\ell)(M_{\ell,2} + b_2),\qquad
      X_{k,3}^{(4)} := \sum_{\ell=1}^k (M_{\ell,3} + b_3).
 \end{align*}
Indeed, for each $k\in\NN$, we have
 \begin{align*}
   &X_{k,1} = \sum_{\ell=1}^k (\bA^{k-\ell})_{1,1} (M_{\ell,1} + b_1)
            = \sum_{\ell=1}^k (M_{\ell,1} + b_1), \\
   &X_{k,2} = \sum_{\ell=1}^k (\bA^{k-\ell})_{2,1} (M_{\ell,1} + b_1)
              + \sum_{\ell=1}^k (\bA^{k-\ell})_{2,2} (M_{\ell,2} + b_2)\\
   &\phantom{X_{k,2}}
            = a_{2,1} \sum_{\ell=1}^k (k-\ell) (M_{\ell,1} + b_1)
              + \sum_{\ell=1}^k (M_{\ell,2} + b_2) \\
  &\phantom{X_{k,2}}
            = a_{2,1} X_{k,2}^{(1)} + X_{k,2}^{(2)},
 \end{align*}
 and
 \begin{align*}
 &X_{k,3} = \sum_{\ell=1}^k (\bA^{k-\ell})_{3,1} (M_{\ell,1} + b_1)
              + \sum_{\ell=1}^k (\bA^{k-\ell})_{3,2} (M_{\ell,2} + b_2)
              + \sum_{\ell=1}^k (\bA^{k-\ell})_{3,3}  (M_{\ell,3} + b_3)\\
 &\phantom{X_{k,3}}
   = \sum_{\ell=1}^k \left( \binom{k-\ell}{2}a_{3,2}a_{2,1} + (k-\ell)a_{3,1} \right) (M_{\ell,1} + b_1)  \\
 &\phantom{X_{k,3}=\;}
          + \sum_{\ell=1}^k (k-\ell) a_{3,2} (M_{\ell,2} + b_2)
          + \sum_{\ell=1}^k (M_{\ell,3} + b_3)\\
 &\phantom{X_{k,3}}
  = a_{3,2}a_{2,1} \sum_{\ell=1}^k \binom{k-\ell}{2} (M_{\ell,1} + b_1)
     + a_{3,1} \sum_{\ell=1}^k (k-\ell)(M_{\ell,1} + b_1)  \\
 &\phantom{X_{k,3}=\;}
     + a_{3,2} \sum_{\ell=1}^k (k-\ell)(M_{\ell,2} + b_2)
     + \sum_{\ell=1}^k (M_{\ell,3} + b_3),
 \end{align*}
and thus we have that
\begin{align*}
X_{k,3} = a_{3,2}a_{2,1} X_{k,3}^{(1)} + a_{3,1} X_{k,2}^{(1)} + a_{3,2}X_{k,3}^{(3)} + X_{k,3}^{(4)}, \qquad k\in\NN.
\end{align*}

Using that $X_{j,1}^{(1)}=\sum_{\ell=1}^j (M_{\ell,1} + b_1)$
 and $X_{j,2}^{(2)}=\sum_{\ell=1}^j (M_{\ell,2} + b_2)$, $j\in\NN$,
 by Lemma \ref{specSums}, for each $k,n\in\NN$, we have that
\begin{align}\label{integrals}
\begin{split}
X_{k,2}^{(1)}=&X_{k,3}^{(2)}
	= \sum_{\ell=1}^k (k - \ell) (M_{\ell,1} + b_1)
	=n\int_0^{\frac{k}{n}}X_{\ns,1}^{(1)}\,\dd s,\\
X_{k,3}^{(3)}
	=& \sum_{\ell=1}^k (k - \ell) (M_{\ell,2} + b_2)
	=n\int_0^{\frac{k}{n}}X_{\ns,2}^{(2)}\,\dd s,\\
X_{k,3}^{(1)}
	=&\sum_{\ell=1}^k \binom{k-\ell}{2} (M_{\ell,1} + b_1)
	=n^2\int_0^{\frac{k}{n}}\left(\int_0^{\frac{\nr}{n}}X_{\ns,1}^{(1)}\,\dd s\right)\,\dd r.
\end{split}
\end{align}

\section{Proof of Theorem \ref{main_1}}\label{Proof1}%

The proof is analogous to the proof of Theorem 2.1 in Barczy et al.\ \cite{BarBezPap2},
 but, for completeness, we give a detailed proof. We divide the proof into several steps.

{\sl Step 1/(a).}
The SDE \eqref{SDE_X_1} has a pathwise unique strong solution $(\bcX_t:=[\cX_{t,1},\cX_{t,2},\cX_{t,3}]^\top )_{t\in\RR_+}$ for all initial values
 $\bcX_0 = \bx \in \RR^3$, and if $\bx \in \RR_+^3$, then $\bcX_t \in \RR_+^3$ almost surely for all $t \in \RR_+$,
  since $b_i$, $v^{(i)}_{i,i}\in \RR_+$, $i=1,2,3$, see, e.g., Ikeda and Watanabe \cite[Chapter IV, Example 8.2]{IkeWat}.
By \eqref{help_Wei_Winnicki1}, we have $(n^{-1} X_{\nt,1})_{t\in\RR_+} \distr (\cX_{t,1})_{t\in\RR_+}$ as $n \to \infty$,
 where $(\cX_{t,1})_{t\in\RR_+}$ satisfies the first equation of the SDE \eqref{SDE_X_1} with initial value $\cX_{0,1}=0$.
Similarly, since $a_{2,1}=a_{3,1}=a_{3,2}=0$,
 as it was explained in Section \ref{Section_conv_results}, the second and third coordinate processes $(X_{k,2})_{k\in\ZZ_+}$
 and $(X_{k,3})_{k\in\ZZ_+}$ are critical single-type GWI processes,
 so $(n^{-1} X_{\nt,2})_{t\in\RR_+} \distr (\cX_{t,2})_{t\in\RR_+}$ as $n \to \infty$,
 and $(n^{-1} X_{\nt,3})_{t\in\RR_+} \distr (\cX_{t,3})_{t\in\RR_+}$ as $n \to \infty$,
 where $(\cX_{t,2})_{t\in\RR_+}$ and $(\cX_{t,3})_{t\in\RR_+}$  satisfy the second and third equations of the SDE \eqref{SDE_X_1}
 with initial value $\cX_{0,2}=0$ and $\cX_{0,3}=0$, respectively.
However, we need to prove joint convergence of $(n^{-1} X_{\nt,1})_{t\in\RR_+}$,
 $(n^{-1} X_{\nt,2})_{t\in\RR_+}$ and $(n^{-1} X_{\nt,3})_{t\in\RR_+}$ as $n\to\infty$.

{\sl Step 1/(b).}
Using $a_{2,1}=a_{3,1}=a_{3,2}=0$ and \eqref{Mk1},
 we obtain that the sequence $(\bM_k)_{k\in\NN}$ of martingale differences with respect
 to the filtration $(\cF_k^\bX)_{k\in\ZZ_+}$ takes the form
 \[
   \bM_k = \bX_k - \bX_{k-1} - \bb , \qquad k \in \NN .
 \]
Consider the random step processes
 \begin{align}\label{help15}
   \bcM_t^{(n)} := \begin{bmatrix} \cM_{t,1}^{(n)} \\ \cM_{t,2}^{(n)} \\ \cM_{t,3}^{(n)}  \end{bmatrix}
   := \frac{1}{n} \sum_{k=1}^\nt \bM_k
   = \frac{1}{n} \bX_\nt - \frac{\nt}{n} \bb , \qquad t \in \RR_+ , \qquad n \in \NN ,
 \end{align}
 where at the third equality we used that \ $\bX_0 = \bzero$.
We show that
 \begin{equation}\label{conv_bM}
  (\bcM_t^{(n)})_{t\in\RR_+} \distr (\bcM_t)_{t\in\RR_+} \qquad
  \text{as \ $n \to \infty$,}
 \end{equation}
 where the limit process $\bcM_t = ([\cM_{t,1}, \cM_{t,2}, \cM_{t,3}]^\top)_{t\in\RR_+}$ is the pathwise unique strong solution of the SDE
 \begin{equation}\label{SDE_bM}
  \begin{cases}
   \dd \cM_{t,1}
   = \sqrt{v^{(1)}_{1,1} (\cM_{t,1} + b_1 t)^+ } \, \dd \cW_{t,1} , \qquad t\in\RR_+,\\[1mm]
   \dd \cM_{t,2}
   = \sqrt{v^{(2)}_{2,2} (\cM_{t,2} + b_2 t)^+ } \, \dd \cW_{t,2} , \qquad t\in\RR_+,\\[1mm]
   \dd \cM_{t,3}
   = \sqrt{v^{(3)}_{3,3} (\cM_{t,3} + b_3 t)^+ } \, \dd \cW_{t,3} , \qquad t\in\RR_+
  \end{cases}
 \end{equation}
 with initial value $ \bcM_0 = \bzero$, where $(\cW_{t,1})_{t\in\RR_+}$, $(\cW_{t,2})_{t\in\RR_+}$ and $(\cW_{t,3})_{t\in\RR_+}$
 are independent standard Wiener processes.

{ \sl Step 1/(c).}
We check that the SDE \eqref{SDE_bM} has a pathwise unique strong solution $(\bcM_t)_{t\in\RR_+}$
 for all initial values $\bcM_0 = \bx \in \RR^3$.
Observe that if $(\bcM_t)_{t\in\RR_+}$ is a strong solution of the SDE \eqref{SDE_bM} with initial value
 $ \bcM_0 = \bx \in \RR^3$, then, by It\^o's formula, the process $[\cP_{t,1}, \cP_{t,2} , \cP_{t,3}]^\top := \bcM_t + \bb t$, $t \in \RR_+$,
 is a pathwise unique strong solution of the SDE
 \begin{equation}\label{SDE_bP}
  \begin{cases}
   \dd \cP_{t,1}
   = b_1 \, \dd t
     + \sqrt{v^{(1)}_{1,1} \, \cP_{t,1}^+} \, \dd \cW_{t,1} , \qquad t\in\RR_+, \\[1mm]
   \dd \cP_{t,2}
   = b_2 \, \dd t
     + \sqrt{v^{(2)}_{2,2} \, \cP_{t,2}^+} \, \dd \cW_{t,2}, \qquad t\in\RR_+,\\[1mm]
  \dd \cP_{t,3}
   = b_3 \, \dd t
     + \sqrt{v^{(3)}_{3,3} \, \cP_{t,3}^+} \, \dd \cW_{t,3}, \qquad t\in\RR_+,
  \end{cases}
 \end{equation}
 with initial value $[\cP_{0,1}, \cP_{0,2}, \cP_{0,3}]^\top = \bx$.
Conversely, if $[\cP_{t,1}, \cP_{t,2}, \cP_{t,3}]^\top$, $t \in \RR_+$, is a strong solution
 of the SDE \eqref{SDE_bP} with initial value $[\cP_{0,1}, \cP_{0,2}, \cP_{0,3}]^\top = \bp \in \RR^3$,
 then, by It\^o's formula, the process $\bcM_t := [\cP_{t,1}, \cP_{t,2}, \cP_{t,3}]^\top - \bb t$, $t \in \RR_+$,
 is a strong solution of the SDE \eqref{SDE_bM} with initial value $\bcM_0 = \bp$.
The equations in \eqref{SDE_bP} are the same as in \eqref{SDE_X_1},
 thus, as it was explained in Step 1/(a), the SDE \eqref{SDE_bP} admits a unique strong solution with arbitrary initial value in $\RR^3$.
Consequently, the SDE \eqref{SDE_bM} (with initial value $\bzero\in\RR^3$) admits a unique strong solution $(\bcM_t)_{t\in\RR_+}$,
 and $(\bcM_t + \bb t)_{t\in\RR_+} \distre (\bcX_t)_{t\in\RR_+}$.

{ \sl Step 2/(a).}
In order to prove \eqref{conv_bM}, we want to apply Theorem \ref{Conv2DiffThm} with the following choices:
$d = r = 3$, $\bcU = \bcM$,
 $\bU_k^{(n)} = n^{-1} \bM_k$, $n,k\in\NN$, $\bU_0^{(n)} =\bzero$, $n\in\NN$,
 $\cF_k^{(n)} = \cF_k^\bX$, $n \in \NN$, $k \in \ZZ_+$ (yielding $\bcU^{(n)} = \bcM^{(n)}$, $n\in\NN$),
 and $\bbeta : \RR_+ \times \RR^3 \to \RR^3$ and $\bgamma : \RR_+ \times \RR^3 \to \RR^{3\times3}$
given by $\bbeta(t, \bx) = \bzero$ and
 \[
   \bgamma(t, \bx)
   = \begin{bmatrix}
      \sqrt{v^{(1)}_{1,1}(x_1 + b_1 t)^+} & 0 & 0 \\
      0 & \sqrt{v^{(2)}_{2,2}(x_2 + b_2 t)^+} & 0 \\
      0 & 0 & \sqrt{v^{(3)}_{3,3}(x_3 + b_3 t)^+} \\
     \end{bmatrix}
 \]
 for $t \in \RR_+$ and $\bx = [x_1, x_2,x_3]^\top \in \RR^3$. With these notations, the SDE \eqref{SDE_bM} takes the form
\begin{equation*}
\dd\bcM_t=\bbeta(t,\bcM_t)+\bgamma(t,\bcM_t)\dd\bcW_t,\qquad t\in\RR_+,
\end{equation*}
with initial value $\bcM_0=\bzero\in\RR^3$.

The convergence $\bU_0^{(n)}\distr \bzero$ as $n\to\infty$, and condition (i) of Theorem \ref{Conv2DiffThm} trivially holds
 (since $\EE(\bM_k \mid \cF_{k-1}^\bX)=\bzero$, $k\in\NN$, and $\bbeta(t, \bx) = \bzero$, $t\in\RR_+$, $\bx\in\RR^3$).
Now, we show that conditions (ii) and (iii) of Theorem \ref{Conv2DiffThm} hold.
We have to check that for all $T \in \RR_{++}$,
 \begin{gather} \label{Condb1}
  \sup_{t\in[0,T]}
   \biggl\Vert\frac{1}{n^2}
          \sum_{k=1}^\nt
           \EE(\bM_k \bM_k^\top \mid \cF_{k-1}^\bX)
          - \int_0^t \bcR^{(n)}_s \,\bV_\bxi \, \dd s \biggr\Vert
  \stoch 0 \qquad \text{as \ $n\to\infty$,} \\
  \frac{1}{n^2}
  \sum_{k=1}^\nT
   \EE(\|\bM_k\|^2 \bbone_{\{\|\bM_k\|>n\theta\}} \mid \cF_{k-1}^\bX)
  \stoch 0  \qquad \text{as \ $n\to\infty$ \ for all \ $\theta \in \RR_{++}$,} \label{Condb2}
 \end{gather}
 where the process $(\bcR^{(n)}_s)_{s \in \RR_+}$ and the matrix $\bV_\bxi$ are defined by
 \begin{gather*}
  \bcR^{(n)}_s
  := \begin{bmatrix}
      (\cM_{s,1}^{(n)} + b_1 s)^+ & 0 & 0 \\
      0 & (\cM_{s,2}^{(n)} + b_2 s)^+ & 0 \\
      0 & 0 & (\cM_{s,3}^{(n)} + b_3 s)^+ \\
     \end{bmatrix} , \qquad s \in \RR_+ , \qquad n \in \NN , \\
  \bV_\bxi
  := \begin{bmatrix}
      v^{(1)}_{1,1} & 0  & 0 \\
      0 & v^{(2)}_{2,2} & 0 \\
      0 & 0 & v^{(3)}_{3,3} \\
     \end{bmatrix} .
 \end{gather*}
Indeed, $\EE( \bM_k\mid \cF_{k-1}^{\bX}) = \bzero$ and thus $\var(\bM_k\mid\cF_{k-1}^{\bX})=\EE(\bM_k\bM_k^\top\mid\cF_{k-1}^{\bX})$, $k\in\NN$,
and, since $\gamma(t,\bx)$ is symmetric for all $t\in\RR_+$ and $\bx\in\RR^3$, we have
 \begin{align*}
  \gamma(t,\bx) \gamma(t,\bx)^\top
  & =  \begin{bmatrix}
      \sqrt{v^{(1)}_{1,1}(x_1 + b_1 t)^+} & 0 & 0 \\
      0 & \sqrt{v^{(2)}_{2,2}(x_2 + b_2 t)^+} & 0 \\
      0 & 0 & \sqrt{v^{(3)}_{3,3}(x_3 + b_3 t)^+} \\
     \end{bmatrix}^2\\
  & = \begin{bmatrix}
       v^{(1)}_{1,1} (x_1 + b_1 t)^+ & 0 & 0 \\
        0 & v^{(2)}_{2,2}(x_2 + b_2 t)^+ & 0 \\
        0 & 0 &  v^{(3)}_{3,3}(x_3 + b_3 t)^+
       \end{bmatrix} \\
  & =  \begin{bmatrix}
        (x_1 + b_1 t)^+ & 0 & 0 \\
        0 & (x_2 + b_2 t)^+ & 0 \\
        0 & 0 & (x_3 + b_3 t)^+
       \end{bmatrix}
       \bV_\bxi, \qquad t\in\RR_+, \quad \bx\in\RR^3.
 \end{align*}

{ \sl Step 2/(b).}
 Now we check \eqref{Condb1}.
For all $s \in \RR_+$ and $n \in \NN$, we have
 \[
   \bcM_s^{(n)} + \bb s
   = \frac{1}{n} \bX_\ns + \frac{ns-\ns}{n} \bb ,
 \]
 thus
 \[
   \bcR^{(n)}_s
   = \begin{bmatrix}
      \cM_{s,1}^{(n)} + b_1 s & 0 & 0 \\
      0 & \cM_{s,2}^{(n)} + b_2 s & 0 \\
      0 & 0 & \cM_{s,3}^{(n)} + b_3 s
     \end{bmatrix} , \qquad s \in \RR_+ , \qquad n \in \NN ,
 \]
 and hence
 \begin{align*}
  \int_0^t \bcR^{(n)}_s \, \dd s
  &= \frac{1}{n^2} \sum_{k=0}^{\nt-1}
      \begin{bmatrix}
        X_{k,1} & 0 & 0 \\
        0 & X_{k,2} & 0 \\
        0 & 0 & X_{k,3} \\
      \end{bmatrix}
     + \frac{nt-\nt}{n^2}
     \begin{bmatrix}
       X_{\nt,1} & 0 & 0 \\
       0 & X_{\nt,2} & 0 \\
       0 & 0 & X_{\nt,3}
     \end{bmatrix} \\
  &\quad
     + \frac{\nt+(nt-\nt)^2}{2n^2}
        \begin{bmatrix}
         b_1 & 0 & 0 \\
         0 & b_2 & 0 \\
         0 & 0 & b_3 \\
        \end{bmatrix},
   \qquad t \in \RR_+ , \qquad n \in \NN ,
 \end{align*}
 see, e.g., the proof of Theorem 1.1 in Barczy et al.\ \cite{BarBezPap1}.
By Lemma \ref{Moments},
 \[
   \frac{1}{n^2}
   \sum_{k=1}^\nt
    \EE(\bM_k \bM_k^\top \mid \cF_{k-1}^\bX)
   = \frac{\nt}{n^2} \bV^{(0)}
     + \frac{1}{n^2}
       \sum_{k=1}^\nt \Big(X_{k-1,1} \bV^{(1)} + X_{k-1,2} \bV^{(2)}  +  X_{k-1,3} \bV^{(3)}  \Big)
 \]
 for all $t \in \RR_+$ and $n \in \NN$.
Since $\xi_{1,1,2,1} \ase 0$, $\xi_{1,1,3,1} \ase 0$, $\xi_{1,1,1,2} \ase 0$, $\xi_{1,1,3,2} \ase 0$, $\xi_{1,1,1,3} \ase 0$
 and $\xi_{1,1,2,3} \ase 0$  (due to \ $a_{1,2} = a_{1,3} =  a_{2,1} = a_{2,3}= a_{3,1} = a_{3,2} = 0$), we have
 $v^{(1)}_{i,j} = 0$, $i,j\in\{1,2,3\}$, $(i,j)\ne (1,1)$,  $v^{(2)}_{i,j} = 0$, $i,j\in\{1,2,3\}$, $(i,j)\ne (2,2)$, and
 $v^{(3)}_{i,j} = 0$, $i,j\in\{1,2,3\}$, $(i,j)\ne (3,3)$.
Consequently, we obtain
 \begin{align*}
  &X_{k-1,1} \bV^{(1)} + X_{k-1,2} \bV^{(2)} + X_{k-1,3} \bV^{(3)}\\
  & = X_{k-1,1} \begin{bmatrix}
                v^{(1)}_{1,1} & 0 & 0 \\
                0 & 0 & 0 \\
                0 & 0 & 0
               \end{bmatrix}
     + X_{k-1,2} \begin{bmatrix}
                  0 & 0 & 0 \\
                  0 & v^{(2)}_{2,2} & 0 \\
                  0 & 0 & 0
                 \end{bmatrix}
     + X_{k-1,3} \begin{bmatrix}
                  0 & 0 & 0 \\
                  0 & 0 & 0 \\
                  0 & 0 & v^{(3)}_{3,3}
                 \end{bmatrix}\\
  &= \begin{bmatrix}
        X_{k-1,1} & 0 & 0 \\
        0 & X_{k-1,2} & 0 \\
        0 & 0 &  X_{k-1,3}
     \end{bmatrix}
     \begin{bmatrix}
      v^{(1)}_{1,1} & 0 & 0\\
      0 & v^{(2)}_{2,2} & 0 \\
     0 & 0 & v^{(3)}_{3,3}
     \end{bmatrix}
   = \begin{bmatrix} X_{k-1,1} & 0 & 0 \\ 0 & X_{k-1,2} & 0 \\ 0 & 0 & X_{k-1,3} \end{bmatrix}
     \bV_\bxi , \qquad k\in\NN.
 \end{align*}
So
 \begin{align*}
  & \frac{1}{n^2}
   \sum_{k=1}^\nt
    \EE(\bM_k \bM_k^\top \mid \cF_{k-1}^\bX)
     - \int_0^t \bcR^{(n)}_s \bV_\bxi \, \dd s   \\
  &\qquad  = \frac{\nt}{n^2}\bV^{(0)}
      - \frac{nt-\nt}{n^2}
       \begin{bmatrix} X_{\nt,1} & 0 & 0 \\ 0 & X_{\nt,2} & 0 \\ 0 & 0 & X_{\nt,3} \end{bmatrix}
       \bV_\bxi
       - \frac{\nt + (nt - \nt)^2}{2n^2} \begin{bmatrix} b_1 & 0 & 0 \\ 0 & b_2 & 0 \\ 0 & 0 & b_3 \end{bmatrix}
         \bV_\bxi
 \end{align*}
 for all $t\in\RR_+$ and $n\in\NN$.
Hence, in order to show \eqref{Condb1}, by Slutsky's lemma and taking into account the facts that
 for all $T\in\RR_{++}$,
 \[
   \sup_{t\in[0,T]}
      \frac{\lfloor nt\rfloor + (nt - \lfloor nt\rfloor)^2}{n^2}
    \leq \sup_{t\in[0,T]}  \frac{\lfloor nt\rfloor + 1}{n^2} \to 0
      \qquad \text{as \ $n\to\infty$,}
 \]
 and $\sup_{t\in[0,T]}\frac{\lfloor nt\rfloor}{n^2}\bV^{(0)}\to\bzero$ as $n\to\infty$, it suffices to prove
 that for all $T\in\RR_{++}$, we have
 \begin{equation}\label{Condb11}
  \frac{1}{n^2}
  \sup_{t \in [0,T]}
   \| (nt-\nt ) \bX_\nt\|
  \leq
  \frac{1}{n^2}
  \sup_{t \in [0,T]}
   \|\bX_\nt\|
  \stoch 0 \qquad \text{as \ $n \to \infty$.}
 \end{equation}
To prove \eqref{Condb11}, it is enough to show that
\begin{equation*}
\EE\left(
  \frac{1}{n^4}
  \sup_{t \in [0,T]}
   \left \|
       \bX_\nt \right \|^2
\right)\to0 \qquad \text{as $n\to\infty$.}
\end{equation*}
By \eqref{Xdeco}, for each $k\in\NN$ we have
\begin{equation*}
  \begin{bmatrix}
   X_{k,1} \\
   X_{k,2}  \\
   X_{k,3}
  \end{bmatrix}
   =
  \begin{bmatrix}
   X_{k,1}^{(1)} \\
   X_{k,2}^{(2)}  \\
   X_{k,3}^{(4)}
  \end{bmatrix}
   =
  \begin{bmatrix}
   \sum_{\ell=1}^k(M_{\ell,1}+b_1) \\
   \sum_{\ell=1}^k(M_{\ell,2}+b_2) \\
   \sum_{\ell=1}^k(M_{\ell,3}+b_3)
  \end{bmatrix},
\end{equation*}
 which together with Lemma \ref{M_sum_growth} and $\eta_1=\eta_2=\eta_3=1$ (following from Lemma \ref{eta_values}) yield that
\begin{align*}
\EE\left(
  \frac{1}{n^4}
  \sup_{t \in [0,T]}
   \left \|
       \bX_{\nt} \right \|^2
\right)
\leq&
  \frac{1}{n^4}
\EE\left(
  \sup_{t \in [0,T]}
   \left(X_{\nt,1}^{(1)}\right)^2
+\sup_{t \in [0,T]}
        \left(X_{\nt,2}^{(2)}\right)^2
+\sup_{t \in [0,T]}
        \left(X_{\nt,3}^{(4)}\right)^2
\right)\\
=&\frac{1}{n^4}(\OO(n^{\eta_1+1})+\OO(n^{\eta_2+1})+\OO(n^{\eta_3+1}))=\OO(n^{-2})\to0\qquad \text{as $n\to\infty$,}
\end{align*}
implying \eqref{Condb11}, and hence \eqref{Condb1}.

{\sl Step 2/(c).}
We check condition \eqref{Condb2}.
We show that for all $T\in\RR_{++}$ and $\theta\in\RR_{++}$,
 \[
   \frac{1}{n^2}
  \sum_{k=1}^\nT
   \EE(\|\bM_k\|^2 \bbone_{\{\|\bM_k\|>n\theta\}} \mid \cF_{k-1}^\bX)
   \mean 0 \qquad \text{as \ $n\to\infty$.}
 \]
 Using the inequalities $\bbone_{\{\|\bM_k\|>\alpha\}}<\frac{\|\bM_k\|^2}{\alpha^2}$, $k\in\NN, \alpha\in\RR_{++}$, and $(a+b+c)^2 \leq 3(a^2+b^2+c^2)$, $a,b,c\in\RR_+$, and Lemma \ref{M4moment} (which yields that $\EE(M_{k,i}^4) = \OO(k^2)$, $k\in\NN$ for $i\in\{1,2,3\}$),
 for all $T\in\RR_{++}$ and $\theta\in\RR_{++}$, we have
 \begin{align*}
  &\EE\left(  \frac{1}{n^2}
        \sum_{k=1}^\nT
         \EE(\|\bM_k\|^2 \bbone_{\{\|\bM_k\|>n\theta\}} \mid \cF_{k-1}^\bX)
     \right)
   = \frac{1}{n^2}
      \sum_{k=1}^\nT \EE(\|\bM_k\|^2 \bbone_{\{\|\bM_k\|>n\theta\}}) \\
  & \leq  \frac{1}{n^2}
          \sum_{k=1}^\nT \EE\left( \frac{\|\bM_k\|^4}{n^2\theta^2} \right)
   \leq \frac{3}{n^4 \theta^2} \sum_{k=1}^\nT \EE(M_{k,1}^4 + M_{k,2}^4 + M_{k,3}^4 )
    =   \frac{1}{n^4\theta^2} \sum_{k=1}^\nT \OO(k^2)
    = \OO(n^{-1})\to 0
 \end{align*}
 as $n\to\infty$.

{\sl Step 3.}
Using \eqref{conv_bM} and Lemma \ref{Conv2Funct}, we can prove \eqref{Conv_X_1}.
For each $n\in\NN$, by \eqref{help15}, we have $(n^{-1}\bX_{\nt})_{t\in\RR_+} = \Psi^{(n)}(\bcM^{(n)})$,
 where the mapping $\Psi^{(n)} : \DD(\RR_+, \RR^3) \to \DD(\RR_+, \RR^3)$ is given by
 \[
   (\Psi^{(n)}(f))(t)
    := f\biggl(\frac{\nt}{n}\biggr) + \frac{\nt}{n} \bb
 \]
 for $f \in \DD(\RR_+, \RR^3)$ and $t \in \RR_+$.
Further, using that $(\bcM_t +\bb t)_{t\in\RR_+} \distre (\bcX_t)_{t\in\RR_+}$, we have
 $\bcX \distre \Psi(\bcM)$, where the mapping
 $\Psi : \DD(\RR_+, \RR^3) \to \DD(\RR_+, \RR^3)$ is given by
 \[
   (\Psi(f))(t) := f(t) + \bb t , \qquad
   f \in \DD(\RR_+, \RR^3) , \qquad t \in \RR_+.
 \]
The mappings $\Psi^{(n)}$, $n\in\NN$, and $\Psi$ are measurable, which can be checked in the same way
 as in Step 4/(a) in Barczy et al.\ \cite{BarBezPap1} replacing $\DD(\RR_+,\RR)$ by $\DD(\RR_+,\RR^3)$
 in the argument given there.
One can also check that the set $C := \CC(\RR_+, \RR^3)$ satisfies $C \in \cB(\DD(\RR_+, \RR^3))$, $\PP(\bcM \in C) = 1$, and
 $\Psi^{(n)}(f^{(n)}) \to \Psi(f)$ in $\DD(\RR_+, \RR^3)$ as $n \to \infty$ if $f^{(n)} \to f$ in $\DD(\RR_+, \RR^3)$ as $n \to \infty$
 with $f \in C$, \ $f^{(n)}\in \DD(\RR_+, \RR^3)$, $n\in\NN$.
Namely, one can follow the same argument as in Step 4/(b) in Barczy et al.\ \cite{BarBezPap1} replacing
 $\DD(\RR_+,\RR)$ by $\DD(\RR_+,\RR^3)$, and $\CC(\RR_+,\RR)$ by $\CC(\RR_+,\RR^3)$, respectively,
 in the argument given there.
So we can apply Lemma \ref{Conv2Funct}, and we obtain
 $(n^{-1} \bX_\nt)_{t\in\RR_+} = \Psi^{(n)}(\bcM^{(n)}) \distr \Psi(\bcM)$ as $n \to \infty$, where
 $( (\Psi(\bcM))(t))_{t\in\RR_+}  = (\bcM_t + \bb t)_{t\in\RR_+} \distre (\bcX_t)_{t\in\RR_+}$, as desired.

\section{Proof of Theorem \ref{main_2}}\label{proof2}

By \eqref{Xdeco}, we have the decomposition
 \begin{align}\label{help_proof_Thm3_3}
     \begin{bmatrix}
        X_{\nt,1} \\
        X_{\nt,2} \\
        X_{\nt,3}
      \end{bmatrix}
   = \begin{bmatrix}
        X_{\nt,1}^{(1)} \\[1mm]
         X^{(2)}_{\nt,2} \\[1mm]
        a_{3,1}X_{\nt,3}^{(2)}+a_{3,2}X_{\nt,3}^{(3)}+X_{\nt,3}^{(4)}
     \end{bmatrix} , \qquad
   t \in \RR_+ , \quad n \in \NN,
 \end{align}
 where
 \begin{gather*}
    X^{(1)}_{\nt,1} = \sum_{j=1}^\nt  (M_{j,1} + b_1), \qquad X^{(2)}_{\nt,2} = \sum_{j=1}^\nt (M_{j,2} + b_2),
                                                  \qquad X_{\nt,3}^{(4)} = \sum_{j=1}^\nt (M_{j,3} + b_3),
\end{gather*}
 and due to $X_{\nt,1}=X_{\nt,1}^{(1)}$, $X_{\nt,2}=X_{\nt,2}^{(2)}$ and \eqref{integrals},
\begin{gather*}
     X^{(2)}_{\nt,3}=n^2\int_0^{ \frac{\nt}{n}} n^{-1}X_{\ns,1} \, \dd s, \qquad     X^{(3)}_{\nt,3} =n^2\int_0^{ \frac{\nt}{n}} n^{-1}X_{\ns,2} \, \dd s
 \end{gather*}
 for $t\in\RR_+$ and $n\in\NN$.

{\sl Step 1/(a).} First, we prove that
\begin{equation}\label{Conv_X_2_Thm2}
  \left(\begin{bmatrix}
          n^{-1} X_{\nt,1} \\
          n^{-1} X_{\nt,2}
         \end{bmatrix}\right)_{t\in\RR_+}
  \distr
  \left(\begin{bmatrix}
          \cX_{t,1} \\
          \cX_{t,2}
         \end{bmatrix}\right)_{t\in\RR_+}
  \qquad \text{as \ $n \to \infty$.}
 \end{equation}
The SDE
 \begin{equation}\label{SDE_X_2_mod}
  \begin{cases}
   \dd\cX_{t,1}
   = b_1 \, \dd t + \sqrt{v^{(1)}_{1,1} \cX_{t,1}^+} \, \dd \cW_{t,1} , \\[2mm]
   \dd\cX_{t,2}
   = b_2 \, \dd t + \sqrt{v^{(2)}_{2,2} \cX_{t,2}^+} \, \dd \cW_{t,2} ,
  \end{cases}
  \qquad t \in\RR_+ ,
 \end{equation}
 has a pathwise unique strong solution $( [\cX_{t,1},\cX_{t,2}]^\top )_{t\in\RR_+}$ for all initial values
 $ [\cX_{0,1},\cX_{0,2}]^\top = [x_1,x_2]^\top \in \RR^2$,
 where $(\cW_{t,1})_{t\in\RR_+}$ and $(\cW_{t,2})_{t\in\RR_+}$ are independent standard Wiener processes.
Further, if $ [x_1,x_2]^\top\in \RR_+^2$, then $ [\cX_{t,1},\cX_{t,2}]^\top \in \RR_+^2$ almost surely for all $t \in \RR_+$,
  since $b_i$, $v^{(i)}_{i,i}\in \RR_+$, $i\in\{1,2\}$, see, e.g., Ikeda and Watanabe \cite[Chapter IV, Example 8.2]{IkeWat}.
Note that the SDE \eqref{SDE_X_2_mod} is nothing else but the first and second equations of the SDE \eqref{SDE_X_2}.
By \eqref{help_Wei_Winnicki1}, we have $(n^{-1} X_{\nt,1})_{t\in\RR_+} \distr (\cX_{t,1})_{t\in\RR_+}$ as $n \to \infty$,
 where $(\cX_{t,1})_{t\in\RR_+}$ satisfies the first equation of the SDE \eqref{SDE_X_2} with initial value $\cX_{0,1}=0$.
Similarly, since $a_{2,1}=0$, as it was explained in Section \ref{Section_conv_results}, the second coordinate process
 $(X_{k,2})_{k\in\ZZ_+}$ is a critical single-type GWI process,
 so $(n^{-1} X_{\nt,2})_{t\in\RR_+} \distr (\cX_{t,2})_{t\in\RR_+}$ as $n \to \infty$,
 where $(\cX_{t,2})_{t\in\RR_+}$  satisfies the second equation of the SDE \eqref{SDE_X_2}
 with initial value $\cX_{0,2}=0$.
In what follows, we will prove the joint convergence of $(n^{-1} X_{\nt,1})_{t\in\RR_+}$
 and $(n^{-1} X_{\nt,2})_{t\in\RR_+}$ as $n\to\infty$, namely, we prove \eqref{Conv_X_2_Thm2}.
The proof goes along the same lines as that of Theorem \ref{main_1}.

{\sl Step 1/(b).}
Using $a_{2,1}=0$ and \eqref{Mk1}, we obtain that
 \[
   \begin{bmatrix}
        M_{k,1} \\
        M_{k,2} \\
   \end{bmatrix}
   = \begin{bmatrix}
        X_{k,1} \\
        X_{k,2} \\
      \end{bmatrix}
      -\begin{bmatrix}
        X_{k-1,1} \\
        X_{k-1,2} \\
      \end{bmatrix}
     - \begin{bmatrix}
        b_1 \\
        b_2 \\
      \end{bmatrix},
    \qquad k\in\NN.
 \]
Consider the random step processes
 \begin{align}\label{help15_Thm2}
    \begin{bmatrix} \cM_{t,1}^{(n)} \\ \cM_{t,2}^{(n)}  \end{bmatrix}
   := \frac{1}{n} \sum_{k=1}^\nt
      \begin{bmatrix}
        M_{k,1} \\
        M_{k,2} \\
      \end{bmatrix}
   = \frac{1}{n}
       \begin{bmatrix}
         X_{\nt,1} \\
         X_{\nt,2} \\
       \end{bmatrix}
      - \frac{\nt}{n}
       \begin{bmatrix}
         b_1 \\
         b_2 \\
       \end{bmatrix}
	, \qquad t \in \RR_+ , \qquad n \in \NN,
 \end{align}
where at the second equality we used that $X_{0,1}=X_{0,2}=0$.
We show that
 \begin{equation}\label{conv_bM_Thm2}
  \left(\begin{bmatrix} \cM_{t,1}^{(n)} \\ \cM_{t,2}^{(n)}  \end{bmatrix}\right)_{t\in\RR_+}
   \distr \left(
       \begin{bmatrix}
        \cM_{t,1} \\
        \cM_{t,2} \\
      \end{bmatrix} \right)_{t\in\RR_+} \qquad
  \text{as \ $n \to \infty$,}
 \end{equation}
 where the limit process $([\cM_{t,1}, \cM_{t,2}]^\top)_{t \in \RR_+}$ is the pathwise unique strong solution of the SDE
 \begin{equation}\label{SDE_bM_Thm2}
  \begin{cases}
   \dd \cM_{t,1}
   = \sqrt{v^{(1)}_{1,1} (\cM_{t,1} + b_1 t)^+ } \, \dd \cW_{t,1} , \qquad t\in\RR_+,\\[1mm]
  \dd \cM_{t,2}
   = \sqrt{v^{(2)}_{2,2} (\cM_{t,2} + b_2 t)^+ } \, \dd \cW_{t,2} , \qquad t\in\RR_+
  \end{cases}
 \end{equation}
 with initial value $[\cM_{0,1},\cM_{0,2}] = [0,0]$.

{\sl Step 1/(c).}
We check that the SDE \eqref{SDE_bM_Thm2} has a pathwise unique strong solution $([\cM_{t,1}, \cM_{t,2}]^\top)_{t\in\RR_+}$
 for all initial values $[\cM_{0,1},\cM_{0,2}]^\top  = \bx = [x_1,x_2]^\top \in \RR^2$.
Observe that if $([\cM_{t,1}, \cM_{t,2}]^\top)_{t\in\RR_+}$ is a strong solution of the SDE \eqref{SDE_bM_Thm2} with initial value
 $[\cM_{0,1},\cM_{0,2}]^\top  = \bx = [x_1,x_2]^\top \in \RR^2$, then, by It\^o's formula, the process
 \[
   \begin{bmatrix}
      \cP_{t,1} \\
      \cP_{t,2} \\
   \end{bmatrix}
   :=  \begin{bmatrix}
         \cM_{t,1} \\
         \cM_{t,2} \\
       \end{bmatrix}
      + \begin{bmatrix}
         b_1 \\
         b_2 \\
       \end{bmatrix} t, \qquad t \in \RR_+,
 \]
  is a pathwise unique strong solution of the SDE
 \begin{equation}\label{SDE_bP_Thm2}
  \begin{cases}
   \dd \cP_{t,1}
   = b_1 \, \dd t
     + \sqrt{v^{(1)}_{1,1} \, \cP_{t,1}^+} \, \dd \cW_{t,1}, \\[1mm]
  \dd \cP_{t,2}
   = b_2 \, \dd t
     + \sqrt{v^{(2)}_{2,2} \, \cP_{t,2}^+} \, \dd \cW_{t,2},
  \end{cases} \qquad t\in\RR_+,
 \end{equation}
 with initial value $[\cP_{0,1}, \cP_{0,2}]^\top = \bx$.
Conversely, if $[\cP_{t,1}, \cP_{t,2}]^\top$, $t \in \RR_+$, is a strong solution
 of the SDE \eqref{SDE_bP_Thm2} with initial value $[\cP_{0,1}, \cP_{0,2}]^\top = \bp \in \RR^2$,
 then, by It\^o's formula, the process $[\cP_{t,1}, \cP_{t,2}]^\top - [b_1, b_2]^\top t$, $t \in \RR_+$,
 is a strong solution of the SDE \eqref{SDE_bM_Thm2} with initial value $[\cP_{0,1}, \cP_{0,2}]^\top - [b_1, b_2]^\top 0 = \bp$.
The equations in \eqref{SDE_bP_Thm2} are the same as in \eqref{SDE_X_2_mod},
thus, as it was explained in Step 1/(a), the SDE \eqref{SDE_bP_Thm2}  admits
a unique strong solution with an arbitrary initial value in $\RR^2$.
Consequently, the SDE \eqref{SDE_bM_Thm2} (with initial value $[0,0]^\top\in\RR^2$)
 admits a unique strong solution $([\cM_{t,1},\cM_{t,2}]^\top)_{t\in\RR_+}$, and
 \[
  \left( \begin{bmatrix}
         \cM_{t,1} \\
         \cM_{t,2} \\
       \end{bmatrix}
      + \begin{bmatrix}
         b_1 \\
         b_2 \\
       \end{bmatrix} t \right)_{t\in\RR_+}
   \distre
  \left(  \begin{bmatrix}
          \cX_{t,1} \\
          \cX_{t,2}
         \end{bmatrix}
   \right)_{t\in\RR_+}.
 \]

{\sl Step 2/(a).}
 In order to prove \eqref{conv_bM_Thm2}, we want to apply Theorem \ref{Conv2DiffThm} with the following choices:
$d = r = 2$, $\cF_k^{(n)} = \cF_k^\bX$, $n \in \NN$, $k \in \ZZ_+$,
 \begin{align*}
  \bU_k^{(n)}
    & = n^{-1}
       \begin{bmatrix}
         M_{k,1} \\
         M_{k,2} \\
       \end{bmatrix}, \quad n,k\in\NN, \qquad \bU_0^{(n)} =\bzero, \quad n\in\NN,\\
  \bcU_t
    & = \begin{bmatrix}
         \cM_{t,1} \\
         \cM_{t,2} \\
       \end{bmatrix}, \qquad t\in\RR_+
 \end{align*}
 (yielding $(\bcU_t^{(n)})_{t\in\RR_+} =  ([\cM_{t,1}^{(n)}, \cM_{t,2}^{(n)}]^\top)_{t\in\RR_+}$, $n\in\NN$),
  and $\bbeta : \RR_+ \times \RR^2 \to \RR^2$ and $\bgamma : \RR_+ \times \RR^2 \to \RR^{2\times2}$
 given by $\bbeta(t, \bx) = \bzero$ and
 \[
   \bgamma(t, \bx)
   =  \begin{bmatrix}
      \sqrt{v^{(1)}_{1,1}(x_1 + b_1 t)^+} & 0 \\
      0 & \sqrt{v^{(2)}_{2,2}(x_2 + b_2 t)^+} \\
     \end{bmatrix}
 \]
 for all $t \in \RR_+$ and $\bx = [x_1, x_2]^\top\in \RR^2$.

The convergence $\bU_0^{(n)}\distr \bzero$ as $n\to\infty$, and condition (i) of Theorem \ref{Conv2DiffThm} trivially holds,
 since
 \[
    \EE\left( \begin{bmatrix}
                  M_{k,1} \\
                  M_{k,2} \\
               \end{bmatrix}
              \; \Big \vert \; \cF_{k-1}^\bX\right)
        = \begin{bmatrix}
            0 \\
            0 \\
        \end{bmatrix}, \quad k\in\NN, \qquad \text{and}\qquad  \bbeta(t, \bx) = \bzero,\quad t\in\RR_+, \;\; \bx\in\RR^2.
 \]
Now, we show that conditions (ii) and (iii) of Theorem \ref{Conv2DiffThm} hold.
We have to check that for all $T \in \RR_{++}$,
 \begin{gather} \label{Condb1_Thm2}
  \sup_{t\in[0,T]}
   \biggl\Vert\frac{1}{n^2}
          \sum_{k=1}^\nt
           \EE\left( \begin{bmatrix}
                  M_{k,1} \\
                  M_{k,2} \\
               \end{bmatrix}
           \begin{bmatrix}
                  M_{k,1} \\
                  M_{k,2} \\
               \end{bmatrix}^\top \; \Big \vert \;  \cF_{k-1}^\bX \right)
          - \int_0^t \bcR^{(n;1,2)}_s \,\bV_\bxi^{(1,2)} \, \dd s \biggr\Vert
  \stoch 0 \qquad \text{as \ $n\to\infty$,} \\
  \frac{1}{n^2}
  \sum_{k=1}^\nT
   \EE\left(\left\|  \begin{bmatrix}
                  M_{k,1} \\
                  M_{k,2} \\
               \end{bmatrix}  \right\|^2
       \bbone_{\big\{ \|  [ M_{k,1}, M_{k,2} ]^\top
                      \|>n\theta \big\} } \; \Big \vert \;  \cF_{k-1}^\bX\right)
  \stoch 0  \qquad \text{as \ $n\to\infty$ \ for all \ $\theta \in \RR_{++}$,} \label{Condb2_Thm2}
 \end{gather}
 where the process $(\bcR^{(n;1,2)}_s)_{s \in \RR_+}$ and the matrix $\bV_\bxi^{(1,2)}$ are defined by
 \begin{align*}
  \bcR^{(n;1,2)}_s
  &:= \begin{bmatrix}
      (\cM_{s,1}^{(n)} + b_1 s)^+ & 0 \\
      0 & (\cM_{s,2}^{(n)} + b_2 s)^+ \\
     \end{bmatrix} , \qquad s \in \RR_+ , \qquad n \in \NN , \\
 \bV_\bxi^{(1,2)}
  &:= \begin{bmatrix}
      v^{(1)}_{1,1} & 0 \\
      0 & v^{(2)}_{2,2} \\
     \end{bmatrix} .
 \end{align*}
Indeed, $\EE([M_{k,1}, M_{k,2}]\mid \cF_{k-1}^{\bX}) = [0,0]$ and thus
\begin{equation*}
\var\left(\begin{bmatrix}M_{k,1}\\M_{k,2}\end{bmatrix}\,\Big\vert\,\cF_{k-1}^{\bX}\right)
	=\EE\left(\begin{bmatrix}M_{k,1}\\M_{k,2}\end{bmatrix}
		\begin{bmatrix}M_{k,1}\\M_{k,2}\end{bmatrix}^\top\,\Big\vert\,\cF_{k-1}^{\bX}\right), \qquad k\in\NN,
\end{equation*}
and, for all $t\in\RR_+$ and $\bx=(x_1,x_2)\in\RR^2$, we have
 \begin{align*}
  \gamma(t,\bx) \gamma(t,\bx)^\top
  & =  \begin{bmatrix}
      \sqrt{v^{(1)}_{1,1}(x_1 + b_1 t)^+} & 0 \\
      0 & \sqrt{v^{(2)}_{2,2}(x_2 + b_2 t)^+} \\
     \end{bmatrix}
     \begin{bmatrix}
      \sqrt{v^{(1)}_{1,1}(x_1 + b_1 t)^+} & 0 \\
      0 & \sqrt{v^{(2)}_{2,2}(x_2 + b_2 t)^+} \\
     \end{bmatrix} \\
  & = \begin{bmatrix}
       v^{(1)}_{1,1} (x_1 + b_1 t)^+ & 0 \\
        0 & v^{(2)}_{2,2}(x_2 + b_2 t)^+
       \end{bmatrix}
   =  \begin{bmatrix}
        (x_1 + b_1 t)^+ & 0 \\
         0 & (x_2 + b_2 t)^+
       \end{bmatrix}
       \bV_\bxi^{(1,2)}.
 \end{align*}

{\sl Step 2/(b).}
 Now we check \eqref{Condb1_Thm2}. By \eqref{help15_Thm2}, for all $s \in \RR_+$ and $n \in \NN$, we have
 \[
   \begin{bmatrix}
         \cM^{(n)}_{s,1} \\[1mm]
         \cM^{(n)}_{s,2} \\
       \end{bmatrix}
   +  \begin{bmatrix}
         b_1 \\
         b_2 \\
       \end{bmatrix} s
   = \frac{1}{n}
       \begin{bmatrix}
         X_{\ns,1} \\
         X_{\ns,2} \\
       \end{bmatrix}
    + \frac{ns-\ns}{n}
      \begin{bmatrix}
         b_1 \\
         b_2 \\
       \end{bmatrix},
 \]
 where both coordinates of both terms of the sum on the right hand side of the above equality are positive,
  and thus
 \[
   \bcR^{(n;1,2)}_s
   = \begin{bmatrix}
      \cM_{s,1}^{(n)} + b_1 s & 0 \\
      0 & \cM_{s,2}^{(n)} + b_2 s
     \end{bmatrix} , \qquad s \in \RR_+ , \qquad n \in \NN .
 \]
Hence
 \begin{align*}
  \int_0^t \bcR^{(n;1,2)}_s \, \dd s
  &= \frac{1}{n^2} \sum_{k=0}^{\nt-1}
      \begin{bmatrix}
        X_{k,1} & 0 \\
        0 & X_{k,2} \\
      \end{bmatrix}
     + \frac{nt-\nt}{n^2}
     \begin{bmatrix}
       X_{\nt,1} & 0 \\
       0 & X_{\nt,2}
     \end{bmatrix} \\
  &\quad
     + \frac{\nt+(nt-\nt)^2}{2n^2}
        \begin{bmatrix}
         b_1 & 0 \\
         0 & b_2 \\
        \end{bmatrix},
   \qquad t \in \RR_+ , \qquad n \in \NN ,
 \end{align*}
 as, e.g., in the proof of Theorem 1.1 in Barczy et al.\ \cite{BarBezPap1}.
By Lemma \ref{Moments},
 \begin{align*}
   &\frac{1}{n^2}
    \sum_{k=1}^\nt
      \EE\left( \begin{bmatrix}
                  M_{k,1} \\
                  M_{k,2} \\
               \end{bmatrix}
           \begin{bmatrix}
                  M_{k,1} \\
                  M_{k,2} \\
               \end{bmatrix}^\top \; \Big \vert \;  \cF_{k-1}^\bX \right) \\
   &\qquad  = \frac{\nt}{n^2} \bV^{(0;1,2)}
                + \frac{1}{n^2}
                  \sum_{k=1}^\nt \Big(X_{k-1,1} \bV^{(1;1,2)} + X_{k-1,2} \bV^{(2;1,2)} +  X_{k-1,3} \bV^{(3;1,2)} \Big)
 \end{align*}
 for all $t \in \RR_+$ and $n \in \NN$, where, for each $i\in\{ 0,1,2,3\}$, $\bV^{(i;1,2)}$ denotes the matrix obtained from $\bV^{(i)}$
 by deleting its third row and third column.
Since $\xi_{1,1,2,1} \ase 0$, $\xi_{1,1,1,2} \ase 0$, $\xi_{1,1,3,1}\ase0$, and $\xi_{1,1,3,2}\ase0$ (due to \ $a_{1,2} = a_{2,1} = a_{1,3} = a_{2,3} = 0$), for each $i,j\in\{1,2\}$ we have
 $v^{(1)}_{i,j} =0$ if $(i,j)\ne(1,1)$, $v^{(2)}_{i,j} = 0$ if $(i,j)\ne(2,2)$, and $v^{(3)}_{i,j}= 0$.
Consequently, we obtain that
 \begin{align*}
  &X_{k-1,1} \bV^{(1;1,2)} + X_{k-1,2} \bV^{(2;1,2)} + X_{k-1,3} \bV^{(3;1,2)}\\
  & = X_{k-1,1} \begin{bmatrix}
                v^{(1)}_{1,1} & 0 \\
                0 & 0
               \end{bmatrix}
     + X_{k-1,2} \begin{bmatrix}
                  0 & 0 \\
                  0 & v^{(2)}_{2,2}
                 \end{bmatrix}
     + X_{k-1,3} \begin{bmatrix}
                  0 & 0 \\
                  0 & 0
                 \end{bmatrix}
  = X_{k-1,1} \begin{bmatrix}
                v^{(1)}_{1,1} & 0 \\
                0 & 0
               \end{bmatrix}
     + X_{k-1,2} \begin{bmatrix}
                  0 & 0 \\
                  0 & v^{(2)}_{2,2}
                 \end{bmatrix} \\
  &= \begin{bmatrix}
        X_{k-1,1} & 0 \\
        0 & X_{k-1,2}
     \end{bmatrix}
     \begin{bmatrix}
      v^{(1)}_{1,1} & 0\\
      0 & v^{(2)}_{2,2}
     \end{bmatrix}
   = \begin{bmatrix}
       X_{k-1,1} & 0 \\
       0 & X_{k-1,2}
     \end{bmatrix}
     \bV^{(1,2)}_\bxi, \qquad k\in\NN.
 \end{align*}
So
 \begin{align*}
  & \frac{1}{n^2}\sum_{k=1}^\nt
           \EE\left( \begin{bmatrix}
                  M_{k,1} \\
                  M_{k,2} \\
               \end{bmatrix}
           \begin{bmatrix}
                  M_{k,1} \\
                  M_{k,2} \\
               \end{bmatrix}^\top \; \Big \vert \;  \cF_{k-1}^\bX \right)
          - \int_0^t \bcR^{(n;1,2)}_s \,\bV_\bxi^{(1,2)} \, \dd s  \\
  &\qquad  = \frac{\nt}{n^2}\bV^{(0;1,2)}
      - \frac{nt-\nt}{n^2}
       \begin{bmatrix} X_{\nt,1} & 0  \\ 0 & X_{\nt,2} \end{bmatrix}
       \bV_\bxi^{(1,2)}
       - \frac{\nt + (nt - \nt)^2}{2n^2} \begin{bmatrix} b_1 & 0 \\ 0 & b_2 \end{bmatrix}
         \bV_\bxi^{(1,2)}
 \end{align*}
 for all $t\in\RR_+$ and $n\in\NN$.
Hence, in order to show \eqref{Condb1_Thm2}, by Slutsky's lemma and taking into account the facts that
 for all $T\in\RR_{++}$,
 \[
   \sup_{t\in[0,T]}
      \frac{\lfloor nt\rfloor + (nt - \lfloor nt\rfloor)^2}{n^2}
    \leq \sup_{t\in[0,T]}  \frac{\lfloor nt\rfloor + 1}{n^2} \to 0
      \qquad \text{as \ $n\to\infty$,}
 \]
 and $\sup_{t\in[0,T]}\frac{\lfloor nt\rfloor}{n^2}\bV^{(0;1,2)}\to\bzero$ as $n\to\infty$, it suffices to prove
 that for all $T\in\RR_{++}$, we have
 \begin{equation}\label{Condb11_Thm2}
  \frac{1}{n^2}
  \sup_{t \in [0,T]}
   \left\| (nt-\nt )
       \begin{bmatrix} X_{\nt,1} & 0  \\ 0 & X_{\nt,2} \end{bmatrix}
      \right\|
  \leq
  \frac{1}{n^2}
  \sup_{t \in [0,T]}
   \left \|
       \begin{bmatrix}
         X_{\nt,1} \\
         X_{\nt,2}  \\
       \end{bmatrix} \right \|
  \stoch 0 \qquad \text{as \ $n \to \infty$.}
 \end{equation}
To prove \eqref{Condb11_Thm2}, it is enough to show that
\begin{equation*}
\EE\left(
  \frac{1}{n^4}
  \sup_{t \in [0,T]}
   \left \|
       \begin{bmatrix}
         X_{\nt,1} \\
         X_{\nt,2}  \\
       \end{bmatrix} \right \|^2
\right)\to0 \qquad \text{as $n\to\infty$.}
\end{equation*}
By \eqref{Xdeco}, for each $k\in\NN$ we have
\begin{equation*}
  \begin{bmatrix}
   X_{k,1} \\
   X_{k,2}  \\
  \end{bmatrix}
   =
  \begin{bmatrix}
   X_{k,1}^{(1)} \\
   X_{k,2}^{(2)}  \\
  \end{bmatrix}
   =
  \begin{bmatrix}
   \sum_{\ell=1}^k(M_{\ell,1}+b_1) \\
   \sum_{\ell=1}^k(M_{\ell,2}+b_2) \\
  \end{bmatrix},
\end{equation*}
which together with Lemma \ref{M_sum_growth} and $\eta_1=\eta_2=1$ (following from Lemma \ref{eta_values}) yield that
\begin{align*}
\EE\left(
  \frac{1}{n^4}
  \sup_{t \in [0,T]}
   \left \|
       \begin{bmatrix}
         X_{\nt,1} \\
         X_{\nt,2}  \\
       \end{bmatrix} \right \|^2
\right)
\leq&
  \frac{1}{n^4}
\EE\left(
  \sup_{t \in [0,T]}
   \left( X_{\nt,1}^{(1)} \right)^2
+\sup_{t \in [0,T]}
        \left( X_{\nt,2}^{(2)}\right)^2
\right)\\
=&\frac{1}{n^4}(\OO(n^{\eta_1+1})+\OO(n^{\eta_2+1}))=\OO(n^{-2})\to0\qquad \text{as $n\to\infty$,}
\end{align*}
implying \eqref{Condb11_Thm2}, and hence \eqref{Condb1_Thm2}.

{\sl Step 2/(c).} Next, we check condition \eqref{Condb2_Thm2}.
We show that for all $T\in\RR_{++}$ and $\theta\in\RR_{++}$,
 \[
  \frac{1}{n^2}
  \sum_{k=1}^\nT
   \EE\left(\left\|  \begin{bmatrix}
                  M_{k,1} \\
                  M_{k,2} \\
               \end{bmatrix}  \right\|^2
       \bbone_{\big\{ \|   [ M_{k,1}, M_{k,2} ]^\top
                      \|>n\theta \big\} } \; \Big \vert \;  \cF_{k-1}^\bX\right)
   \mean 0 \qquad \text{as \ $n\to\infty$.}
 \]
 Using the inequalities
$\bbone_{\{\|  [ M_{k,1}, M_{k,2} ]^\top \|>\alpha\}}
<\frac{\|  [ M_{k,1}, M_{k,2} ]^\top\|^2}{\alpha^2}$,
 $k\in\NN, \alpha\in\RR_{++}$, and $(a+b)^2 \leq 2(a^2+b^2)$, $a,b\in\RR_+$, and Lemma \ref{M4moment} (which yields that $\EE(M_{k,i}^4) = \OO(k^2)$, $k\in\NN$ for $i\in\{1,2\}$),
 for all $T\in\RR_{++}$ and $\theta\in\RR_{++}$, we have
 \begin{align*}
  &\EE\left(  \frac{1}{n^2}
        \sum_{k=1}^\nT
         \EE\left(\left\|  \begin{bmatrix}
                  M_{k,1} \\
                  M_{k,2} \\
               \end{bmatrix}  \right\|^2
       \bbone_{\big\{ \|  [ M_{k,1}, M_{k,2} ]^\top
                      \|>n\theta \big\} } \; \Big \vert \;  \cF_{k-1}^\bX\right)
     \right)\\
  & = \frac{1}{n^2}
      \sum_{k=1}^\nT \EE\left(  \left\|  \begin{bmatrix}
                  M_{k,1} \\
                  M_{k,2} \\
               \end{bmatrix}  \right\|^2
       \bbone_{\big\{ \|   [ M_{k,1},M_{k,2} ]^\top
                      \|>n\theta \big\} } \right) \\
  & \leq  \frac{1}{n^2}
          \sum_{k=1}^\nT \EE\left( \frac{1}{n^2\theta^2} \left\|
               \begin{bmatrix}
                  M_{k,1} \\
                  M_{k,2} \\
               \end{bmatrix}  \right\|^4   \right)
   \leq \frac{2}{n^4 \theta^2} \sum_{k=1}^\nT \EE(M_{k,1}^4 + M_{k,2}^4)
    =   \frac{1}{n^4\theta^2} \sum_{k=1}^\nT \OO(k^2)
    = \OO(n^{-1})\to 0
 \end{align*}
 as $n\to\infty$.

{\sl Step 3.} Using \eqref{conv_bM_Thm2} and Lemma \ref{Conv2Funct}, we can prove \eqref{Conv_X_2_Thm2}.
For each $n\in\NN$, by \eqref{help15_Thm2}, we have
 \[
  \left(n^{-1}
       \begin{bmatrix}
         X_{\nt,1} \\
         X_{\nt,2}  \\
       \end{bmatrix} \right)_{t\in\RR_+}
     = \Psi^{(n)}\left(
       \left(\begin{bmatrix}
        \cM_{t,1}^{(n)} \\
        \cM_{t,2}^{(n)} \\
      \end{bmatrix}\right)_{t\in\RR_+}
       \right),
  \]
 where the mapping $\Psi^{(n)} : \DD(\RR_+, \RR^2) \to \DD(\RR_+, \RR^2)$ is given by
 \[
   (\Psi^{(n)}(f))(t)
    := f\biggl(\frac{\nt}{n}\biggr)
         + \frac{\nt}{n}
         \begin{bmatrix}
          b_1 \\
          b_3 \\
      \end{bmatrix}
 \]
 for $f \in \DD(\RR_+, \RR^2)$ and $t \in \RR_+$.
Further, using that
 \[
   \left( \begin{bmatrix}
        \cM_{t,1} + b_1t\\
        \cM_{t,2} + b_2t\\
      \end{bmatrix} \right )_{t\in\RR_+}
   \distre \left(
       \begin{bmatrix}
        \cX_{t,1} \\
        \cX_{t,2} \\
      \end{bmatrix} \right)_{t\in\RR_+},
 \]
 we have
 \[
   \left(
       \begin{bmatrix}
        \cX_{t,1} \\
        \cX_{t,2} \\
   \end{bmatrix} \right)_{t\in\RR_+}
   \distre
   \Psi\left( \begin{bmatrix}
        \cM_{t,1} \\
        \cM_{t,2} \\
      \end{bmatrix}\right),
 \]
 where the mapping
 $\Psi : \DD(\RR_+, \RR^2) \to \DD(\RR_+, \RR^2)$ is given by
 \[
   (\Psi(f))(t)
   := f(t) +  \begin{bmatrix}
               b_1 \\
               b_2 \\
              \end{bmatrix} t , \qquad
   f \in \DD(\RR_+,\RR^2) , \qquad t \in \RR_+.
 \]
The mappings $\Psi^{(n)}$, $n\in\NN$, and $\Psi$ are measurable, which can be checked in the same way
 as in Step 4/(a) in Barczy et al.\ \cite{BarBezPap1} replacing $\DD(\RR_+,\RR)$ by $\DD(\RR_+,\RR^2)$
 in the argument given there.
One can also check that the set $C := \CC(\RR_+, \RR^2)$ satisfies $C \in \cB(\DD(\RR_+, \RR^2))$,
 $\PP( ( [\cM_{t,1}, \cM_{t,2}]^\top )_{t\in\RR_+} \in C) = 1$, and
 $\Psi^{(n)}(f^{(n)}) \to \Psi(f)$ in $\DD(\RR_+, \RR^2)$ as $n \to \infty$ if $f^{(n)} \to f$ in $\DD(\RR_+, \RR^2)$ as $n \to \infty$
 with $f \in C$, \ $f^{(n)}\in \DD(\RR_+, \RR^2)$, $n\in\NN$.
Namely, one can follow the same argument as in Step 4/(b) in Barczy et al.\ \cite{BarBezPap1} replacing
 $\DD(\RR_+,\RR)$ by $\DD(\RR_+,\RR^2)$, and $\CC(\RR_+,\RR)$ by $\CC(\RR_+,\RR^2)$, respectively,
 in the argument given there.
So we can apply Lemma \ref{Conv2Funct}, and we obtain
\begin{align*}
  \left( n^{-1}
         \begin{bmatrix}
         X_{\nt,1} \\
         X_{\nt,2}  \\
        \end{bmatrix}
       \right)_{t\in\RR_+}
  = \Psi^{(n)}\left(\left( \begin{bmatrix}
        \cM_{t,1}^{(n)}  \\
        \cM_{t,2}^{(n)}  \\
      \end{bmatrix}\right)_{t\in\RR_+} \right)
  \distr
  \Psi\left(
     \left(\begin{bmatrix}
        \cM_{t,1}\\
        \cM_{t,2} \\
      \end{bmatrix}\right)_{t\in\RR_+} \right)
  \qquad \text{as $n \to \infty$,}
 \end{align*}
 where
 \[
  \Psi\left(
     \left(\begin{bmatrix}
        \cM_{t,1} \\
        \cM_{t,2} \\
      \end{bmatrix}\right)_{t\in\RR_+} \right)
    = \left(\begin{bmatrix}
        \cM_{t,1} \\
        \cM_{t,2} \\
      \end{bmatrix}
      + \begin{bmatrix}
        b_1 \\
        b_2 \\
      \end{bmatrix}t
      \right)_{t\in\RR_+}
    \distre \left(
     \left(\begin{bmatrix}
        \cX_{t,1} \\
        \cX_{t,2} \\
      \end{bmatrix}\right)_{t\in\RR_+} \right),
 \]
 yielding \eqref{Conv_X_2_Thm2}.

{\sl Step 4/(a).}
Next, we turn to prove the main assertion. By \eqref{integrals}, we have
 \begin{align}\label{help22}
   n^{-2} X^{(2)}_{\nt,3}=\int_0^{ \frac{\nt}{n}} n^{-1}X_{\ns,1}^{(1)} \, \dd s, \qquad
    n^{-2} X^{(3)}_{\nt,3}=\int_0^{ \frac{\nt}{n}} n^{-1}X_{\ns,2}^{(2)} \, \dd s
 \end{align}
 for all $t \in \RR_+$ and $n \in \NN$.
Using \eqref{Conv_X_2_Thm2} and Lemmas \ref{ConvMatrixMult} and \ref{ConvIteratedIntegral}, we check that
 \begin{align}\label{Conv_I_2}
  \left(\begin{bmatrix}
          n^{-1} X_{\nt,1} \\
          n^{-1} X_{\nt,2} \\
          a_{3,1}n^{-2} X^{(2)}_{\nt,3} + a_{3,2}n^{-2} X^{(3)}_{\nt,3}
         \end{bmatrix}\right)_{t\in\RR_+}
  \distr
  \left(\begin{bmatrix}
          \cX_{t,1} \\
          \cX_{t,2} \\
          a_{3,1}\int_0^t\cX_{s,1}\, \dd s+a_{3,2}\int_0^t\cX_{s,2} \, \dd s
         \end{bmatrix}\right)_{t\in\RR_+}
 \end{align}
 as $n \to \infty$.
Namely, let us apply Lemma \ref{ConvIteratedIntegral} with $d=2$, and the functions $\Phi,\Phi_n$, $n\in\NN$, given there, and Lemma \ref{ConvMatrixMult} with
\begin{equation*}
\bB:=\begin{bmatrix}
1 & 0 & 0 & 0 & 0 & 0\\
0 & 1 & 0 & 0 & 0 & 0\\
0 & 0 & a_{3,1} & a_{3,2} & 0 & 0\\
\end{bmatrix}.
\end{equation*}
Then, by \eqref{help22}, we have
\begin{align}\label{help20}
 \begin{split}
	\begin{bmatrix}
          \cX_{t,1} \\
          \cX_{t,2} \\
          a_{3,1}\int_0^t\cX_{s,1}\, \dd s+a_{3,2}\int_0^t\cX_{s,2} \, \dd s
         \end{bmatrix}
  &=\bB \left(\Phi\left(
	\left(\begin{bmatrix}
	\cX_{u,1}\\
	\cX_{u,2}
	\end{bmatrix}\right)_{u\in\RR_+}\right)\right)(t)
	,\\[1mm]
	\begin{bmatrix}
          n^{-1} X_{\nt,1} \\
          n^{-1} X_{\nt,2} \\
          a_{3,1}n^{-2} X^{(2)}_{\nt,3} + a_{3,2}n^{-2} X^{(3)}_{\nt,3}
         \end{bmatrix}&
   = \bB \left(\Phi_n\left(
    \left(n^{-1}
	\begin{bmatrix}
	X_{\lfloor nu\rfloor,1}\\
	X_{\lfloor nu\rfloor,2}
	\end{bmatrix}\right)_{u\in\RR_+}\right)\right)(t)
 \end{split}
\end{align}
for $t\in\RR_+$ and $n\in\NN$.
Since $([\cX_{u,1},\cX_{u,2}]^\top)_{u\in\RR_+}$ has continuous sample paths almost surely, Lemma \ref{ConvIteratedIntegral} yields that
 \[
  \Phi_n\left(
    \left(n^{-1}
	\begin{bmatrix}
	X_{\lfloor nu\rfloor,1}\\
	X_{\lfloor nu\rfloor,2}
	\end{bmatrix}\right)_{u\in\RR_+}\right)
   \distr
   \Phi\left(
	\left(\begin{bmatrix}
	\cX_{u,1}\\
	\cX_{u,2}
	\end{bmatrix}\right)_{u\in\RR_+}\right) \qquad \text{as \ $n\to\infty$.}
 \]
Taking into account \eqref{help20}, an application of Lemma \ref{ConvMatrixMult} implies \eqref{Conv_I_2}, as desired.

{\sl Step 4/(b).}
Next, we shall show that
\begin{equation}\label{case2side0}
\sup_{t\in[0,T]}\left\Vert\begin{bmatrix}
0\\
0\\
n^{-2}X_{\nt,3}^{(4)}
\end{bmatrix}\right\Vert\stoch
0\qquad
\text{as \ $n \to \infty$ \ for all \ $T \in \RR_{++}$.}
\end{equation}
Note that for all $T\in\RR_+$ and $n\in\NN$, the left hand side of \eqref{case2side0} is equal to $n^{-2}\sup_{t\in[0,T]}\left|X_{\nt,3}^{(4)}\right|$.
By \eqref{Xdeco} and Lemma \ref{M_sum_growth}, for all $T\in\RR_+$ we have
\begin{equation*}
 \EE\left(n^{-4}\sup_{t\in[0,T]}|X_{\nt,3}^{(4)}|^2\right)
   =n^{-4}\EE\left(\sup_{t\in[0,T]}\left(\sum_{j=1}^{\nt}(M_{j,3}+b_3)\right)^2 \right)=\OO(n^{-4+\eta_3+1}),\qquad n\in\NN.
\end{equation*}
 By Lemma \ref{eta_values}, we have $\eta_3=2$, and hence
 $\OO(n^{-4+ \eta_3+1})\to 0$ as $n\to\infty$, yielding that
 $\left(n^{-2}\sup_{t\in[0,T]}|X_{\nt,3}^{(4)}|\right)^2\stoch0$ as $n\to\infty$ for all $T\in\RR_+$.
This implies \eqref{case2side0}.

{\sl Step 4/(c).}
Consequently, by Lemma \ref{Lem_JacShi} (a kind of Slutsky's lemma for stochastic processes with trajectories
 in $\DD(\RR_+, \RR^d)$), the decomposition \eqref{help_proof_Thm3_3} and the convergences \eqref{Conv_I_2} and \eqref{case2side0} yield that
 \[
 \left( \begin{bmatrix}
          n^{-1} X_{\nt,1} \\
          n^{-1} X_{\nt,2} \\
          n^{-2} X_{\nt,3}
         \end{bmatrix}  \right)_{t\in\RR_+}
   \distr \left(\begin{bmatrix}
             \cX_{t,1} \\
	   \cX_{t,2} \\
             \int_0^t a_{3,1}\cX_{s,1}+a_{3,2}\cX_{s,2} \, \dd s
            \end{bmatrix}\right)_{t\in\RR_+}
 \qquad \text{as \ $n \to \infty$.}
 \]
By It\^o's formula, the limit process  above is the pathwise unique strong solution of the SDE \eqref{SDE_X_2}
 with initial value $[\cX_{0,1},\cX_{0,2},\cX_{0,3}]^\top =\bzero$,
 thus we get the statement of Theorem \ref{main_2}.

\section{Proof of Theorem \ref{main_3}}\label{Proof3}

By \eqref{Xdeco}, we have the decomposition
 \begin{align}\label{help_proof_Thm3_4}
  \begin{bmatrix}
       n^{-1} X_{\nt,1} \\
       n^{-2} X_{\nt,2} \\
       n^{-2} X_{\nt,3}
      \end{bmatrix}
   = \begin{bmatrix}
      n^{-1} X_{\nt,1}^{(1)} \\[1mm]
      a_{2,1} \, n^{-2} X^{(1)}_{\nt,2} + n^{-2} X^{(2)}_{\nt,2} \\[1mm]
      a_{3,1} \, n^{-2} X_{\nt,3}^{(2)}  + n^{-2} X_{\nt,3}^{(4)}
     \end{bmatrix} , \qquad
   t \in \RR_+ , \quad n \in \NN ,
 \end{align}
 where
 \begin{equation*}
    X^{(2)}_{\nt,2}  = \sum_{j=1}^{\lfloor nt\rfloor} (M_{j,2} + b_2), \qquad X_{\nt,3}^{(4)} = \sum_{j=1}^{\nt} (M_{j,3} + b_3),
 \end{equation*}
 and, by \eqref{integrals},
 \begin{equation}\label{help21}
    X^{(1)}_{\nt,2} =  X^{(2)}_{\nt,3} = n^2 \int_0^{\frac{\nt}{n}} n^{-1}X_{\ns,1}^{(1)} \, \dd s,  \qquad t\in\RR_+, \quad n\in\NN.
 \end{equation}

{ \sl Step 1.}
By \eqref{help_Wei_Winnicki1}, we have
\begin{equation}\label{Conv_Base_3}
(n^{-1} X_{\nt,1})_{t\in\RR_+} \distr (\cX_{t,1})_{t\in\RR_+}\qquad \text{as $n \to \infty$,}
\end{equation}
where $(\cX_{t,1})_{t\in\RR_+}$ is the pathwise unique strong solution of
 the first equation of the SDE \eqref{SDE_X_3} with $\cX_{0,1}=0$.

{\sl Step 2/(a).}
We show that
\begin{equation}\label{Conv_I_3}
 \begin{split}
    &\left(\begin{bmatrix}
     n^{-1}X_{\nt,1}\\
     a_{2,1}n^{-2} X_{\nt,2}^{(2)}\\
     a_{3,1}n^{-2} X_{\nt,3}^{(4)}
    \end{bmatrix}\right)_{t\in\RR_+}
 \distr
 \left(\begin{bmatrix}
             \cX_{t,1} \\
	  a_{2,1} \int_0^t \cX_{s,1} \, \dd s \\
             a_{3,1}\int_0^t \cX_{s,1} \, \dd s
            \end{bmatrix}\right)_{t\in\RR_+} \qquad \text{ as $n\to\infty$.}
\end{split}
\end{equation}
Let us apply Lemma \ref{ConvIteratedIntegral} with $d=1$ and the functions $\Phi,\Phi_n$, $n\in\NN$, given there, and Lemma \ref{ConvMatrixMult} with
\begin{equation*}
\bB:=\begin{bmatrix}
1 & 0 & 0 \\
0 & a_{2,1} & 0 \\
0 & a_{3,1} & 0
\end{bmatrix}.
\end{equation*}
Then, by \eqref{help21}, we have
\begin{equation*}
	\begin{bmatrix}
          \cX_{t,1} \\
          a_{2,1}\int_0^t \cX_{s,1} \, \dd s \\
          a_{3,1}\int_0^t \cX_{s,1} \, \dd s
         \end{bmatrix}
  = \bB \left(\Phi((\cX_{u,1})_{u\in\RR_+})\right)(t),
	\qquad
	\begin{bmatrix}
          n^{-1} X_{\nt,1} \\
          a_{2,1} n^{-2} X^{(1)}_{\nt,2} \\
          a_{3,1} n^{-2} X^{(1)}_{\nt,2}
         \end{bmatrix}
	= \bB\left(\Phi_n( (n^{-1}X_{\lfloor nu \rfloor,1})_{u\in\RR_+} )\right)(t)
\end{equation*}
 for $t\in\RR_+$ and $n\in\NN$. Since $(\cX_{u,1})_{u\in\RR_+}$ has continuous sample paths almost surely,
  Lemma \ref{ConvIteratedIntegral} and \eqref{Conv_Base_3} imply that
  $\Phi_n( (n^{-1}X_{\lfloor nu\rfloor,1})_{u\in\RR_+} )\distr \Phi( (\cX_{u,1})_{u\in\RR_+} )$
  as $n\to\infty$.
Then an application of Lemma \ref{ConvMatrixMult} yields \eqref{Conv_I_3}, as desired.

{\sl Step 2/(b).}
The aim of the following discussion is to show that
\begin{equation}\label{case3side0}
\sup_{t\in[0,T]}\left\Vert\begin{bmatrix}
0\\
n^{-2}X_{\nt,2}^{(2)}\\
n^{-2}X_{\nt,3}^{(4)}
\end{bmatrix}\right\Vert\stoch
0\qquad
\text{as \ $n \to \infty$ \ for all \ $T \in \RR_{++}$.}
\end{equation}
Note that for all $T\in\RR_+$ and $n\in\NN$, we have
\begin{equation*}
\sup_{t\in[0,T]}\left\Vert\begin{bmatrix}
0\\
n^{-2}X_{\nt,2}^{(2)}\\
n^{-2}X_{\nt,3}^{(4)}
\end{bmatrix}\right\Vert
\leq  	n^{-2}\sup_{t\in[0,T]}|X_{\nt,2}^{(2)}|
	+n^{-2}\sup_{t\in[0,T]}|X_{\nt,3}^{(4)}|.
\end{equation*}
By \eqref{Xdeco} and Lemma \ref{M_sum_growth}, for all $T\in\RR_+$ we have
\begin{gather*}
\EE\left(n^{-4}\sup_{t\in[0,T]}|X_{\nt,2}^{(2)}|^2\right)=n^{-4}\EE\left(\sup_{t\in[0,T]}\left(\sum_{j=1}^{\nt}(M_{j,2}+b_2)\right)^2 \right)=\OO(n^{-4+\eta_2+1}),\qquad n\in\NN,\\
\EE\left(n^{-4}\sup_{t\in[0,T]}|X_{\nt,3}^{(4)}|^2\right)=n^{-4}\EE\left(\sup_{t\in[0,T]}\left(\sum_{j=1}^{\nt}(M_{j,3}+b_3)\right)^2 \right)=\OO(n^{-4+\eta_3+1}),\qquad n\in\NN.
\end{gather*}
By Lemma \ref{eta_values}, we have $\eta_2=\eta_3=2$, and hence $\OO(n^{-4+\eta_i+1})\to 0$ as $n\to\infty$, $i\in\{2,3\}$,
 yielding that
\begin{gather*}
\left(n^{-2}\sup_{t\in[0,T]}|X_{\nt,2}^{(2)}|\right)^2\stoch0\qquad \text{as $n\to\infty$ for all $T\in\RR_+$,}\\
\left(n^{-2}\sup_{t\in[0,T]}|X_{\nt,3}^{(4)}|\right)^2\stoch0\qquad \text{as $n\to\infty$ for all $T\in\RR_+$.}
\end{gather*}
This implies \eqref{case3side0}.

{\sl Step 2/(c).}
Consequently, by Lemma \ref{Lem_JacShi}, (a kind of Slutsky's lemma for stochastic processes with trajectories
 in $\DD(\RR_+, \RR^d)$), the decomposition \eqref{help_proof_Thm3_4} and the
 convergences \eqref{Conv_I_3} and \eqref{case3side0} yield that
 \[
    \left(\begin{bmatrix}
n^{-1} X_{\nt,1}\\
n^{-2} X_{\nt,2}\\
n^{-2} X_{\nt,3}
\end{bmatrix}\right)_{t\in\RR_+}
\distr
\left(\begin{bmatrix}
             \cX_{t,1} \\
             a_{2,1} \int_0^t \cX_{s,1} \, \dd s \\[1mm]
             a_{3,1} \int_0^t \cX_{s,1} \, \dd s
            \end{bmatrix}\right)_{t\in\RR_+}
\qquad\text{as \ $n \to \infty$.}
 \]
By It\^o's formula, the limit process above is the pathwise unique strong solution of
 the SDE \eqref{SDE_X_3} with initial value $[\cX_{0,1},\cX_{0,2},\cX_{0,3}]^\top = \bzero$,
 thus we get the statement of Theorem \ref{main_3}.

\section{Proof of Theorem \ref{main_4}}\label{Proof4}

By \eqref{Xdeco}, we have the decomposition
 \begin{align}\label{help_proof_Thm3_5}
   \begin{bmatrix}
       n^{-1} X_{\nt,1} \\
       n^{-2} X_{\nt,2} \\
       n^{-3} X_{\nt,3}
      \end{bmatrix}
   = \begin{bmatrix}
      n^{-1} X_{\nt,1} \\[1mm]
      a_{2,1} \, n^{-2} X^{(1)}_{\nt,2} + n^{-2} X^{(2)}_{\nt,2} \\[1mm]
      a_{3,2}a_{2,1}\, n^{-3} X_{\nt,3}^{(1)} + a_{3,1} \, n^{-3} X_{\nt,3}^{(2)} + a_{3,2} \, n^{-3} X_{\nt,3}^{(3)} + n^{-3} X_{\nt,3}^{(4)}
     \end{bmatrix}
 \end{align}
for all $t\in\RR_+$ and $n\in\NN$, where
\begin{gather*}
X^{(2)}_{\nt,2}  = \sum_{j=1}^{\lfloor nt\rfloor} (M_{j,2} + b_2), \qquad X_{\nt,3}^{(4)} = \sum_{j=1}^{\nt} (M_{j,3} + b_3),\\
X_{\nt,3}^{(2)} = \sum_{\ell=1}^\nt (\nt-\ell) (M_{\ell,1} + b_1),\qquad X_{\nt,3}^{(3)} = \sum_{\ell=1}^\nt (\nt-\ell) (M_{\ell,2} + b_2),
\end{gather*}
and, by \eqref{integrals},
\begin{gather}\label{help23}
X^{(1)}_{\nt,2} = {n^2}\int_0^{\frac{\nt}{n}} n^{-1} X_{\ns,1}^{(1)} \, \dd s,\qquad
X_{\nt,3}^{(1)}=  n^3\int_0^{\frac{\nt}{n}}\left(\int_0^{\frac{\nr}{n}}  n^{-1} X_{\ns,1}^{(1)}\,\dd s\right)\,\dd r.
\end{gather}

{\sl Step 1.}
By \eqref{help_Wei_Winnicki1}, we have
 \begin{equation}\label{Conv_Base_4}
 (n^{-1} X_{\nt,1})_{t\in\RR_+} \distr (\cX_{t,1})_{t\in\RR_+}\qquad \text{as $n \to \infty$,}
 \end{equation}
 where $(\cX_{t,1})_{t\in\RR_+}$ is the pathwise unique strong solution of
 the first equation of the SDE \eqref{SDE_X_5} with $\cX_{0,1}=0$.

{\sl Step 2./(a)}
 We show that
 \begin{align}\label{Conv_I_4}
  \begin{split}
   \left(\begin{bmatrix}
          n^{-1} X_{\nt,1} \\
          a_{2,1} n^{-2} X^{(1)}_{\nt,2} \\
          a_{3,2}a_{2,1} n^{-3} X^{(1)}_{\nt,3}
         \end{bmatrix}\right)_{t\in\RR_+}
  \distr
  \left(\begin{bmatrix}
          \cX_{t,1} \\
          a_{2,1}\int_0^t \cX_{s,1} \, \dd s \\
          a_{3,2}a_{2,1}\int_0^t\left(\int_0^r \cX_{s,1} \, \dd s\right)\,\dd r
         \end{bmatrix}\right)_{t\in\RR_+}\qquad \text{as $n \to \infty$.}
 \end{split}
 \end{align}
Let us apply Lemma \ref{ConvIteratedIntegral} with $d=1$ and the functions $\Phi,\Phi_n$, $n\in\NN$, given there, and Lemma \ref{ConvMatrixMult} with
\begin{equation*}
\bB:=\begin{bmatrix}
1 & 0 & 0 \\
0 & a_{2,1} & 0 \\
0 & 0 & a_{3,2}a_{2,1}
\end{bmatrix}.
\end{equation*}
Then, by \eqref{help23}, we have
\begin{equation*}
	\begin{bmatrix}
          \cX_{t,1} \\
          a_{2,1}\int_0^t \cX_{s,1} \, \dd s \\
          a_{3,2}a_{2,1}\int_0^t( \int_0^r\cX_{s,1} \, \dd s)\,\dd r
         \end{bmatrix}
        = \bB\left(\Phi((\cX_{u,1})_{u\in\RR_+})\right)(t)
\end{equation*}
and
\begin{equation*}
	\begin{bmatrix}
          n^{-1} X_{\nt,1} \\
          a_{2,1}n^{-2} X^{(1)}_{\nt,2} \\
          a_{3,2}a_{2,1}n^{-3} X^{(1)}_{\nt,3}
         \end{bmatrix}
      = \bB \left(\Phi_n( (n^{-1}X_{\lfloor nu \rfloor,1})_{u\in\RR_+}) \right)(t)
\end{equation*}
for $t\in\RR_+$ and $n\in\NN$.
 Since $(\cX_{u,1})_{u\in\RR_+}$ has continuous sample paths almost surely,
 Lemma \ref{ConvIteratedIntegral} and \eqref{Conv_Base_4} yield that
 $\Phi_n( (n^{-1}X_{\lfloor nu\rfloor,1})_{u\in\RR_+} )\distr \Phi((\cX_{u,1})_{u\in\RR_+})$ as $n\to\infty$.
Then an application of  Lemma \ref{ConvMatrixMult} implies \eqref{Conv_I_4}, as desired.

{\sl Step 2/(b).}
The aim of the following discussion is to show that
\begin{equation}\label{case4side0}
\sup_{t\in[0,T]}\left\Vert
\begin{bmatrix}
n^{-1}(X_{\nt,1}-X_{\nt,1}^{(1)})\\
n^{-2}(X_{\nt,2}-a_{2,1}X_{\nt,2}^{(1)})\\
n^{-3}(X_{\nt,3}-a_{3,2}a_{2,1}X_{\nt,3}^{(1)})
\end{bmatrix}
\right\Vert
\stoch0\qquad \text{as \ $n \to\infty$ for all $T\in\RR_{++}$.}
\end{equation}
Note that for all $T\in\RR_+$ and $n\in\NN$ we have
\begin{align*}
\sup_{t\in[0,T]}\left\Vert
\begin{bmatrix}
n^{-1}(X_{\nt,1}-X_{\nt,1}^{(1)})\\
n^{-2}(X_{\nt,2}-a_{2,1}X_{\nt,2}^{(1)})\\
n^{-3}(X_{\nt,3}-a_{3,2}a_{2,1}X_{\nt,3}^{(1)})
\end{bmatrix}
\right\Vert
=&\sup_{t\in[0,T]}\left\Vert\begin{bmatrix}
0\\
n^{-2}X_{\nt,2}^{(2)}\\
a_{3,1} \, n^{-3} X_{\nt,3}^{(2)} + a_{3,2} \, n^{-3} X_{\nt,3}^{(3)} + n^{-3} X_{\nt,3}^{(4)}
\end{bmatrix}\right\Vert\\
\leq&  	n^{-2}\sup_{t\in[0,T]}|X_{\nt,2}^{(2)}|
	+a_{3,1}n^{-3}\sup_{t\in[0,T]}|X_{\nt,3}^{(2)}|\\
	&+a_{3,2}n^{-3}\sup_{t\in[0,T]}|X_{\nt,3}^{(3)}|
	+n^{-3}\sup_{t\in[0,T]}|X_{\nt,3}^{(4)}|.
\end{align*}

By \eqref{Xdeco} and Lemma \ref{M_sum_growth}, for all $T\in\RR_+$ we have
\begin{gather*}
\EE\left(n^{-4}\sup_{t\in[0,T]}|X_{\nt,2}^{(2)}|^2\right)=n^{-4}\EE\left(\sup_{t\in[0,T]}\left(\sum_{j=1}^{\nt}(M_{j,2}+b_2)\right)^2 \right)=\OO(n^{-4+\eta_2+1}), \qquad n\in\NN,\\
\EE\left(n^{-6}\sup_{t\in[0,T]}|X_{\nt,3}^{(2)}|^2\right)=n^{-6}\EE\left(\sup_{t\in[0,T]}\left(\sum_{j=1}^{\nt}(\nt-j)(M_{j,1}+b_1)\right)^2 \right)=\OO(n^{-6+\eta_1+3}), \qquad n\in\NN,
\end{gather*}
\begin{gather*}
\EE\left(n^{-6}\sup_{t\in[0,T]}|X_{\nt,3}^{(3)}|^2\right)=n^{-6}\EE\left(\sup_{t\in[0,T]}\left(\sum_{j=1}^{\nt}(\nt-j)(M_{j,2}+b_2)\right)^2 \right)=\OO(n^{-6+\eta_2+3}), \qquad n\in\NN,\\
\EE\left(n^{-6}\sup_{t\in[0,T]}|X_{\nt,3}^{(4)}|^2\right)=n^{-6}\EE\left(\sup_{t\in[0,T]}\left(\sum_{j=1}^{\nt}(M_{j,3}+b_3)\right)^2 \right)=\OO(n^{-6+\eta_3+1}), \qquad n\in\NN.
\end{gather*}
By Lemma \ref{eta_values}, we have $\eta_1=1$, $\eta_2=2$, and $\eta_3=3$, and hence
 $\OO(n^{-4+\eta_2+1})\to 0$ as $n\to\infty$, $\OO(n^{-6+\eta_i+3})\to 0$ as $n\to\infty$, $i\in\{1,2\}$, and $\OO(n^{-6+\eta_3+1})\to 0$ as $n\to\infty$,
  yielding that
\begin{gather*}
\left(n^{-2}\sup_{t\in[0,T]}|X_{\nt,2}^{(2)}|\right)^2\stoch0\qquad \text{as $n\to\infty$ for all $T\in\RR_+$,}\\
\left(n^{-3}\sup_{t\in[0,T]}|X_{\nt,3}^{(2)}|\right)^2\stoch0\qquad \text{as $n\to\infty$ for all $T\in\RR_+$,}\\
\left(n^{-3}\sup_{t\in[0,T]}|X_{\nt,3}^{(3)}|\right)^2\stoch0\qquad \text{as $n\to\infty$ for all $T\in\RR_+$,}\\
\left(n^{-3}\sup_{t\in[0,T]}|X_{\nt,3}^{(4)}|\right)^2\stoch0\qquad \text{as $n\to\infty$ for all $T\in\RR_+$.}
\end{gather*}
This implies \eqref{case4side0}.

{\sl Step 2/(c).}
Consequently, by Lemma \ref{Lem_JacShi} (a kind of Slutsky's lemma for stochastic processes with trajectories
 in $\DD(\RR_+, \RR^d)$), the decomposition \eqref{help_proof_Thm3_5} and the convergences
 \eqref{Conv_I_4} and \eqref{case4side0} yield that%
 \[
   \left(\begin{bmatrix}
       n^{-1} X_{\nt,1} \\
       n^{-2} X_{\nt,2} \\
       n^{-3} X_{\nt,3}
      \end{bmatrix}\right)_{t\in\RR_+} \distr
\left(\begin{bmatrix}
             \cX_{t,1} \\
             a_{2,1} \int_0^t \cX_{s,1} \, \dd s \\[1mm]
             a_{3,2}a_{2,1} \int_0^t\left(\int_0^r \cX_{s,1} \, \dd s\right)\,\dd r
            \end{bmatrix}\right)_{t\in\RR_+}
 \qquad \text{as \ $n \to \infty$.}
 \]
By It\^o's formula, the limit process above is the pathwise unique strong solution of the SDE \eqref{SDE_X_5} with initial value
 $[\cX_{0,1},\cX_{0,2},\cX_{0,3}]^\top=\bzero$, thus we get the statement of Theorem \ref{main_4}.

\vspace*{5mm}

\appendix

\vspace*{5mm}

\noindent{\bf\Large Appendices}

\section{Some special sums}\label{appInt}

We present some auxiliary results for sums formed from the values of a function defined on $\ZZ_+$.
These results are used throughout the proofs.

\begin{Lem}\label{specSums}
Let $f:\ZZ_+\to\RR$.
Then, for each $n,k\in\NN$, we have
\begin{align}\label{specSums_formula1}
 &\sum_{\ell=0}^k f(\ell) = n \int_{0}^{\frac{k+1}{n}}f(\ns)\,\dd s,\\ \label{specSums_formula2}
 &\sum_{\ell=1}^k(k-\ell)f(\ell)=n\int_{0}^{\frac{k}{n}}\sum_{\ell=1}^\ns f(\ell)\,\dd s,
\end{align}
and
\begin{equation}\label{specSums_formula3}
 \sum_{\ell=1}^k\binom{k-\ell}{2}f(\ell)=n^2\int_{0}^{\frac{k}{n}}\left(\int_{0}^{\frac{\nr}{n}}\sum_{\ell=1}^\ns f(\ell)\,\dd s\right)\,\dd r.
\end{equation}
\end{Lem}

\noindent\textbf{Proof.}
Let $n,k\in\NN$ be fixed.
Then we have that
 \begin{align*}
 \sum_{\ell=0}^k f(\ell)
   & =n \sum_{\ell=0}^k f(\ell)\cdot \frac{1}{n}
    =n \sum_{\ell=0}^k f(\ell) \int_{\frac{\ell}{n}}^{\frac{\ell+1}{n}} 1\,\dd s\\
   &=n \sum_{\ell=0}^k \int_{\frac{\ell}{n}}^{\frac{\ell+1}{n}} f(\ns)\,\dd s
   = n\int_{0}^{\frac{k+1}{n}} f(\ns)\,\dd s,
 \end{align*}
 yielding \eqref{specSums_formula1}.

 Let $F(k):=\sum_{\ell=1}^kf(\ell)$, $k\in\ZZ_+$ (and recall the convention $F(0)=0$).
Using \eqref{specSums_formula1}, we have
 \begin{align*}
  \sum_{\ell=1}^k (k-\ell)f(\ell)
	&= \sum_{\ell=1}^k \sum_{\ell_1=1}^{k-\ell} f(\ell)
	= \sum_{\ell_1=1}^{k-1} \sum_{\ell=1}^{k-\ell_1} f(\ell)
	= \sum_{\ell_1=1}^{k-1} F(k-\ell_1)\\
	&= \sum_{\ell=1}^{k-1} F(\ell)
     = \sum_{\ell=0}^{k-1} F(\ell)
	=n\int_{0}^{\frac{k}{n}}F(\ns)\,\dd s,
 \end{align*}
 yielding \eqref{specSums_formula2}.

Furthermore, using \eqref{specSums_formula2} and \eqref{specSums_formula1}, we have that
\begin{align*}
   \sum_{\ell=1}^k \binom{k-\ell}{2}f(\ell)
	&	=\sum_{\ell=1}^k\left(\sum_{h=0}^{k-\ell-1}h\right)f(\ell)
		=\sum_{\ell=1}^k\left(\sum_{h=\ell}^{k-1}(h-\ell)\right)f(\ell)
		=\sum_{h=1}^{k-1}\sum_{\ell=1}^h(h-\ell)f(\ell)\\
	&	=n\sum_{h=0}^{k-1}\int_{0}^{\frac{h}{n}}F(\ns)\,\dd s
		=n^2\int_0^{\frac{k}{n}}\left(\int_{0}^{\frac{\nr}{n}}F(\ns)\,\dd s\right)\,\dd r,
\end{align*}
 yielding \eqref{specSums_formula3}.
\proofend

\section{Moments}\label{GWI_moments}

In this appendix, we collect some facts about the moments of 
 $\bX_k$, $\bM_k$, $k\in\ZZ_+$, and some related random variables.
\begin{Lem}\label{Moments}
Let $(\bX_k)_{k\in\ZZ_+}$ be a $p$-type GWI process such that $\bX_0 = \bzero$ and the moment condition \eqref{baseconds} holds.
Then for each $k \in \NN$, we have $\EE(\bX_k \mid \cF_{k-1}^\bX) = \bA \bX_{k-1} + \bb$
 and
 \begin{gather}
  \EE(\bX_k)
   = \sum_{j=0}^{k-1} \bA^j \bb, \label{mean2} \\
  \var\bigl(\bX_k \mid \cF_{k-1}^\bX\bigr)
    = \var\bigl(\bM_k \mid \cF_{k-1}^\bX\bigr)
	= \EE(\bM_k \bM_k^\top \mid \cF_{k-1}^\bX)
     = \bV^{(0)} + \sum_{i=1}^p X_{k-1,i} \bV^{(i)}, \label{condVa2r}\\
   \var\bigl(\bX_k\bigr)
      = \sum_{j=0}^{k-1} \bA^j \EE(\bM_{k-j} \bM_{k-j}^\top) (\bA^\top)^j, \label{Var2}\\
  \EE(\bM_k \bM_k^\top)
  = \bV^{(0)} + \sum_{i=1}^p \EE(X_{k-1,i}) \bV^{(i)}, \label{Cov2}
 \end{gather}
 where $\bA$, $\bb$ and $\bV^{(i)}$, $i\in\{0,1,\ldots,p\}$, are introduced in \eqref{help_jelolesek}.
\end{Lem}

Lemma \ref{Moments} is a special case of Lemma A.1 in Ispány and Pap \cite{IspPap2} for $p$-type GWI processes starting from $\bzero$.
For completeness, we note that Lemma A.1 in Ispány and Pap \cite{IspPap2} is stated only for critical
 $p$-type GWI processes, but its proof readily shows that it holds not only in the critical case.

Recall that, for a lower triangular matrix $\bA\in\RR_+^{p\times p}$ such that $a_{j,j}=1$, $j\in\{1,\dots,p\}$,
 we introduced the notations $\bC=\bA-\bI_p$ and
 \begin{equation}
\eta_i = \eta_i(\bA) =\max_{m\in\{1,\dots,p\}}\{m \mid \exists\, j\in\{1,\dots,p\}: c_{i,j}^{[m-1]}>0\},
 \qquad i\in\{1,\dots,p\},
 \end{equation}
 where $c_{i,j}^{[m-1]}$ denotes the $(i,j)$-th entry of $\bC^{m-1} =(\bA-\bI_p)^{m-1}$,
 see \eqref{etadef}.

\begin{Lem}\label{EEX1}%previously meanEstimate
Let $(\bX_k)_{k\in\ZZ_+}$ be a strongly critical $p$-type GWI process such that $\bX_0 = \bzero$,
 the moment condition \eqref{baseconds} holds, and suppose that the offspring mean matrix
 $\bA$ is lower triangular such that $a_{i,i}=1$, $i\in\{1,\dots,p\}$.
Then for each $i,j\in\{1,\dots,p\}$ we have
 \begin{gather*}
  \EE(X_{k,i})=\OO(k^{\eta_i}), \qquad |\EE(M_{k,i}M_{k,j})|=\OO(k^{\eta_i\wedge\eta_j}),\qquad k\in\NN.
 \end{gather*}
\end{Lem}

\noindent
\textbf{Proof.}
By \eqref{mean2} and \eqref{mPow}, using the notation $\bC=\bA-\bI_p$, for each $k\in\NN$ we have
\begin{align*}
\EE(\bX_k)&=\sum_{\ell=0}^{k-1}\bA^\ell\bb=\sum_{\ell=0}^{k-1}\sum_{m=0}^{p-1}\binom{\ell}{m}\bC^m\bb=\sum_{m=0}^{p-1}\bC^m\bb\sum_{\ell=0}^{k-1}\binom{\ell}{m}=\sum_{m=0}^{p-1}\bC^m\bb\binom{k}{m+1}.
\end{align*}
This implies that for each $i\in\{1,\dots,p\}$, we have
 \begin{align}\label{help17}
 \EE(X_{k,i})=\sum_{m=0}^{p-1}\left(\sum_{r=1}^pc_{i,r}^{[m]}b_r\right) \binom{k}{m+1}=\sum_{m=1}^{p}\left(\sum_{r=1}^pc_{i,r}^{[m-1]}b_r\right)\binom{k}{m},
    \qquad k\in\NN.
 \end{align}
We note that \eqref{help17} holds for $k=0$ as well, since $X_{0,i}=0$, $i\in\{1,\ldots,p\}$, and $\binom{0}{m}=0$, $m\in\{1,\ldots,p\}$.
Then $\EE(X_{k,i})$, $k\in\NN$, is a polynomial (in $k$) with constant term zero, so it is either the zero polynomial, or it has degree at least 1.
If it is the zero polynomial, then, since $\eta_i\geq1$, $i\in\{1,\dots,p\}$, we have $\EE(X_{k,i})=0=\OO(1)=\OO(k^{\eta_i})$, $k\in\NN$.

Otherwise, using the non-negativity of the coefficients $c^{[q]}_{\ell,k}$, $q\in\ZZ_+$, $\ell,k\in\{1,\ldots,p\}$, and $b_\ell$, $\ell\in\{1,\ldots,p\}$,
 the degree of the polynomial $\EE(X_{k,i})$, $k\in\NN$, equals the largest $m\in\{1,\ldots,p\}$ for which
 \begin{equation*}
 \sum_{r=1}^pc_{i,r}^{[m-1]}b_r>0.
 \end{equation*}
Due to the non-negativity of $c_{i,r}^{[m-1]}b_r$, $r\in\{1,\ldots,p\}$, such a sum is positive if and only if at least one of its terms is positive,
 thus it is enough to find the largest $m \in\{1,\ldots,p\}$ for which there exists some $j\in\{1,\ldots,p\}$ such that $c_{i,j}^{[m-1]}b_j>0$.
For this largest $m$, there exists a $j\in\{1,\dots,p\}$ such that $c_{i,j}^{[m-1]}b_j>0$,
 which implies that $c_{i,j}^{[m-1]}>0$, and thus $\eta_i\geq m$.
So $\eta_i$ is at least as large as the degree of the polynomial $\EE(X_{k,i})$, $k\in\NN$,
 yielding that $\EE(X_{k,i})=\OO(k^{\eta_i})$, $k\in\NN$ holds.

By \eqref{Cov2}, for each $i,j\in\{1,\dots,p\}$, we have
\begin{equation}\label{help18}
|\EE(M_{k,i}M_{k,j})|\leq|v_{i,j}^{(0)}|+\sum_{r=1}^p|v_{i,j}^{(r)}|\EE(X_{k-1,r}),\qquad k\in\NN.
\end{equation}
Next, we check that if $v_{i,j}^{(r)}\neq0$ for some $i,j,r\in\{1,\ldots,p\}$,
 then $\eta_r\leq\eta_i\wedge\eta_j$.
If $a_{i,r}=0$ for some $i,r\in\{1,\dots,p\}$, then $\xi_{1,1,r,i}\ase0$, and thus
 $v_{i,j}^{(r)}=\cov(0,\xi_{1,1,r,j})=0$, $j\in\{1,\ldots,p\}$.
Hence, if $v_{i,j}^{(r)}\ne 0$ for some $i,j,r\in\{1,\ldots,p\}$, then $a_{i,r}>0$,
 and we check that $\eta_i\geq \eta_r$ holds in this case.
Indeed, since $\bA$ is lower triangular, $a_{i,r}>0$ yields that $r\leq i$.
 If $r=i$, then $\eta_i=\eta_r\geq\eta_r$ trivially holds.
Otherwise, if $r<i$, then $c_{i,r}^{[1]}=a_{i,r}>0$ and thus $\eta_i\geq\eta_r+1>\eta_r$ by Lemma \ref{eta_ineq}.
A similar argument (replacing $i$ by $j$) shows that $\eta_j\geq \eta_r$.
This yields that $\eta_r\leq \eta_i\wedge\eta_j$ provided that $v_{i,j}^{(r)}\neq0$ for some $i,j,r\in\{1,\ldots,p\}$, as desired.
Since $a_{i,r}a_{j,r}=0$ and thus $v_{i,j}^{(r)}=0$ whenever $r>i\wedge j$,
 we may sum at the right-hand side of \eqref{help18} only up to $i\wedge j$, and,
 using that $\EE(X_{k,r})=\OO(k^{\eta_r})$, $k\in\NN$ (which we already proved in the present proof above),
  we have
 \begin{align*}
  |\EE(M_{k,i}M_{k,j})|
  &\leq |v_{i,j}^{(0)}|+\sum_{r=1}^{i\wedge j}|v_{i,j}^{(r)}|\EE(X_{k-1,r}) \\
  & = \OO(1) + \left( \sum_{r=1}^{i\wedge j}|v_{i,j}^{(r)}|+1 \right)
      \sum_{ \big\{r\in\{1,\ldots, i\wedge j\}:\,v_{i,j}^{(r)}\neq 0 \big\} }\OO(k^{\eta_r})
  = \OO(k^{\eta_i\wedge\eta_j}),
  \qquad k\in\NN.
 \end{align*}
\proofend

In the next remark, we compare Theorem 3 of Foster and Ney \cite{FosNey} and our Lemma \ref{EEX1}.

\begin{Rem}\label{FNmoment}
Foster and Ney \cite[Theorem 3]{FosNey} provides asymptotic growth rates for the moments of some critical decomposable GWI processes.
We recall their result specialized to our setup, that is,
 when the offspring mean matrix $\bA$ is lower triangular, $a_{i,i}=1$ for each $i\in\{1,\dots,p\}$, and  $a_{j+1,j}>0$ for each $j\in\{1,\dots,p-1\}$.
Let $(\bX_k)_{k\in\ZZ_+}$ be a strongly critical decomposable $p$-type GWI process starting from $\bzero\in\RR^p$
 and with an offspring mean matrix $\bA$ given above,
 and further suppose that all the moments of the coordinates of $\bxi_i$, $i\in\{1,\ldots,p\}$, and of $\bvare$ are finite.
Then we check that Theorem 3 of Foster and Ney \cite{FosNey} implies that
\begin{equation}\label{FosNeyMoment}
  \lim_{k\to\infty}k^{-i}\EE(X_{k,i}) =  b_1\frac{\prod_{\ell=1}^{i-1}a_{\ell+1,\ell}}{i!} \qquad \text{for each \ $i\in\{1,\dots,p\}$,} 
\end{equation}
 where $\prod_{\ell=1}^0:=1$.
Let $i\in\{1,\ldots,p\}$ be fixed.
Then, by Theorem 3 of Foster and Ney \cite{FosNey}, using their notations, with $\bw_i=1$ and $N=p$, we have
\begin{equation*}
  \lim_{k\to\infty}k^{-i} \EE(X_{k,i})
   =\lim_{k\to\infty}k^{-i} w_k(\be_i)
   = \sum_{r=1}^\infty \frac{1}{r!}A_r(\be_i),
\end{equation*}
where $\be_i:=[\delta_{i,j}]_{j\in\{1,\dots,p\}}\in\RR^p$ and for $\bzero<\bh\in\ZZ_+^p$,
\begin{equation*}
A_r(\bh)=\sum_{\substack{\bm_1+\dots+\bm_r=\bh\\\bm_j>\bzero,\forall j\in\{1,\dots,r\}}}\binom{\bh}{\bm_1,\dots,\bm_r}\prod_{j=1}^r\frac{\ba_1\cdot\bc(\bm_j)}{\langle\bm_j,\bN\rangle}, \qquad r\in\NN,
\end{equation*}
where $\bN:=[1,\dots,N]^\top\in\ZZ_+^N$,
\begin{equation*}
\bc(\bh)=\frac{\bu}{2}\sum_{\alpha=1}^{p_1}v_\alpha\sum_{i,j=1}^{p_1}q_{i,j}^{(\alpha)}\sum_{\bzero<\bm<\bh}\binom{\bh}{\bm}\frac{c_i(\bm)c_j(\bh-\bm)}{\langle\bh,\bN\rangle-1}
 \qquad \text{for $\bh\in\ZZ_+^N$ with $h_1+\cdots +h_N\geq 2$,}
\end{equation*}
 and
\begin{equation*}
\bc(\be_j)=\frac{\beta_j\bu(\bv_j\cdot\bw_j)}{(j-1)!}, \qquad j\in\{1,\dots,p\}.
\end{equation*}
In our case, for each $j\in\{1,\dots,p\}$, we have
\begin{gather*}
\ba_1=b_1,\quad \bu=1, \quad \bv_j=1,\quad \bw_j=\delta_{i,j},\quad N=p,
\end{gather*}
and therefore
\begin{equation*}
\bc(\be_j)=\frac{\beta_j \delta_{i,j}}{(j-1)!}, \qquad j\in\{1,\dots,p\},
\end{equation*}
and thus
\begin{equation*}
\sum_{r=1}^\infty \frac{1}{r!}A_r(\be_i)=A_1(\be_i)=\frac{b_1\beta_i}{\langle\be_i,\bN\rangle(i-1)!}=\frac{b_1\beta_i}{i!}.
\end{equation*}
By the definition of $\beta_j$, we have
\begin{equation*}
\beta_j=\prod_{\ell=1}^{j-1}(\bv_\ell\bA_{\ell+1,\ell})\cdot\bu_{\ell+1}=\prod_{\ell=1}^{j-1}(1 a_{\ell+1,\ell}) \cdot 1 = \prod_{\ell=1}^{j-1}a_{\ell+1,\ell}=c_{j,1}^{[j-1]},
 \qquad j\in\{ 1,\dots,p\},
\end{equation*}
 where we recall that $c_{j,1}^{[j-1]}$ denotes the $(j,1)$-th entry of $\bC^{j-1} = (\bA-\bI_p)^{j-1}$,
 and the last equality holds by part (ii) of Remark \ref{matrix_power_remark}.
Thus
\begin{equation*}
 \lim_{k\to\infty} k^{-i}\EE(X_{k,i})=\lim_{k\to\infty} k^{-i}w_k(\be_i)
    = b_1\frac{\prod_{\ell=1}^{i-1}a_{\ell+1,\ell}}{i!},
\end{equation*}
 yielding \eqref{FosNeyMoment}.
Note that \eqref{FosNeyMoment} can be rewritten as $\binom{k}{i}^{-1}\EE(X_{k,i})\to b_1c_{i,1}^{[i-1]}$ as $k\to\infty$.

On the other hand, using the proof of Lemma \ref{EEX1}, we check
 \begin{equation*}
 \lim_{k\to\infty} \binom{k}{i-r_i+1}^{-1} \EE(X_{k,i})
    =  b_{r_i}c_{i,r_i}^{[i-r_i]}
   \qquad\text{and}\qquad
    \lim_{k\to\infty}  k^{-(i-r_i+1)} \EE(X_{k,i})
     =b_{r_i}\frac{c_{i,r_i}^{[i-r_i]}}{(i-r_i+1)!},
 \end{equation*}
 where
  \begin{equation}\label{r_i_def}
   r_i:=\min\{r\in\{1,\dots,i\}:b_r>0\}\vee 1,\qquad i\in\{1,\dots,p\}.
\end{equation}
By \eqref{help17}, we have
 \begin{equation*}
 \EE(X_{k,i})=\sum_{m=1}^{p}\left(\sum_{r=1}^pc_{i,r}^{[m-1]}b_r\right)\binom{k}{m},
    \qquad k\in\NN,
 \end{equation*}
 where we recall the convention $\binom{k}{m}=0$ for $k,m\in\ZZ_+$ with $k<m$.
Using that $c_{i,r}^{[m-1]}=0$ if $m>i-r$ (see part (i) of Remark \ref{matrix_power_remark}),
we get
 \begin{equation*}
 \EE(X_{k,i})=\sum_{m=1}^{p}\sum_{r=1}^{i}c_{i,r}^{[m-1]}b_r\binom{k}{m},
    \qquad k\in\NN.
 \end{equation*}
By changing the order of the summations, we get
\begin{align*}
 \EE(X_{k,i})=&\sum_{r=1}^{i}\sum_{m=1}^{p}c_{i,r}^{[m-1]}b_r\binom{k}{m}=\sum_{r=1}^{i}b_r\sum_{m=1}^{i-r+1}c_{i,r}^{[m-1]}\binom{k}{m}\\
 =&\sum_{r=1}^ib_r\left(c_{i,r}^{[i-r]}\binom{k}{i-r+1}+\sum_{m=1}^{i-r}c_{i,r}^{[m-1]}\binom{k}{m}\right)\\
 =&\sum_{r=1}^ib_r\left(c_{i,r}^{[i-r]}\binom{k}{i-r+1}+\OO(k^{i-r})\right),
    \qquad k\in\NN.
\end{align*}
Since $c_{i,i}^{[0]}=1>0$ and $c_{i,r}^{[i-r]}=\prod_{j=r}^{i-1}a_{j+1,j}>0$ for $r\in\{1,\dots,i-1\}$
 (see part (ii) of Remark \ref{matrix_power_remark}), with $r_i$ defined in \eqref{r_i_def},
 we obtain that
\begin{equation}\label{help24}
 \EE(X_{k,i})= b_{r_i}c_{i,r_i}^{[i-r_i]}\binom{k}{i-r_i+1}+\OO(k^{i-r_i}),
    \qquad k\in\NN.
\end{equation}
Note that if $i\in\{1,\ldots,p\}$ and $b_r=0$, $r\in\{1,\ldots,p\}$, then $\EE(X_{k,i})=0$, $k\in\NN$, and $r_i=1$,
 and in this case \eqref{help24} reads as $\EE(X_{k,i}) = \OO(k^{i-1})$, $k\in\NN$, which holds trivially.
Consequently, if $b_1>0$, then $r_i=1$, $i\in\{1,\ldots,p\}$, and hence our Lemma \ref{EEX1} is in accordance with Theorem 3
 of Foster and Ney \cite{FosNey}.

If $p=3$, then the case considered in this remark is nothing else but the case \textup{(4)} of \eqref{tablazat_esetek}.
Note also that \eqref{FosNeyMoment} with $p=3$ is the same as \eqref{help_Xk_exp_4} in Proposition \ref{Pro_Xk_exp}.
\proofend
\end{Rem}

\begin{Lem}\label{M4moment}
Let $(\bX_k)_{k\in\ZZ_+}$ be a strongly critical $p$-type GWI process such that
  $\bX_0 = \bzero$, the moment condition \eqref{baseconds} holds, and suppose that the offspring mean matrix $\bA$ is lower triangular
  such that $a_{r,r}=1$, $r\in\{1,\dots,p\}$.
Suppose that $i\in\{1,\dots,p\}$ is such that $a_{i,r}=\delta_{i,r}$ for $r\in\{1,\dots,p\}$, $\EE(\xi_{1,1,i,i}^4)<\infty$, and $\EE(\vare_i^4)<\infty$.
Then we have
\begin{equation*}
\EE(M_{k,i}^4)=\OO(k^2),\qquad k\in\NN.
\end{equation*}
\end{Lem}

\noindent
\textbf{Proof.}
By the assumption $a_{i,r}=\delta_{i,r}$, $ r\in\{1,\ldots,p\}$,
 we have $\xi_{k,j,r,i}\ase0$ for each $k,j\in\NN$ and  $r\neq i$, $r\in\{1,\dots,p\}$, and
 thus, by \eqref{Mk1} and \eqref{MBPI(d)}, we have
 \begin{align*}
 M_{k,i}=X_{k,i}-X_{k-1,i}-b_i=\sum_{j=1}^{X_{k-1,i}}\xi_{k,j,i,i}+\vare_{k,i}-\sum_{j=1}^{X_{k-1,i}}1-b_i
        =\sum_{j=1}^{X_{k-1,i}}\left(\xi_{k,j,i,i}-1\right)+(\vare_{k,i}-b_i)
 \end{align*}
 for each $k\in\NN$.
From this, using the power mean inequality, we get
\begin{equation*}
\EE(M_{k,i}^4)\leq 2^3\EE\left[\left(\sum_{j=1}^{X_{k-1,i}}\left(\xi_{k,j,i,i}-1\right)\right)^4\right] +2^3\EE\left[\left(\vare_{k,i}-b_i\right)^4\right],
   \qquad k\in\NN.
\end{equation*}
Using that, for each $k\in\NN$, the random variables $\xi_{k,j,i,i}$, $j\in\NN$, are independent
 and $\EE(\xi_{k,j,i,i})=1$, $j\in\NN$, by the multinomial theorem, we obtain that
\begin{align*}
 \EE\left[\left(\sum_{j=1}^{X_{k-1,i}}\left(\xi_{k,j,i,i}-1\right)\right)^4\right]
  &=\EE\left[\EE\left[\left(\sum_{j=1}^{X_{k-1,i}}\left(\xi_{k,j,i,i}-1\right)\right)^4 \,\Big\vert\, \cF_{k-1}^\bX\right]\right]\\
  &=\EE\left[X_{k-1,i}\EE\left[\left(\xi_{1,1,i,i}-1\right)^4\right]\right]
     +\EE\left[6\binom{X_{k-1,i}}{2} (\EE\left[\left(\xi_{1,1,i,i}-1\right)^2\right])^2\right]\\
  &\leq  \OO(\EE(X_{k-1,i}) ) + \OO(\EE(X_{k-1,i}^2)), \qquad k\in\NN,
\end{align*}
 and thus arrive at
 \begin{equation*}
 \EE(M_{k,i}^4)=\OO(\EE(X_{k-1,i}))+\OO(\EE(X_{k-1,i}^2))+\OO(1), \qquad  k\in\NN.
 \end{equation*}

We turn to calculate the growth rate of $\EE(X_{k,i}^2)=\var(X_{k,i})+\EE(X_{k,i})^2$  in $k\in\NN$.
Note that, due to $a_{i,r}=\delta_{i,r}$, $r\in\{1,\ldots,p\}$,
 we have $c_{i,r}^{[m]}=0$  for $r\in\{1,\ldots,p\}$ and $m\in\NN$.
Indeed, in the considered case, the $i$-th row of the matrix $\bC$ is identically zero,
 and consequently the $i$-th row of $\bC^m$, $m\in\NN$, is identically zero as well.
Thus $\eta_i=1$ and, by Lemma \ref{matrixlemma}, we have $a_{i,r}^{[m]}=\delta_{i,r}$ for each $r\in\{1,\dots,p\}$ and $m\in\ZZ_+$.
By these facts, \eqref{Var2} and Lemma \ref{EEX1}, we have
 \begin{align*}
 \var(X_{k,i})
    & = \sum_{j=0}^{k-1} \Big( \bA^j \EE(\bM_{k-j} \bM_{k-j}^\top) (\bA^\top)^j \Big)_{i,i}
          = \sum_{j=0}^{k-1} \sum_{r=1}^p \sum_{q=1}^p a_{i,r}^{[j]} \EE(M_{k-j,r} M_{k-j,q}) a_{i,q}^{[j]} \\
     &=\sum_{j=0}^{k-1}\EE(M_{k-j,i}^2)= \sum_{j=0}^{k-1}\OO(k^{\eta_i}) =  \sum_{j=0}^{k-1}\OO(k)=\OO(k^2), \qquad k\in\NN.
 \end{align*}
Consequently, using again Lemma \ref{EEX1}, we get
 \begin{align*}
  \EE(M_{k,i}^4)&=\OO(\EE(X_{k-1,i}))+\OO(\var(X_{k-1,i}))+\OO(\EE(X_{k-1,i})^2)+\OO(1) \\
                & = \OO(k^{\eta_i}) + \OO(k^2) + \OO(k^{2\eta_i}) + \OO(1) =\OO(k^2),
                \qquad  k\in\NN,
 \end{align*}
 since $\eta_i=1$.
\proofend

\begin{Lem}\label{M_sum_growth}
Let $(\bX_k)_{k\in\ZZ_+}$ be a strongly critical $p$-type GWI process such that
 $\bX_0 = \bzero$, the moment condition \eqref{baseconds} holds, and suppose that the offspring mean matrix
 $\bA$ is lower triangular such that $a_{i,i}=1$, $i\in\{1,\dots,p\}$.
Then for $i\in\{1,\dots,p\}$ and $T\in\RR_{++}$, we have
\begin{gather}\label{doobBound1}
\EE\left(\sup_{t\in[0,T]}\left(\sum_{\ell=1}^{\nt}(M_{\ell,i}+b_i)\right)^2\right)=\OO(n^{\eta_i+1}),\qquad n\in\NN,\\ \label{doobBound2}
\EE\left(\sup_{t\in[0,T]}\left(\sum_{\ell=1}^{\nt}(\nt-\ell)(M_{\ell,i}+b_i)\right)^2\right)=\OO(n^{\eta_i+3}),\qquad n\in\NN.
\end{gather}
\end{Lem}

\noindent\textbf{Proof.} Let $T\in\RR_{++}$ and $i\in\{1,\ldots,p\}$ be fixed.
By the power mean inequality, we have
 \begin{equation*}
 \sup_{t\in[0,T]}\left(\sum_{\ell=1}^{\nt}(M_{\ell,i}+b_i)\right)^2\leq2\nT^2 b_i^2 + 2\sup_{t\in[0,T]}\left(\sum_{\ell=1}^\nt M_{\ell,i}\right)^2,
    \qquad n\in\NN.
 \end{equation*}
Using Doob's maximal inequality (see, e.g., Revuz and Yor \cite[Chapter II, Corollary (1.6)]{RevYor})
for the martingale $\left(\sum_{\ell=1}^k M_{\ell,i}\right)_{k\in\NN}$ (with the filtration \ $(\cF_k^\bX)_{k\in\NN}$), we obtain
\begin{equation*}
\EE\left(\sup_{t\in[0,T]}\left(\sum_{\ell=1}^\nt M_{\ell,i}\right)^2\right)
	=\EE\left(\max_{k\in\{0,\dots,\nT\}}\left(\sum_{\ell=1}^k M_{\ell,i}\right)^2\right)
	\leq4\EE\left(\left(\sum_{\ell=1}^\nT M_{\ell,i}\right)^2\right),\qquad n\in\NN.
\end{equation*}
 Since for each $\ell=1,\ldots,\nT-1$ and $j=\ell+1,\ldots,\nT$, $n\in\NN$, we have
 \begin{align*}
  \EE(M_{\ell,i}M_{j,i})
    & = \EE\big( \EE(M_{\ell,i}M_{j,i} \mid \cF_{j-1}^\bX) \big)
      = \EE\big( M_{\ell,i} \EE(M_{j,i} \mid \cF_{j-1}^\bX) \big)\\
    & = \EE\big( M_{\ell,i} \EE(X_{j,i}  - \EE(X_{j,i} \mid \cF_{j-1}^\bX ) \mid \cF_{j-1}^\bX) \big)
     = \EE(M_{\ell,i}\cdot 0)
     =0,
 \end{align*}
by Lemma \ref{EEX1}, we get that
\begin{align*}
\EE\left(\left(\sum_{\ell=1}^\nT M_{\ell,i}\right)^2\right)
	&=\EE\left( \sum_{\ell=1}^\nT M_{\ell,i}^2 + 2 \sum_{\ell=1}^{\nT-1} \sum_{j=\ell+1}^\nT M_{\ell,i} M_{j,i}\right)\\
	&=\sum_{\ell=1}^\nT\EE(M_{\ell,i}^2) + 2\sum_{\ell=1}^{\nT-1}\sum_{j=\ell+1}^\nT\EE(M_{\ell,i}M_{j,i})\\
	&=\sum_{\ell=1}^\nT\EE(M_{\ell,i}^2)=\sum_{\ell=1}^\nT\OO(\ell^{\eta_i})= \OO(n^{\eta_i+1}),
    \qquad n\in\NN.
\end{align*}
This together with $2\nT^2b_i^2=\OO(n^2)=\OO(n^{\eta_i+1})$, $n\in\NN$, since $\eta_i\geq1$, yield \eqref{doobBound1}.

 Next, we prove \eqref{doobBound2}.
By Lemma \ref{specSums}, we have
\begin{equation}\label{help19}
\sum_{\ell=1}^{\nt}(\nt-\ell)(M_{\ell,i}+b_i)=n\int_{0}^{\frac{\nt}{n}}\sum_{\ell=1}^{\ns}(M_{\ell,i}+b_i)\,\dd s, \qquad t\in\RR_+,\quad  n\in\NN.
\end{equation}
For each $n\in\NN$ and $t\in[0,n^{-1})$, we integrate on the set $\{0\}$ in the integral on the right hand side of \eqref{help19},
 so the integral is $0$ in this case.
Otherwise, for each $n\in\NN$ and all $t\in[n^{-1},\infty)$, we have $\nt>0$, and, by (the measure-theoretic) Jensen's inequality,
\begin{align*}
   \left(\int_{0}^{\frac{\nt}{n}}\sum_{\ell=1}^{\ns}(M_{\ell,i}+b_i)\,\dd s\right)^2
    &=\frac{\nt^2}{n^2}\left(\frac{n}{\nt}\int_{0}^{\frac{\nt}{n}}\sum_{\ell=1}^{\ns}(M_{\ell,i}+b_i)\,\dd s\right)^2\\
	&\leq \frac{\nt^2}{n^2}\frac{n}{\nt}\int_{0}^{\frac{\nt}{n}}\left(\sum_{\ell=1}^{\ns}(M_{\ell,i}+b_i)\right)^2\,\dd s\\
	&\leq t\int_{0}^{t}\sup_{ r\in[0,t]}\left(\sum_{\ell=1}^{ \lfloor nr\rfloor}(M_{\ell,i}+b_i)\right)^2\,\dd s
	\leq t^2\sup_{ r\in[0,t]}\left(\sum_{\ell=1}^{ \lfloor nr\rfloor}(M_{\ell,i}+b_i)\right)^2.
\end{align*}
 Consequently, by \eqref{help19} and \eqref{doobBound1}, we get that
\begin{align*}
 \EE\left(\sup_{t\in[0,T]}\left(\sum_{\ell=1}^{\nt}(\nt-\ell)(M_{\ell,i}+b_i)\right)^2\right)
     &\leq n^2T^2\EE\left(\sup_{t\in[0,T]}\left(\sum_{\ell=1}^{\nt}(M_{\ell,i}+b_i)\right)^2\right)\\
     & = n^2T^2 \OO(n^{\eta_i+1}) = \OO(n^{\eta_i+3}),\qquad n\in\NN.
\end{align*}
\proofend

\section{A version of the continuous mapping theorem}
\label{CMT}

In this appendix, we recall a version of the continuous mapping theorem
for $\RR^d$-valued stochastic processes with c\`adl\`ag paths,
 and then we present some auxiliary lemmas that make the application of
this continuous mapping theorem easier throughout the proofs.

A function $f : \RR_+ \to \RR^d$ is called c\`adl\`ag if it is right
 continuous with left limits.
Let $\DD(\RR_+, \RR^d)$ and $\CC(\RR_+, \RR^d)$ denote the space of
 all $\RR^d$-valued c\`adl\`ag and continuous functions on $\RR_+$, respectively.
Let $\cB(\DD(\RR_+, \RR^d))$ denote the Borel $\sigma$-algebra on
 $\DD(\RR_+, \RR^d)$ for the metric defined in Jacod and Shiryaev
 \cite[Chapter VI, (1.26)]{JacShi} (see also Billingsley \cite[Chapter 3, Section 16]{Bil} or the proof of Lemma \ref{ConvMatrixMult}).
With this metric $\DD(\RR_+, \RR^d)$ is a
 complete and separable metric space and the topology induced by this metric is
 the so-called Skorokhod $J_1$ topology.
Recall that  $\skor{d}$ denotes the convergence in this Skorokhod
 topology for $\DD(\RR_+,\RR^d)$ (see Section \ref{Section_multi-type_GWI}).
Note that $\CC(\RR_+,\RR^d)\in\cB(\DD(\RR_+,\RR^d))$, see, e.g., Ethier and Kurtz \cite[Problem 3.11.25]{EthKur}.
For $\RR^d$-valued stochastic processes $(\bcY_t)_{t \in \RR_+}$ and
 $(\bcY^{(n)}_t)_{t \in \RR_+}$, $n \in \NN$, with c\`adl\`ag paths we write
 $\bcY^{(n)} \distr \bcY$ if the distribution of $\bcY^{(n)}$ on the
 space $(\DD(\RR_+,  \RR^d), \cB(\DD(\RR_+, \RR^d)))$ converges weakly to the
 distribution of $\bcY$ on the space
 $(\DD(\RR_+,  \RR^d), \cB(\DD(\RR_+, \RR^d)))$ as $n \to \infty$.
If $\xi$ and $\xi_n$, $n \in \NN$, are random elements with values in a metric space $(E,d)$,
 then we denote by $\xi_n \distr \xi$ the weak convergence of the
 distribution of $\xi_n$ on the space $(E, \cB(E))$ towards the
 distribution of $\xi$ on the space $(E, \cB(E))$ as $n \to \infty$,
 where $\cB(E)$ denotes the Borel $\sigma$-algebra on $E$ induced by
 the given metric $d$.

The following version of the continuous mapping theorem can be found for
 example in Theorem 3.27 of Kallenberg \cite{Kal}.

\begin{Lem}\label{lemma:kallenberg}
Let $(S, d_S)$ and $(T, d_T)$ be metric spaces and
 $(\xi_n)_{n \in \NN}$, $\xi$ be random elements with values in $S$
 such that $\xi_n \distr \xi$ as $n \to \infty$.
Let $f : S \to T$ and $f_n : S \to T$, $n \in \NN$, be measurable
 mappings and $C \in \cB(S)$ such that $\PP(\xi \in C) = 1$ and
 $\lim_{n \to \infty} d_T(f_n(s_n), f(s)) = 0$ if
 $\lim_{n \to \infty} d_S(s_n,s) = 0$ and $s \in C$, $s_n\in S$, $n\in\NN$.
Then $f_n(\xi_n) \distr f(\xi)$ as $n \to \infty$.
\end{Lem}

For functions $f$ and $f_n$, $n \in \NN$, in $\DD(\RR_+, \RR^d)$,
 we write $f_n \lu f$ if $(f_n)_{n \in \NN}$ converges to $f$
 locally uniformly, i.e., if $\sup_{t \in [0,T]} \|f_n(t) - f(t)\| \to 0$
 as $n \to \infty$ for all $T \in \RR_{++}$.
For measurable mappings $\Phi : \DD(\RR_+, \RR^d) \to \DD(\RR_+, \RR^q)$
 and $\Phi_n : \DD(\RR_+, \RR^d) \to \DD(\RR_+, \RR^q)$, $n \in \NN$,
 we will denote by $C_{\Phi, (\Phi_n)_{n \in \NN}}$ the set of all functions
 $f \in \CC(\RR_+, \RR^d)$ \ for which \ $\Phi_n(f_n) \lu \Phi(f)$ whenever
 $f_n \lu f$ \ with \ $f_n \in \DD(\RR_+, \RR^d)$, $n \in \NN$.

One can formulate the following consequence of Lemma \ref{lemma:kallenberg}.

\begin{Lem}\label{Conv2Funct}
Let $d,\,q\in\NN$.
Let $(\bcU_t)_{t \in \RR_+}$ and $(\bcU^{(n)}_t)_{t \in \RR_+}$, $n \in \NN$,
 be $\RR^d$-valued stochastic processes with c\`adl\`ag paths such that
 $\bcU^{(n)} \distr \bcU$ as $n\to\infty$.
Let $\Phi : \DD(\RR_+, \RR^d) \to \DD(\RR_+, \RR^q)$ and
  $\Phi_n : \DD(\RR_+, \RR^d) \to \DD(\RR_+, \RR^q)$, $n \in \NN$, be
 measurable mappings such that there exists
 $C \subset C_{\Phi,(\Phi_n)_{n\in\NN}}$ with $C \in \cB(\DD(\RR_+, \RR^d))$ and $\PP(\bcU \in C) = 1$.
Then $\Phi_n(\bcU^{(n)}) \distr \Phi(\bcU)$ as $n\to\infty$.
\end{Lem}

 Next, we present some auxiliary lemmas which are used throughout the proofs.
The forthcoming Lemma \ref{ConvMatrixMult} may be known in the literature, but we could not find it anywhere,
 so, for completeness, we present a proof as well.

\begin{Lem}\label{ConvMatrixMult}
Let $d,q\in\NN$, $\bB\in\RR^{q\times d}$, and let $\Psi_\bB:\DD(\RR_+,\RR^d)\to\DD(\RR_+,\RR^q)$,
 $(\Psi_\bB(f))(t):=\bB f(t)$ for $f\in\DD(\RR_+,\RR^d)$ and $t\in\RR_+$.
Then $ \Psi_\bB$ is continuous  (in particular, measurable).
Further, if $(\bcU_t)_{t \in \RR_+}$ and $(\bcU^{(n)}_t)_{t \in \RR_+}$, $n \in \NN$,
 are $\RR^d$-valued stochastic processes with c\`adl\`ag paths such that
 $\bcU^{(n)} \distr \bcU$ as $n\to\infty$, then $\bB\bcU^{(n)} \distr \bB\bcU$ as $n\to\infty$.
\end{Lem}

\noindent\textbf{Proof.}
 We will use the notations of Billingsley \cite[Chapter 3, Section 16]{Bil}.
Recall that the Skorokhod metric on $\DD(\RR_+,\RR^d)$ is given by
 \begin{equation*}
  d_\infty^\circ(f,g):=\sum_{m=1}^\infty2^{-m}\left(1\wedge d_m^\circ( f^{(m)},g^{(m)})\right),\qquad  f,g\in \DD(\RR_+,\RR^d),
 \end{equation*}
 where for $f\in\DD(\RR_+,\RR^d)$ and $m\in\NN$, define $f^{(m)}(t):=(1\wedge(m-t))^+f(t)$, $t\in\RR_+$, and
 \begin{equation*}
  d_m^\circ(f,g):=\inf_{\lambda\in\Lambda_m}\Big\{\Vert\lambda\Vert_m^\circ\vee\sup_{t\in[0,m]}\Vert f(t)-g(\lambda(t))\Vert\Big\},
  \qquad f,g\in\DD(\RR_+,\RR^d).
 \end{equation*}
 Here, for each $m\in\NN$, $\Lambda_m$ denotes the class of strictly increasing continuous functions of $[0,m]$ onto itself,
 and
  \[
    \Vert\lambda\Vert_m^\circ:= \sup_{0\leq s<t\leq m} \left\vert \log\left(\frac{\lambda(t) - \lambda(s)}{t-s} \right)\right\vert,
         \qquad \lambda\in\Lambda_m.
  \]
For each $m\in\NN$, we have
 \begin{align*}
 d_m^\circ(\Psi_\bB(f),\Psi_\bB(g))
 	&=\inf_{\lambda\in\Lambda_m}\Big\{\Vert\lambda\Vert_m^\circ\vee\sup_{t\in[0,m]}\Vert \bB f(t)-\bB g(\lambda(t))\Vert\Big\}\\
	&\leq\inf_{\lambda\in\Lambda_m}\Big\{\Vert\lambda\Vert_m^\circ\vee\Vert\bB\Vert\sup_{t\in[0,m]}\Vert f(t)-g(\lambda(t))\Vert\Big\}\\
	&\leq(1+\Vert\bB\Vert)d_m^\circ(f,g),\qquad f,g\in\DD(\RR_+,\RR^d).
 \end{align*}
Note that
 \begin{align*}
   (\Psi_\bB(f))^{(m)}(t)
    & = (\bB f)^{(m)}(t)=(1\wedge(m-t))^+ (\bB f)(t)
    =\bB\left((1\wedge(m-t))^+f(t)\right)\\
    &=\bB(f^{(m)}(t)) =(\bB(f^{(m)}))(t)=(\Psi_\bB(f^{(m)}))(t), \qquad t\in\RR_+, \; m\in\NN.
 \end{align*}
Thus we have
 \begin{align*}
  d_\infty^\circ(\Psi_\bB(f),\Psi_\bB(g))
	&=\sum_{m=1}^\infty2^{-m}\left(1\wedge d_m^\circ( (\Psi_\bB(f))^{(m)},(\Psi_\bB(g))^{(m)})\right)\\
	&=\sum_{m=1}^\infty2^{-m}\left(1\wedge d_m^\circ(\Psi_\bB(f^{(m)}),\Psi_\bB(g^{(m)}))\right)\\
	&\leq\sum_{m=1}^\infty2^{-m}\left(1\wedge\big((1+\Vert\bB\Vert)d_m^\circ(f^{(m)},g^{(m)}) \big)\right)\\
	&\leq(1+\Vert\bB\Vert)d_\infty^\circ(f,g), \qquad f,g\in\DD(\RR_+,\RR^d).
 \end{align*}
 This shows that $\Psi_\bB$ is (Lipschitz) continuous.

Let $(\bcU_t)_{t \in \RR_+}$ and $(\bcU^{(n)}_t)_{t \in \RR_+}$, $n \in \NN$,
 be $\RR^d$-valued stochastic processes with c\`adl\`ag paths such that
 $\bcU^{(n)} \distr \bcU$ as $n\to\infty$.
Since $\Psi_\bB$ is continuous, the continuous mapping theorem (see, e.g., Billingsley \cite[Theorem 2.7]{Bil}) 
 yields that $\bB\bcU^{(n)} = \Psi_\bB(\bcU^{(n)}) \distr \Psi_\bB(\bcU)= \bB\bcU$ as $n\to\infty$.
We note that, in the special case when $\bcU$ has continuous sample paths,
 Lemma \ref{Conv2Funct} also implies that $\bB\bcU^{(n)} \distr \bB\bcU$ as $n\to\infty$.
Indeed, using the notations of Lemma \ref{Conv2Funct}, let $\Phi:=\Psi_\bB$ $\Phi_n:=\Psi_\bB$, $n\in\NN$.
Then $C_{\Phi, (\Phi_n)_{n \in \NN}} = \CC(\RR_+,\RR^d)$, since for all $T>0$ and $f,f_n\in\DD(\RR_+,\RR^d)$, $n\in\NN$,
 we have
 \[
  \sup_{t\in[0,T]} \Vert \bB f_n(t) - \bB f(t)\Vert \leq \Vert \bB\Vert \sup_{t\in[0,T]} \Vert f_n(t) - f(t)\Vert.
 \]
Since $ \CC(\RR_+,\RR^d)\in\cB(\DD(\RR_+, \RR^d))$ (see, e.g., Ethier
and Kurtz \cite[Problem 3.11.25]{EthKur}) and $\PP(\cU\in \CC(\RR_+,\RR^d))=1$, one can apply Lemma \ref{Conv2Funct}.
\proofend

\begin{Lem}\label{coordMeasurable}
Let $d,p,q\in\NN$, $\Phi_1:\DD(\RR_+,\RR^d)\to\DD(\RR_+,\RR^p)$ and $\Phi_2:\DD(\RR_+,\RR^d)\to\DD(\RR_+,\RR^q)$ be measurable mappings.
Then $\Phi:\DD(\RR_+,\RR^d)\to\DD(\RR_+,\RR^{p+q})$, $\Phi(f):=[\Phi_1(f),\Phi_2(f)]^\top$, $f\in\DD(\RR_+,\RR^d)$ is measurable as well.
\end{Lem}
\noindent\textbf{Proof.} By Jacod and Shiryaev \cite[Chapter VI, Theorem 1.14/(c) and 1.21]{JacShi},
we have that $\mathcal{B}(\DD(\RR_+,\RR^d))=(\mathcal{B}(\DD(\RR_+,\RR)))^{\otimes d}$, $d\in\NN$,
where $(\mathcal{B}(\DD(\RR_+,\RR)))^{\otimes d}$ denotes
the product $\sigma$-algebra $\sigma(\{B_1\times\dots\times B_d: B_i\in\mathcal{B}(\DD(\RR_+,\RR)), i=1,\dots,d\})$.
Since
\begin{align*}
\mathcal{B}(\DD(\RR_+,\RR^{p+q}))
	&=(\mathcal{B}(\DD(\RR_+,\RR)))^{\otimes (p+q)}
	=(\mathcal{B}(\DD(\RR_+,\RR)))^{\otimes p}\otimes(\mathcal{B}(\DD(\RR_+,\RR)))^{\otimes q}\\
	&=\mathcal{B}(\DD(\RR_+,\RR^p))\otimes\mathcal{B}(\DD(\RR_+,\RR^q)),
\end{align*}
it is enough to check that $\Phi^{-1}(B_1\times B_2)\in\mathcal{B}(\DD(\RR_+,\RR^d))$
for all $B_1\in\mathcal{B}(\DD(\RR_+,\RR^p))$ and $B_2\in\mathcal{B}(\DD(\RR_+,\RR^q))$.
Indeed, it is sufficient to check that the preimages of those sets in $\mathcal{B}(\DD(\RR_+,\RR^{p+q}))$
 which generate $\mathcal{B}(\DD(\RR_+,\RR^{p+q}))$ belong to $\mathcal{B}(\DD(\RR_+,\RR^d))$,
 see, e.g., Klenke \cite[Theorem 1.81]{Kle}.
For all $B_1\in\mathcal{B}(\DD(\RR_+,\RR^{p}))$ and $B_2\in\mathcal{B}(\DD(\RR_+,\RR^{q}))$, we have that
\begin{equation*}
\Phi^{-1}(B_1\times B_2)=\Phi_1^{-1}(B_1)\cap\Phi_2^{-1}(B_2).
\end{equation*}
Since $\Phi_1,\Phi_2$ are measurable, we get $\Phi_i^{-1}(B_i)\in\mathcal{B}(\DD(\RR_+,\RR^{d}))$, $i=1,2$,
yielding that $\Phi_1^{-1}(B_1)\cap\Phi_2^{-1}(B_2)\in\mathcal{B}(\DD(\RR_+,\RR^{d}))$ (due to the fact that $\mathcal{B}(\DD(\RR_+,\RR^{d}))$ is a $\sigma$-algebra).
\proofend

The following lemma is a corollary of a result in Jacod and Shiryaev \cite[Chapter VI, Proposition 2.2 (b)]{JacShi}.

\begin{Lem}\label{CoordConv}
Let $p,q\in\NN$, $f\in\DD(\RR_+,\RR^p)$, $f_n\in\DD(\RR_+,\RR^p)$, $n\in\NN$,
and $g\in\DD(\RR_+,\RR^q)$, $g_n\in\DD(\RR_+,\RR^q)$, $n\in\NN$, be
such that $f_n\skor{p} f$ as $n\to\infty$ and $g_n\skor{q} g$ as $n\to\infty$.
If $g\in\CC(\RR_+,\RR^q)$, $g_n\in\CC(\RR_+,\RR^q)$, $n\in\NN$,
then we have
\begin{equation*}
\begin{bmatrix}
f_n\\
g_n
\end{bmatrix}\skor{p+q}
\begin{bmatrix}
f\\
g
\end{bmatrix}\qquad \text{as $n\to\infty$.}
\end{equation*}
\end{Lem}

\noindent\textbf{Proof.} By Jacod and Shiryaev \cite[Chapter VI, Proposition 2.1 (a)]{JacShi},
for all $t\in\RR_{++}$ there exists a sequence $(t_n)_{n\in\NN}$ in $\RR_{++}$ such that
$t_n\to t$ as $n\to\infty$ and $\Delta f_n(t_n)\to\Delta f(t)$ as $n\to\infty$.
Since the functions $g$, $g_n$, $n\in\NN$ are continuous, we have
$\bzero=\Delta g_n(t_n)\to\Delta g(t)=\bzero$ as $n\to\infty$.
By Jacod and Shiryaev \cite[Chapter VI, Proposition 2.2 (b)]{JacShi}
this implies the assertion.
\proofend

\begin{Cor}\label{ContMapVector}
Let $d,p,q\in\NN$, $\Phi_1:\DD(\RR_+,\RR^d)\to\DD(\RR_+,\RR^p)$ and $\Phi_2:\DD(\RR_+,\RR^d)\to\DD(\RR_+,\RR^q)$ be continuous mappings.
If $\Phi_1(\DD(\RR_+,\RR^d))\subset\CC(\RR_+,\RR^p)$ and $\Phi_2(\DD(\RR_+,\RR^d))\subset\CC(\RR_+,\RR^q)$,
then $\Phi:\DD(\RR_+,\RR^d)\to\DD(\RR_+,\RR^{p+q})$, $\Phi(f):=[\Phi_1(f),\Phi_2(f)]^\top$, $f\in\DD(\RR_+,\RR^d)$ is continuous as well, and $\Phi(\DD(\RR_+,\RR^d))\subset\CC(\RR_+,\RR^{p+q})$.
\end{Cor}

\noindent\textbf{Proof.}
 The fact that $\Phi$ maps $\DD(\RR_+,\RR^d)$ into $\CC(\RR_+,\RR^{p+q})$ is obvious.
Let $f\in\DD(\RR_+,\RR^d)$, $f_n\in\DD(\RR_+,\RR^d)$, $n\in\NN$, be such that $f_n\skor{d} f$ as $n\to\infty$.
By the continuity of $\Phi_1$ and $\Phi_2$, we have $\Phi_1(f_n)\skor{p}\Phi_1(f)$ as $n\to\infty$ and $\Phi_2(f_n)\skor{q}\Phi_2(f)$ as $n\to\infty$.
By the assumption, the functions $\Phi_1(f), \Phi_2(f)$, and $\Phi_1(f_n), \Phi_2(f_n)$, $n\in\NN$, are continuous,
and thus, by Lemma \ref{CoordConv}, we have $\Phi(f_n)\skor{p+q}\Phi(f)$ as $n\to\infty$, meaning $\Phi$ is a 
 sequentially continuous map.
 Hence, using that $\DD(\RR_+,\RR^d)$ and $\DD(\RR_+,\RR^{p+q})$ are metric spaces, it yields that $\Phi$ is continuous.
\proofend

\begin{Lem}\label{IntegralCont}
For each $d\in\NN$, the map $\varphi:\DD(\RR_+,\RR^d)\to\DD(\RR_+,\RR^d)$
defined by
\begin{equation*}
(\varphi(f))(t):=\int_0^tf(s)\,\dd s,\qquad f\in\DD(\RR_+,\RR^d),\; t\in\RR_+,
\end{equation*}
is continuous, and maps $\DD(\RR_+,\RR^d)$ into $\CC(\RR_+,\RR^d)$.
\end{Lem}

\noindent\textbf{Proof.}
The fact that $\varphi$ maps $\DD(\RR_+,\RR^d)$ into $\CC(\RR_+,\RR^d)$ is obvious.
For $i\in\{1,\dots,d\}$, let $\bB_i:=[\delta_{1,i},\dots,\delta_{d,i}]\in\RR^{1\times d}$.
Let $\widetilde\varphi:\DD(\RR_+,\RR)\to\DD(\RR_+,\RR)$, $(\widetilde\varphi(g))(t):=\int_0^tg(s)\,\dd s$, $g\in\DD(\RR_+,\RR)$, $t\in\RR_+$.
For each $i\in\{1,\dots,d\}$, the mapping $\varphi_i:\DD(\RR_+,\RR^d)\to\DD(\RR_+,\RR)$,
 $(\varphi_i(h))(t):=(\widetilde\varphi(\bB_ih))(t)=\int_0^t \bB_ih(s)\,\dd s$, $h\in\DD(\RR_+,\RR^d)$, $t\in\RR_+$,
 is continuous, since $\widetilde\varphi$ is continuous (see, e.g., Ethier and Kurtz \cite[Problem 3.11.26]{EthKur}),
 the multiplication by the matrix $\bB_i$ is continuous (due to Lemma \ref{ConvMatrixMult}),
 and thus their composition is continuous as well.
Obviously, $\varphi_i$ maps $\DD(\RR_+,\RR^d)$ into $\CC(\RR_+,\RR)$ for each $i\in\{1,\ldots,d\}$.

Then we have
\begin{equation*}
(\varphi(f))(t)
=\begin{bmatrix}
\bB_1\int_0^tf(s)\,\dd s\\
\vdots\\
\bB_d\int_0^tf(s)\,\dd s
\end{bmatrix}
=\begin{bmatrix}
\int_0^t\bB_1f(s)\,\dd s\\
\vdots\\
\int_0^t\bB_df(s)\,\dd s
\end{bmatrix}
=\begin{bmatrix}
(\varphi_1(f))(t)\\
\vdots\\
(\varphi_d(f))(t)
\end{bmatrix},\qquad  f\in\DD(\RR_+,\RR^d),\;\; t\in\RR_+.
\end{equation*}
Thus we get $\varphi(f)=[\varphi_1(f),\dots,\varphi_d(f)]^\top$, $f\in\DD(\RR_+,\RR^d)$,
and by applying Corollary \ref{ContMapVector} $d-1$ times, we have that $\varphi$ is continuous.

\noindent\textbf{ Second proof.}
This proof is a more direct one.
Let $\Lambda'$ be the set of strictly monotone increasing, Lipschitz continuous functions $\lambda$ defined on $\RR_+$
such that $\lambda(0)=0$ and $\lim_{t\to\infty}\lambda(t)=\infty$. It is well known that for $\lambda\in\Lambda'$,
the derivative $\lambda'$ exists Lebesgue almost everywhere, see, e.g., Wang \cite[Theorem 3.5.6 and Exercise 3.25]{Wan}.
For $\lambda\in\Lambda'$, let
 \begin{equation*}
 \gamma(\lambda):=\sup_{0\leq s<t}\left|\ln\frac{\lambda(t)-\lambda(s)}{t-s}\right|=\esssup_{t\in\RR_+}\left|\ln\lambda'(t)\right|,
 \end{equation*}
 and let $\Lambda:=\{\lambda \in\Lambda' : \gamma(\lambda)<\infty\}$, see Ethier and Kurtz \cite[page 117]{EthKur}.
It is well-known that
 \begin{equation*}
 \esssup_{t\in\RR_+}|\lambda'(t)-1|\leq1-e^{-\gamma(\lambda)}\leq\gamma(\lambda),\qquad \lambda\in\Lambda,
 \end{equation*}
 see Ethier and Kurtz \cite[formula (5.5), page 117]{EthKur}.
Note that if $\lambda\in\Lambda$, then $\lambda$ is injective and absolutely continuous on each compact interval of $\RR_+$
 (following from the Lipschitz continuity of $\lambda$, see, e.g., Wang \cite[Exercise 3.25]{Wan}).
Thus one can apply integration by substitution for Lebesgue integrals and we get that
 \begin{equation}\label{int_hely}
 \int_0^{\lambda(t)}f(s)\,\dd s=\int_0^tf(\lambda(s))\lambda'(s)\,\dd s,\qquad f\in\DD(\RR_+,\RR^d),\quad t\in\RR_{+}.
 \end{equation}
Let $f\in\DD(\RR_+,\RR^d)$, $f_n\in\DD(\RR_+,\RR^d)$, $n\in\NN$, be such that $f_n\skor{d}f$ as $n\to\infty$.
By Ethier and Kurtz \cite[Proposition 5.3 (b), page 119]{EthKur}, there exists
  a sequence $(\lambda_n)_{n\in\NN}\subset\Lambda$ such that $\lim_{n\to\infty}\gamma(\lambda_n)=0$ and
\begin{equation*}
\lim_{n\to\infty}\sup_{t\in[0,T]}\Vert f_n(t)-f(\lambda_n(t))\Vert=0,\qquad T\in\RR_{++}.
\end{equation*}
Consequently, by \eqref{int_hely} and the triangle inequality, for all $T\in\RR_{++}$, we have
\begin{align*}
\sup_{t\in[0,T]}\left\Vert(\varphi(f_n))(t)-(\varphi(f))(\lambda_n(t))\right\Vert
	=&\sup_{t\in[0,T]}\left\Vert\int_0^tf_n(s)\,\dd s-\int_0^{\lambda_n(t)}f(s)\,\dd s\right\Vert\\
	=&\sup_{t\in[0,T]}\left\Vert\int_0^tf_n(s)\,\dd s-\int_0^tf(\lambda_n(s))\lambda_n'(s)\,\dd s\right\Vert\\
	\leq&\sup_{t\in[0,T]}\left\Vert\int_0^tf_n(s)\,\dd s-\int_0^tf(\lambda_n(s))\,\dd s\right\Vert\\
	&+\sup_{t\in[0,T]}\int_0^t\left\Vert f(\lambda_n(s))\right\Vert\esssup_{ r\in\RR_+}\left\Vert1-\lambda_n'(r)\right\Vert\,\dd s\\
	\leq&\,T\sup_{s\in[0,T]}\left\Vert f_n(s)-f(\lambda_n(s))\right\Vert
	+ T \gamma(\lambda_n)\sup_{s\in[0,T]}\left\Vert f(\lambda_n(s))\right\Vert.
\end{align*}
By Ethier and Kurtz \cite[formulae (5.6) and (5.7), page 117]{EthKur} and $\lim_{n\to\infty}\gamma(\lambda_n)=0$,
 we get
 \[
   \lim_{n\to\infty} \sup_{t\in[0,T]} \vert \lambda_n(t) - t \vert =0, \qquad T\in\RR_{++}.
 \]
Consequently, for all $T\in\RR_{++}$, we have that $\lambda_n(T)\in(T-1,T+1)$ for sufficiently large values of $n\in\NN$,
 thus for such $n\in\NN$, using the monotonity of $\lambda_n$, we have $\lambda_n([0,T])\subset[0,T+1)$.
Therefore, for all $T\in\RR_{++}$ we get that
 \begin{equation*}
 \sup_{t\in[0,T]}\left\Vert(\varphi(f_n))(t)-(\varphi(f))(\lambda_n(t))\right\Vert
	\leq T \sup_{s\in[0,T]}\left\Vert f_n(s)-f(\lambda_n(s))\right\Vert
	+ T \gamma(\lambda_n)\sup_{s\in[0,T+1]}\left\Vert f(s)\right\Vert,
 \end{equation*}
 where $\sup_{s\in[0,T+1]}\left\Vert f(s)\right\Vert<\infty$ (due to the fact that a c\`{a}dl\`{a}g function is bounded,
  see, e.g., Billingsley \cite[Chapter 3, (12.5)]{Bil}),
and $\sup_{s\in[0,T]}\left\Vert f_n(s)-f(\lambda_n(s))\right\Vert$ and $\gamma(\lambda_n)$ both converge to $0$ as $n\to\infty$.
Using again Ethier and Kurtz \cite[Proposition 5.3, page 119]{EthKur}, we obtain that
  $\varphi(f_n)\skor{d}\varphi(f)$ as $n\to\infty$, thus $\varphi$ is continuous.
\proofend

\begin{Lem}\label{ConvIteratedIntegral}
Let $d\in\NN$, $\Phi:\DD(\RR_+,\RR^d)\to \DD(\RR_+,\RR^{3d})$ and
$\Phi_n:\DD(\RR_+,\RR^d)\to \DD(\RR_+,\RR^{3d})$, $n\in\NN$, be mappings defined by
\begin{equation*}
(\Phi(f))(t):=\begin{bmatrix}
                     f(t) \\[1mm]
                    \int_0^t f(s)\,\dd s\\
		\int_0^t\left(\int_0^rf(s)\,\dd s\right)\,\dd r
                  \end{bmatrix},
                  \qquad
 (\Phi_n(f))(t):=\begin{bmatrix}
                    f(t) \\[1mm]
                    \int_0^{\frac{\nt}{n}}  f(s)\,\dd s\\
		\int_0^{\frac{\nt}{n}}\left(\int_0^{\frac{\nr}{n}}f(s)\,\dd s\right)\,\dd r
                  \end{bmatrix}
\end{equation*}
 for $n\in\NN$, $t\in\RR_+$, and $f\in \DD(\RR_+,\RR^d)$.
If $(\bcU_t)_{t \in \RR_+}$ and $(\bcU^{(n)}_t)_{t \in \RR_+}$, $n \in \NN$,
 are $\RR^d$-valued stochastic processes with c\`adl\`ag paths such that $\PP(\bcU\in\CC(\RR_+,\RR^d))=1$ and
 $\bcU^{(n)} \distr \bcU$ as $n\to\infty$, then $\Phi_n(\bcU^{(n)}) \distr \Phi(\bcU)$ as $n\to\infty$.
\end{Lem}

\noindent\textbf{Proof.} 
We will apply Lemma \ref{Conv2Funct}, thus, we first check that the conditions of this lemma hold.
Let $\varphi:\DD(\RR_+,\RR^d)\to\DD(\RR_+,\RR^d)$
 and $\varphi_n:\DD(\RR_+,\RR^d)\to\DD(\RR_+,\RR^d)$, $n\in\NN$, be mappings defined by
\begin{equation*}
 (\varphi(f))(t):=\int_0^tf(s)\,\dd s,\qquad (\varphi_n(f))(t):=\int_0^{\frac{\nt}{n}}f(s)\,\dd s,\qquad f\in\DD(\RR_+,\RR^d),\;  t\in\RR_+,\;n\in\NN.
\end{equation*}
Then we have
\begin{equation*}
(\Phi(f))(t)=\begin{bmatrix}
f(t)\\
(\varphi(f))(t)\\
((\varphi\circ\varphi)(f))(t)
\end{bmatrix},\qquad
(\Phi_n(f))(t)=\begin{bmatrix}
f(t)\\
(\varphi_n(f))(t)\\
((\varphi_n\circ\varphi_n))(f)(t)
\end{bmatrix}
\end{equation*}
for $f\in\DD(\RR_+,\RR^d)$, $t\in\RR_+$ and $n\in\NN$.

We check that $\Phi$ is continuous (in particular, it is measurable).
Let $f\in\DD(\RR_+,\RR^d)$, $f_n\in\DD(\RR_+,\RR^d)$, $n\in\NN$, be
such that $f_n\skor{d} f$ as $n\to\infty$.
By Lemma \ref{IntegralCont}, the mappings $\varphi$ and $\varphi\circ\varphi$
are continuous and map into $\CC(\RR_+,\RR^d)$.
 These imply that $\varphi(f_n)\skor{d}\varphi(f)$ as $n\to\infty$
and $(\varphi\circ\varphi)(f_n)\skor{d}(\varphi\circ\varphi)(f)$ as $n\to\infty$.
By applying Lemma \ref{CoordConv} twice, we get that $\Phi(f_n)\skor{3d}\Phi(f)$ as $n\to\infty$,
 yielding that $\Phi$ is continuous.

Now we turn to prove that $\Phi_n$, $n\in\NN$, are measurable.
 Let $n\in\NN$ be fixed in this paragraph.
By Lemma \ref{coordMeasurable}, it is sufficient
 to check the measurability of $\varphi_n$ and $\varphi_n\circ\varphi_n$,
 since the identity map is trivially measurable.
Using that the finite dimensional sets in $\DD(\RR_+,\RR^d)$ generate the Borel $\sigma$-algebra on $\DD(\RR_+,\RR^d)$
(see, e.g., Jacod and Shiryaev \cite[Chapter VI, Theorem 1.14, part c)]{JacShi}), to check the Borel measurability of $\varphi_n$
 it is enough to verify that the mappings $\pi_t\circ \varphi_n:\DD(\RR_+,\RR^d)\to\RR^d$
 are measurable for  all $t\in\RR_+$,
 where $\pi_t: \DD(\RR_+,\RR^d)\to \RR^d$, $\pi_t(g):=g(t)$, $g\in \DD(\RR_+,\RR^d)$, is the natural projection onto $t$.
Notice that $\pi_t\circ \varphi_n=\pi_{\nt/n}\circ\varphi$ for all $t\in\RR_+$.
Since $\varphi$ is continuous (in particular, measurable)
 and the natural projection onto $\frac{\nt}{n}$ is measurable for all $t\in\RR_+$
 (see, e.g., Billingsley \cite[Theorem 16.6, part (i)]{Bil}),
we get that $\pi_t \circ \varphi_n$ is measurable as well, and thus so is $\varphi_n$.
The measurability of $\varphi_n\circ\varphi_n$ easily follows, since it is a composition of (Borel) measurable maps.
Consequently, as we already explained, using Lemma \ref{coordMeasurable}, we get that $\Phi_n$ is measurable.

To show that $C_{\Phi,(\Phi_n)_{n\in\NN}}=\CC(\RR_+,\RR^d)$, we verify that
$\Phi_n(f_n)\lu\Phi(f)$ as $n\to\infty$ whenever $f_n\lu f$ as $n\to\infty$ with $f\in\CC(\RR_+,\RR^d)$ and $f_n\in\DD(\RR_+,\RR^d)$, $n\in\NN$.
Using that $\Vert[\bz_1,\bz_2,\bz_3]^\top\Vert\leq\Vert\bz_1\Vert+\Vert\bz_2\Vert+\Vert\bz_3\Vert$ for all $\bz_1,\bz_2,\bz_3\in\RR^d$,
we have for all $n\in\NN$ and $t\in\RR_+$,
\begin{align*}
 \Vert (\Phi_n(f_n))(t)-(\Phi(f))(t)\Vert
    &\leq\Vert f_n(t)-f(t)\Vert
         +\Vert (\varphi_n(f_n))(t)-(\varphi(f))(t)\Vert \\
    &\phantom{\leq}+\Vert((\varphi_n\circ\varphi_n)(f_n))(t)- ((\varphi\circ\varphi)(f))(t)\Vert.
\end{align*}
Hence, for all $t\in\RR_+$, we have
\begin{align*}
\Vert (\Phi_n(f_n))(t)- (\Phi(f))(t) \Vert
	\leq&\left\Vert f_n(t)-f(t)\right\Vert +\left\Vert \int_0^{\frac{\nt}{n}}f_n(s)\,\dd s-\int_0^tf(s)\,\dd s\right\Vert \\
	&+ \left\Vert \int_0^{\frac{\nt}{n}}\left(\int_0^{\frac{\nr}{n}}f_n(s)\,\dd s\right)\,\dd r-\int_0^t\left(\int_0^rf(s)\,\dd s\right)\,\dd r\right\Vert \\
	\leq&\left\Vert f_n(t)-f(t)\right\Vert +\int_0^{\frac{\nt}{n}}\left\Vert f_n(s)-f(s)\right\Vert \,\dd s+\int_{\frac{\nt}{n}}^t\left\Vert f(s)\right\Vert \,\dd s\\
	&+\int_0^{\frac{\nt}{n}}\left(\int_0^{\frac{\nr}{n}}\left\Vert f_n(s)-f(s)\right\Vert \,\dd s\right)\,\dd r
	+\int_{\frac{\nt}{n}}^t\left(\int_0^r\left\Vert f(s)\right\Vert \,\dd s\right)\,\dd r\\
	&+\int_0^{\frac{\nt}{n}}\left(\int_{\frac{\nr}{n}}^r\left\Vert f(s)\right\Vert \,\dd s\right)\,\dd r.
\end{align*}
Using also that
 \[
   \sup_{t\in[0,T]} \left( \frac{nt-\nt}{n} \right)
      \leq \frac{1}{n}, \qquad T\in\RR_{++}, \; n\in\NN,
 \]
 for all $T\in\RR_{++}$ and $n\in\NN$, we get that
\begin{align*}
   \sup_{t\in[0,T]}\Vert (\Phi_n(f_n))(t)-&(\Phi(f))(t)\Vert
	\leq \sup_{t\in[0,T]}\left\Vert f_n(t)-f(t)\right\Vert +\frac{\nT}{n}\sup_{t\in[0,T]}\left\Vert f_n(t)-f(t)\right\Vert\\
            &+  \sup_{t\in[0,T]} \left( \frac{nt-\nt}{n} \right)  \sup_{t\in[0,T]}\left\Vert f(t)\right\Vert
      +\frac{\nT^2}{n^2}\sup_{t\in[0,T]}\left\Vert f_n(t)-f(t)\right\Vert\\
       &+  T  \sup_{t\in[0,T]} \left( \frac{nt-\nt}{n} \right) \sup_{t\in[0,T]}\left\Vert f(t)\right\Vert 
+\frac{\nT}{n}   \sup_{t\in[0,T]} \left( \frac{nt-\nt}{n} \right) \sup_{t\in[0,T]}\left\Vert f(t)\right\Vert \\
	\leq&\left(1+T+T^2\right)\sup_{t\in[0,T]}\left\Vert f_n(t)-f(t)\right\Vert +n^{-1}\left(1+2T\right)\sup_{t\in[0,T]}\left\Vert f(t)\right\Vert .
\end{align*}
For all $T\in\RR_{++}$, we have $\sup_{t\in[0,T]}\left\Vert f(t)\right\Vert<\infty$, since $f\in\CC(\RR_+,\RR^d)$,
thus $n^{-1}\sup_{t\in[0,T]}\left\Vert f(t)\right\Vert\to0$ as $n\to\infty$.
Since $f_n\lu f$ as $n\to\infty$, we have $\sup_{t\in[0,T]}\left\Vert f_n(t)-f(t)\right\Vert\to0$ as $n\to\infty$ for all $T\in\RR_+$.
Thus we get that $\Phi_n(f_n)\lu\Phi(f)$ as $n\to\infty$,
and so $C_{\Phi,(\Phi_n)_{n\in\NN}}=\CC(\RR_+,\RR^d)$.
Consequently, one can apply Lemma \ref{Conv2Funct} by choosing $C:=\CC(\RR_+,\RR^d)$,
 which yields that $\Phi_n(\bcU^{(n)}) \distr \Phi(\bcU)$ as $n\to\infty$, as desired.
\proofend

Finally, we recall a result on weak convergence for the sum of stochastic processes with c\`adl\`ag paths 
 due to Jacod and Shiryaev \cite[Lemma VI.3.31]{JacShi}. 

\begin{Lem}[Jacod and Shiryaev {\cite[Lemma VI.3.31]{JacShi}}]\label{Lem_JacShi}
Let $(\bcY^{(n)}_t)_{t\in\RR_+}$, $n\in\NN$,  $(\bcY_t)_{t\in\RR_+}$, and  $(\bcZ^{(n)}_t)_{t\in\RR_+}$, $n\in\NN$,
 be $\RR^d$-valued stochastic processes with c\`adl\`ag paths on a probability space $(\Omega,\cA,\PP)$.
Suppose that $\bcY^{(n)}\distr \bcY$ as $n\to\infty$ and 
 \[
   \lim_{n\to\infty} \PP\Big( \sup_{t\in[0,T]} \Vert \bcZ^{(n)}_t\Vert > \vare \Big) = 0
     \qquad \text{for all \ $T\in\RR_{++}$ \ and \ $\vare\in\RR_{++}$.}
 \] 
Then $\bcY^{(n)} + \bcZ^{(n)} \distr \bcY$ as $n\to\infty$. 
\end{Lem}

\section{Convergence of random step processes}
\label{section_conv_step_processes}

We recall a result about convergence of random step processes towards a diffusion process, see Isp\'any and Pap \cite{IspPap}.

\begin{Thm}\label{Conv2DiffThm}
Let $\bbeta : \RR_+ \times \RR^d \to \RR^d$ and $\bgamma : \RR_+ \times \RR^d \to \RR^{d \times r}$
 be continuous functions.
Assume that uniqueness in the sense of probability law holds for the SDE
 \begin{equation}\label{SDE}
  \dd \, \bcU_t
  = \bbeta (t, \bcU_t) \, \dd t + \bgamma (t, \bcU_t) \, \dd \bcW_t ,
  \qquad t \in \RR_+,
 \end{equation}
 with initial value $\bcU_0 = \bu_0$ for all $\bu_0 \in \RR^d$, where
 $(\bcW_t)_{t\in\RR_+}$ is an $r$-dimensional standard Wiener process.
Let $(\bcU_t)_{t\in\RR_+}$ be a solution of \eqref{SDE} with initial value
 $\bcU_0 = \bzero \in \RR^d$.

For each $n \in \NN$, let $(\bU^{(n)}_k)_{k\in\ZZ_+}$ be a sequence of $d$-dimensional random vectors adapted to a filtration
  $(\cF^{(n)}_k)_{k\in\ZZ_+}$ (i.e., $\bU^{(n)}_k$ is $\cF^{(n)}_k$-measurable) such that
  $\EE(\Vert \bU^{(n)}_k\Vert^2)<\infty$ for each $n,k \in \NN$.
Let
 \[
   \bcU^{(n)}_t := \sum_{k=0}^{\nt}  \bU^{(n)}_k \, ,
   \qquad t \in \RR_+, \quad n \in \NN .
 \]
Suppose that $\bcU^{(n)}_0 = \bU^{(n)}_0 \distr \bzero$ as $n\to\infty$ and that for all $T \in \RR_{++}$,
 \begin{enumerate}
  \item[\textup{(i)}]
   $\sup\limits_{t\in[0,T]}
     \biggl\|\sum\limits_{k=1}^{\nt}
              \EE\bigl(\bU^{(n)}_k \mid \cF^{(n)}_{k-1}\bigr)
             - \int_0^t \bbeta(s,\bcU^{(n)}_s) \dd s\biggr\|
   \stoch 0$ \ as \ $n\to\infty$,  \\
  \item[\textup{(ii)}]
   $\sup\limits_{t\in[0,T]}
     \biggl\|\sum\limits_{k=1}^{\nt}
              \var\bigl(\bU^{(n)}_k \mid \cF^{(n)}_{k-1}\bigr)
             - \int_0^t
                \bgamma(s,\bcU^{(n)}_s) \bgamma(s,\bcU^{(n)}_s)^\top
                \dd s\biggr\|
   \stoch 0$ \ as \ $n\to\infty$, \\
  \item[\textup{(iii)}]
   $\sum\limits_{k=1}^{\lfloor nT \rfloor}
     \EE\bigl(\|\bU^{(n)}_k\|^2 \bbone_{\{\|\bU^{(n)}_k\| > \theta\}}
             \bmid \cF^{(n)}_{k-1}\bigr)
   \stoch 0$ \ as \ $n\to\infty$ \ for all \ $\theta \in \RR_{++}$.
 \end{enumerate}
Then $\bcU^{(n)} \distr \bcU$ as $n \to \infty$.
\end{Thm}

Note that in (ii) of Theorem \ref{Conv2DiffThm}, $\|\cdot\|$ denotes
 an operator norm, while in (i) it denotes a vector norm.

The following result is about the asymptotic behavior of a single-type Galton-Watson process with immigration in the critical case,
 it is due to Wei and Winnicki \cite[Theorem 2.1]{WW}.

\begin{Thm}\label{thm:critical}
Let $(X_k)_{k\in\ZZ_+}$ be a single-type Galton-Watson process with immigration such that
  $\EE(\xi^2) < \infty$, $\EE(\vare^2) < \infty$, $\EE(\xi) = 1$ (critical case) and $\EE(X_0^2) < \infty$.
Then
 \begin{equation*}%\label{convX}
  (n^{-1} X_\nt)_{t\in\RR_+} \distr (\cX_t)_{t\in\RR_+} \qquad
  \text{as \ $n \to \infty$,}
 \end{equation*}
 where the limit process $(\cX_t)_{t\in\RR_+}$ is the pathwise unique strong solution of the SDE
 \begin{equation*}%\label{SDE_single}
  \dd \cX_t
  = \EE(\vare) \, \dd t
    + \sqrt{\var(\xi) \, \cX_t^+} \, \dd \cW_t , \qquad t\in\RR_+,
  \qquad \cX_0 = 0 ,
 \end{equation*}
 where $(\cW_t)_{t\in\RR_+}$ is a standard Wiener process.
\end{Thm}

\section*{Acknowledgements}
This paper was initiated by our longtime co-author, mentor and dear friend Gyula Pap who passed away in October 2019.
We would like to thank P\'eter Kevei for the heuristic argument that we present in Remark \ref{Heuristic_Remark}.
We acknowledge the valuable suggestions both from the referees and the associate editor.


\addcontentsline{toc}{section}{References}

\begin{thebibliography}{99}

\bibitem{AN}
\textsc{Athreya, K. B.} and \textsc{Ney, P. E.} (1972).
\textit{Branching Processes},
Springer-Verlag, New York, Heidelberg.

\bibitem{BalMilWai}
\textsc{Balelli, I., Mili\v{s}i\`c, V.} and \textsc{Wainrib, G.} (2019).
Multi-type {Galton}-{Watson} processes with affinity-dependent selection applied to antibody affinity maturation.
\textit{Bulletin of Mathematical Biology}
\textbf{81(3)} 830--868.


\bibitem{BarBezPap1}
\textsc{Barczy, M., Bezdány, D.} and \textsc{Pap, G.} (2021).
A note on asymptotic behavior of critical Galton-Watson processes with immigration.
\textit{Involve: A Journal of Mathematics}
\textbf{14(5)} 871--891.

\bibitem{BarBezPap2}
\textsc{Barczy, M., Bezdány, D.} and \textsc{Pap, G.} (2023).
Asymptotic behaviour of critical decomposable 2-type Galton-Watson processes with immigration.
\textit{Stochastic Processes and their Applications}
\textbf{160} 318--350.

\bibitem{BarIspPap0}
\textsc{Barczy, M., Isp\'any, M.} and \textsc{Pap, G.} (2011).
Asymptotic behavior of unstable \INARp\ processes.
\textit{Stochastic Processes and their Applications}
\textbf{121(3)} 583--608.

\bibitem{BarNedPap}
\textsc{Barczy, M., Ned\'enyi, F. K.} and \textsc{Pap, G.} (2018).
On aggregation of multitype Galton-Watson branching processes with immigration.
\textit{Modern Stochastics: Theory and Applications}
\textbf{5(1)} 53--79.

\bibitem{Bil} {\sc Billingsley, P.} (1999)
  {\sl Convergence of Probability Measures, 2nd ed.}
 Wiley-Interscience Publication, New York.

\bibitem{ColGol}
\textsc{Coldman, A. J.} and \textsc{Goldie, J. H.} (1986).
A stochastic model for the origin and treatment of tumors containing drug-resistant cells.
\textit{Bulletin of Mathematical Biology}
\textbf{48} 279--292.

\bibitem{EthKur}
\textsc{Ethier, S. N.} and \textsc{Kurtz, T. G.} (1986).
\textit{Markov Processes. Characterization and Convergence}.
Wiley, New York.

\bibitem{FosNey}
\textsc{Foster, J.} and \textsc{Ney, P.} (1978).
Limit laws for decomposable critical branching processes.
\textit{Zeitschrift für Wahrscheinlichkeitstheorie und verwandte Gebiete}
\textbf{46} 13--43.


\bibitem{HJ}
\textsc{Horn, R. A.} and \textsc{Johnson, Ch.\ R.} (1985).
\textit{Matrix Analysis}.
Cambridge University Press, Cambridge.

\bibitem{IkeWat}
\textsc{Ikeda, N.} and \textsc{Watanabe, S.} (1989).
\textit{Stochastic Differential Equations and Diffusion Processes, 2nd ed.}
North-Holland, Kodansha, Amsterdam, Tokyo.

\bibitem{IspPap}
\textsc{Isp\'any, M.} and \textsc{Pap, G.} (2010).
A note on weak convergence of random step processes.
\textit{Acta Mathematica Hungarica}
\textbf{126(4)} 381--395.

\bibitem{IspPap2}
\textsc{Isp\'any, M.} and \textsc{Pap, G.} (2014).
Asymptotic behavior of critical primitive multi-type branching processes with
 immigration.
\textit{Stochastic Analysis and Applications}
\textbf{32(5)} 727--741.

\bibitem{JacShi}
\textsc{Jacod, J.} and \textsc{Shiryaev, A. N.} (2003).
\textit{Limit Theorems for Stochastic Processes, 2nd ed.}
Springer-Verlag, Berlin.

\bibitem{Jag}
\textsc{Jagers, P.} (1969).
The proportions of individuals of different kinds in two-type populations.
A branching process problem arising in biology.
\textit{Journal of Applied Probability}
\textbf{6(2)} 249--260.

\bibitem{Kal}
\textsc{Kallenberg, O.} (1997).
\textit{Foundations of Modern Probability}.
Springer, New York, Berlin, Heidelberg.

\bibitem{Kle}
\textsc{Klenke, A.} (2020).
\textit{Probability Theory. A Comprehensive Course, 3rd ed.}
Springer, Cham.

\bibitem{Latour}
\textsc{Latour, A.} (1998).
Existence and stochastic structure of a non-negative integer-valued autoregressive process.
\textit{Journal of Time Series Analysis}
\textbf{19(4)} 439--455.

\bibitem{LeeYan}
\textsc{Lee, C. S. } and \textsc{Yang, G. L.}  (1995).
A multitype decomposable age-dependent branching process and its applications.
\textit{Journal of Applied Probability}
\textbf{32(3)} 591--608.

\bibitem{Q}
\textsc{Quine, M. P.} (1970).
The multi-type Galton-Watson process with immigration.
\textit{Journal of Applied Probability}
\textbf{7(2)} 411--422.

\bibitem{RevYor}
\textsc{Revuz, D.} and \textsc{Yor, M.} (2001).
\textit{Continuous Martingales and Brownian Motion},
 3rd ed., corrected 2nd printing.
Springer-Verlag, Berlin.

\bibitem{SagStaahl}
\textsc{Sagitov, S.} and \textsc{St{\aa}hlberg, A.} (2023).
Counting unique molecular identifiers in sequencing using a multi-type branching process with immigration.
\textit{Journal of Theoretical Biology}
\textbf{558}, 111365.

\bibitem{Wan}
\textsc{Wang, X.} (2018).
\textit{Lecture Notes in Real Analysis}.
 Birkh\"auser, Cham.

\bibitem{WW}
\textsc{Wei, C. Z.} and \textsc{Winnicki, J.} (1989).
Some asymptotic results for the branching process with immigration.
\textit{Stochastic Processes and their Applications}
\textbf{31(2)} 261--282.

\end{thebibliography}
\end{document}